\documentclass[reqno,fleqn]{amsart}
\usepackage{amsfonts,amsthm,amssymb}
\usepackage{lipsum}
\usepackage[utf8]{inputenc}
\usepackage[english]{babel}

\usepackage{lmodern,textcomp}
\usepackage{amsmath}
\usepackage{mathtools}
\usepackage{latexsym}
\usepackage{tikz}
\usepackage{fancyvrb}
\usepackage{epsfig}
\usepackage{pstricks,slashbox,multirow}
\usepackage{graphicx}
\usepackage{rotating}
\usetikzlibrary{positioning}
\usepackage{subcaption}
\usepackage{datetime}
\newtheorem{theorem}{Theorem}

\newtheorem{definition}[theorem]{Definition}

\newtheorem{prm}{Problem}
\newtheorem{oprm}{Open Problem}

\newtheorem{rem}[theorem]{Remark}

\title[A study on Type-2 isomorphic $C_n(R)$: Part 4: 960 triples of Type-2 isomorphic $C_{54}(R)$]{A study on Type-2 isomorphic circulant graphs. \\ Part 4: 960 triples of Type-2 isomorphic circulant graphs $C_{54}(R)$}

\author{\sc Vilfred Kamalappan} 
\address{Department of Mathematics, Central University of Kerala, Periye, Kasaragod, Kerala, India - 671 316.}
\email{vilfredkamal@gmail.com}



\date{}

\begin{document}

\begin{abstract} This article is the $4^{th}$ part of a detailed study on Type-2 isomorphic circulant graphs having ten parts \cite{v2-1}-\cite{v2-10} and is a continuation of Part 3. Here, we obtain all 960 triples of Type-2 isomorphic circulant graphs of order 54 and each triple of isomorphic circulant graphs is of Type-2 isomorphic w.r.t. $m$ = 3. 
\end{abstract}

\maketitle

	
\section{Introduction}

In \cite{v2-2} and \cite{v2-3}, we studied Type-2 isomorphic circulant graphs of orders 16, 24 and 32 and shown that the number of pairs of Type-2 isomorphic circulant graphs of orders 16, 24 and 32 are 8, 32 and 384, respectively.  This paper is a continuation of \cite{v2-2} and \cite{v2-3}. Using modified definition \ref{d3}, we study Type-2 isomorphic circulant graphs of order 54 and obtain all the 960 triples of isomorphic circulant graphs of order 54 and of Type-2 w.r.t. $m$ = 3. For basic definitions and results on isomorphic circulant graphs, refer \cite{v2-1, v2-3}.

\begin{definition}{\rm\cite{ad67}} \quad \label{d1} For $R =$ $\{r_1$, $r_2$, $\ldots$, $r_k\}$ and $S$ = $\{s_1$, $s_2$, $\ldots$, $s_k\}$, circulant graphs $C_n(R)$ and $C_n(S)$ are {\it Adam's isomorphic} if there exists a positive integer $x$ $\ni$ $\gcd(n, x)$ = 1 and $S$ = $\{xr_1$, $xr_2$, $\ldots$, $xr_k\}_n^*$ where $<r_i>_n^*$, the {\it reflexive modular reduction} of a sequence $< r_i >$, is the sequence obtained by reducing each $r_i$ under modulo $n$ to yield $r_i'$ and then replacing all resulting terms $r_i'$ which are larger than $\frac{n}{2}$ by $n-r_i'.$  
\end{definition}

A circulant graph $C_n(R)$ is said to have {\em Cayley Isomorphism (CI) property} if whenever $C_n(S)$ is isomorphic to $C_n(R)$, they are Adam’s isomorphic \cite{ad67}.

\begin{theorem} \cite{v24} \label{t2} Let $Ad_n(C_n(R))$ = $\{\varphi_{n,x}(C_n(R)) = C_n(xR): x\in\varphi_n \}$. Then, $C_n(S)\in Ad_n(C_n(R))$ if and only if $Ad_n(C_n(R))$ = $Ad_n(C_n(S))$ if and only if $C_n(R)\in Ad_n(C_n(S))$. \hfill $\Box$
\end{theorem}

In \cite{v2-1}, Vilfred modified the definition of Type-2 isomorphism of $C_n(R)$ w.r.t. $m$  given in \cite{v2-2-arX} as follows and hereafter we use the definition \ref{d3}. 

\begin{definition} \cite{v2-1} \quad  \label{d3} Let $V(K_n) = \{u_0,u_1,u_2,...,u_{n-1}\}$, $V(C_n(R))$ = $\{v_0, v_1, v_2, ... , v_{n-1}\}$, $|R| \geq 3$, $r\in R$ and $m > 1$ and $m^3$ be divisors of $\gcd(n, r)$ and $n$, respectively. Define 1-1 mapping $\theta_{n,m,t} :$ $V(C_n(R)) \rightarrow V(K_n)$ such that $\theta_{n,m,t}(v_x) = u_{x+jtm}$,  $\theta_{n,m,t}((v_x, v_{x+s}))$ = $(\theta_{n,m,t}(v_x),$ $\theta_{n,m,t}(v_{x+s}))$ under subscript arithmetic modulo $n$ and $\theta_{n,m,t}(C_n(R))$ = $C_n(\theta_{n,m,t}(R))$ where $\theta_{n,m,t}(R)$ in $C_n(\theta_{n,m,t}(R))$ is calculated under the reflexive modulo $n$, $\forall$ $x \in \mathbb{Z}_n$, $x = qm+j,$ $0 \leq j \leq m-1$, $s\in R$ and $0 \leq q,t \leq \frac{n}{m} -1$. And for a particular value of $t,$ if  $\theta_{n,m,t}(C_n(R))$ = $C_n(S)$ for some $S$  and  $S \neq yR$ for all $y\in \varphi_n$ under reflexive modulo $n,$ then $C_n(R)$ and $C_n(S)$ are called {\em isomorphic circulant graphs of Type-2 w.r.t. $m$.} 
	
When $C_n(R)$ and $C_n(S)$ are Type-2 isomorphic w.r.t. $m$, then we also say that $C_{kn}(kR)$  and $C_{kn}(kS)$ are Type-2 isomorphic w.r.t. $m$, $k\in\mathbb{N}$. Here, $k.C_n(T)$ = $C_{kn}(kT)$, $k\in\mathbb{N}$. 	 
\end{definition}

\begin{rem} \cite{v2-1}  \label{r4} Following steps are used to establish Type-2 isomorphism w.r.t. $m$ between circulant graphs $C_n(R)$ and $C_n(S)$. (i) $R$ $\neq$ $S$ and $|R| = |S| \geq 3$; (ii) $\exists$ $r\in R,S$ and $m > 1$ $\ni$ $m$ is a divisor of $\gcd(n, r)$, $m^3$ is a divisor of $n$ and for some $t$ $\ni$ $1 \leq t \leq \frac{n}{m} -1$, $\theta_{n,m,t}(C_n(R))$ = $C_n(S)$ and (iii) $S$ $\neq$ $xR$ for all $x\in\varphi_n$ under arithmetic reflexive modulo $n$. 

Thus, if $C_n(R)$ and $C_n(S)$ are Type-2 isomorphic circulant graphs w.r.t. $m$, then there exist $r\in R,S$, $m > 1$ and some $t$ $\ni$ $m$ is a divisor of $\gcd(n, r)$, $m^3$ is a divisor of $n$, $1 \leq t \leq \frac{n}{m} -1$, $\theta_{n,m,t}(C_n(R))$ = $C_n(S)$ and $S$ $\neq$ $xR$ for all $x\in\varphi_n$ under arithmetic reflexive modulo $n$.
\end{rem} 

\begin{rem}  \label{r5} \quad The calculation on jump sizes $r_i$s which are integer multiples of $m$ need not be done under the transformation $\theta_{n,m,t}$, while searching for possible value(s) of $t$ for which the transformed graph $\theta_{n,m,t}(C_n(R))$ is circulant of the form $C_n(S)$ for some $S \subseteq [1, \frac{n}{2}]$, as there is no change in these $r_i$s where $r\in R$ and $m > 1$ and $m^3$ are divisors of $\gcd(n, r)$ and $n$, respectively. 

Thus, if $\theta_{n,m,t}(C_n(R))$ = $C_n(S)$ for some $S$ and thereby $C_n(R)$ $\cong$ $C_n(S)$, then $\theta_{n,m,t}(C_n(R \cup mT))$ = $C_n(S \cup mT)$ for any $T$ and thereby $C_n(R \cup mT)$ $\cong$ $C_n(S \cup mT)$.

Also, for a given $C_n(R)$, w.r.t. different values of $m$, we may get different Type-2 isomorphic circulant graphs.
\end{rem}

\begin{rem} \cite{v2-1} \label{r6} \quad {\rm For given $C_n(R)$ and $C_n(S)$ when either $\theta_{n,m,t}(C_n(R))$ = $C_n(S)$ for some $t$ or $C_n(xR)$ = $C_n(S)$ for some $x$, then $C_n(R)$ and $C_n(S)$ are isomorphic, $0 \leq t \leq \frac{n}{m} -1$ and $x\in\varphi_n$. }
\end{rem}

\begin{theorem}{\rm \cite{v24}}\quad \label{t7} {\rm For $n \geq 2$, $1 \leq 2s-1 \leq 2n-1$, $n \neq 2s-1$, $R$ = $\{2,2s-1, 4n-(2s-1)\}$ and $S$ = $\{ 2,$ $2n-(2s-1)$, $2n+2s-1 \}$, $\theta_{8n,2,n}(C_{8n}(R))$ = $C_{8n}(S)$ = $\theta_{8n,2,3n}(C_{8n}(R)),$ $\theta_{8n,2,n}(C_{8n}(S))$ = $C_{8n}(R)$ = $\theta_{8n,2,3n}(C_{8n}(S))$ and circulant graphs $C_{8n}(R)$ and $C_{8n}(S)$ are Type-2 isomorphic  w.r.t. $m$ = 2. When $n$ = $2s-1$, the two circulant graphs are the same. \hfill $\Box$}
\end{theorem}

\begin{theorem} \cite{v2-6} \label{t8} {\rm Let $p$ be an odd prime number, $1 \leq i \leq p$, $1 \leq x \leq p-1$, $y\in\mathbb{N}_0$, $0 \leq y \leq np-1$, $1 \leq x+yp \leq np^2-1$, $d^{np^3, x+yp}_i = (i-1)xpn+x+yp$,  $R^{np^3, x+yp}_i$ $=$ $\{p$, $d^{np^3, x+yp}_i$, $np^2-d^{np^3, x+yp}_i$, $np^2+d^{np^3, x+yp}_i$, $2np^2-d^{np^3, x+yp}_i$, $2np^2+$ $d^{np^3, x+yp}_i,$ $3np^2-d^{np^3, x+yp}_i$, $3np^2+d^{np^3, x+yp}_i$, . . . , $(p-1)np^2$ - $d^{np^3, x+yp}_i$, $(p-1)np^2+d^{np^3, x+yp}_i$, $np^3-d^{np^3, x+yp}_i$, $np^3-p\}$ and $i,j,n,x\in\mathbb{N}$. Then, for a given set of values of $n$, $p$, $x$ and $y$, $\theta_{np^3,p,jn} (C_{np^3}(R^{np^3, x+yp}_i))$ = $C_{np^3}(R^{np^3, x+yp}_{i+j})$ and the $p$ circulant graphs $C_{np^3}(R^{np^3, x+yp}_i)$ are isomorphic of Type-2 w.r.t.  $p$, $1 \leq i,j \leq p$ where $i+j$ in $R^{np^3, x+yp}_{i+j}$ is calculated under addition modulo $p$ and $C_{np^3}(R^{np^3, x+yp}_0)$ = $C_{np^3}(R^{np^3, x+yp}_p)$. \hfill $\Box$}
\end{theorem}

\section{Main result}

Given a circulant graph $C_n(R)$ having isomorphic circulant graphs of Type-2 w.r.t. $m$, remark \ref{r5} helps us to obtain more  isomorphic graphs which covers Type-2 w.r.t. $m$ as well as some Type-1 isomorphic graphs of $C_n(R)$. In this section, we obtain all the 960 triples of Type-2 isomorphic circulant graphs of order 54. We start with a problem retaled to Type-1 isomorphism of circulant graphs $C_{54}(r_1,r_2,r_3)$.

\begin{prm} \label{p1} {\rm For $s$ = 9,18,27 and $\gcd(54, s)$ = $s$, show that the following statements are true.   
\begin{enumerate}
 \item [\rm (a)] $C_{54}(1,s,17,19)$, $C_{54}(5,s,13,23)$, $C_{54}(s,7,11,25)$ are Type-1 isomorphic; 
 \item [\rm (b)] $C_{54}(2,s,16,20)$, $C_{54}(4,s,14,22)$, $C_{54}(s,8,10,26)$ are Type-1 isomorphic.  
\end{enumerate} }	
\end{prm}
\noindent
{\bf Solution.}\quad Given triples of circulant graphs are Type-1 isomorphic by the following. 
\begin{enumerate}
\item [\rm (a1)] $C_{54}(5,9,13,23)$ = $C_{54}(5(1,9,17,19))$ and $C_{54}(7,9,11,25)$ = $C_{54}(7(1,9,17,19))$. 

$\Rightarrow$ $C_{54}(1,9,17,19)$, $C_{54}(5,9,13,23)$ and $C_{54}(7,9,11,25)$  are Type-1 isomorphic.

\item [\rm (a2)] $C_{54}(5,13,18,23)$ = $C_{54}(5(1,17,18,19))$ and $C_{54}(7,11,18,25)$ = $C_{54}(7(1,17,18,19))$. 

$\Rightarrow$ $C_{54}(1,17,18,19)$, $C_{54}(5,13,18,23)$ and $C_{54}(7,11,18,25)$ are Type-1 isomorphic.

\item [\rm (a3)] $C_{54}(5,13,23,27)$ = $C_{54}(5(1,17,19,27))$ and $C_{54}(7,11,25,27)$ = $C_{54}(7(1,17,19,27))$. 

$\Rightarrow$ $C_{54}(1,17,19,27)$, $C_{54}(5,13,23,27)$ and $C_{54}(7,11,25,27)$ are Type-1 isomorphic.

\item [\rm (b1)] $C_{54}(4,9,14,22)$ = $C_{54}(7(2,9,16,20))$ and $C_{54}(8,9,10,26)$ = $C_{54}(5(2,9,16,20))$. 

$\Rightarrow$ $C_{54}(2,9,16,20)$, $C_{54}(4,9,14,22)$ and $C_{54}(8,9,10,26)$  are Type-1 isomorphic.

\item [\rm (b2)] $C_{54}(4,14,18,22)$ = $C_{54}(7(2,16,18,20))$ and $C_{54}(8,10,18,26)$ = $C_{54}(5(2,16,18,20))$. 

$\Rightarrow$ $C_{54}(2,16,18,20)$, $C_{54}(4,14,18,22)$ and $C_{54}(8,10,18,26)$ are Type-1 isomorphic.

\item [\rm (b3)] $C_{54}(4,14,22,27)$ = $C_{54}(7(2,16,20,27))$ and $C_{54}(8,10,26,27)$ = $C_{54}(5(2,16,20,27))$. 

$\Rightarrow$ $C_{54}(2,16,20,27)$, $C_{54}(4,14,22,27)$ and $C_{54}(8,10,26,27)$ are Type-1 isomorphic. 	  \hfill $\Box$
\end{enumerate} 

\begin{prm} \label{p2} {\rm Show that the following triples of circulant graphs are Type-2 isomorphic w.r.t. $m$ = 3.  
\begin{enumerate}
\item [\rm (a)] $C_{54}(1,3,17,19)$, $C_{54}(3,7,11,25)$, $C_{54}(3,5,13,23)$;  and
\item [\rm (b)] $C_{54}(2,3,16,20)$, $C_{54}(3,4,14,22)$, $C_{54}(3,8,10,26)$.   
	\end{enumerate} }	
\end{prm}
\noindent
{\bf Solution.}\quad Let $R_1$ = $\{1,3,17,19\}$, $S_1$ = $\{3,7,11,25\}$, $T_1$ = $\{3,5,13,23\}$, 

\hspace{2cm} $R_2$ = $\{2,3,16,20\}$, $S_2$ = $\{3,4,14,22\}$, $T_2$ = $\{3,8,10,26\}$. 

Here, $n$ = 54 = $2\times 3^3$ and $r$ = $s$ so that $s\in R_i,S_i,T_i$ and $m$ = 3 = $\gcd(54, 3)$, $i$ = 1,2. 

We have $\theta_{54,3,t}(3)$ = 3 and $\theta_{54,3,t}(54-3)$ = $54-3$, $0 \leq t \leq 17$. 
\begin{enumerate}
	\item [\rm (a)]  $\theta_{54,3,2}(C_{54}(1,3,17,19))$ = $\theta_{54,3,2}(C_{54}(1,3,17,19, 35,37,51,53))$ 
	
	\hspace{3cm} = $C_{54}(7,3,29,25, 47,43,51,11)$ = $C_{54}(3,7,11,25)$ and 
	\\
	$\theta_{54,3,4}(C_{54}(1,3,17,19))$ = $\theta_{54,3,4}(C_{54}(1,3,17,19, 35,37,51,53))$ 
	
	\hspace{3cm} = $C_{54}(13,3,41,31, 5,49,51,23)$ = $C_{54}(3,5,13,23)$.
		
	$\Rightarrow$ $C_{54}(1,3,17,19)$ $\cong$ $C_{54}(3,7,11,25)$ and $C_{54}(1,3,17,19)$ $\cong$ $C_{54}(3,5,13,23)$. 
\\	
$Ad_{54}(C_{54}(1,3,17,19))$ = $Ad_{54}(C_{54}(1,3,17,19, 35,37,51,53))$ 
	
	\hfill = $\{C_{54}(x(1,3,17,19, 35,37,51,53)): x = 1,5,7,11,13,17,19,23,25\}$
	
	= $\{C_{54}(1,3,17,19)$, $C_{54}(5,13,15,23)$, $C_{54}(7,11,21,25)\}$. 
	
	$\Rightarrow$ $C_{54}(3,7,11,25),C_{54}(3,5,13,23)\notin Ad_{54}(C_{54}(1,3,17,19))$.
	
	$\Rightarrow$ $C_{54}(1,3,17,19)$, $C_{54}(3,7,11,25)$ and  $C_{54}(3,5,13,23)$ are Type-2 isomorphic w.r.t. $m$ = 3 = $\gcd(54,3)$.
	
\item [\rm (b)]  $\theta_{54,3,2}(C_{54}(2,3,16,20))$ = $\theta_{54,3,2}(C_{54}(2,3,16,20, 34,38,51,52))$ 

\hspace{3.2cm} 	= $C_{54}(14,3,22,32, 40,50,51,4)$ = $C_{54}(3,4,14,22)$ and 
\\
$\theta_{54,3,4}(C_{54}(2,3,16,20))$ = $\theta_{54,3,4}(C_{54}(2,3,16,20, 34,38,51,52))$ 

\hspace{3.2cm} 	= $C_{54}(26,3,28,44, 46,8,51,10)$ = $C_{54}(3,8,10,26)$.

$\Rightarrow$ $C_{54}(2,3,16,20)$ $\cong$ $C_{54}(3,4,14,22)$ and $C_{54}(2,3,16,20)$ $\cong$ $C_{54}(3,8,10,26)$.

$\Rightarrow$ $C_{54}(2,3,16,20)$ $\cong$ $C_{54}(3,4,14,22)$ $\cong$ $C_{54}(3,8,10,26)$.
\\
 $Ad_{54}(C_{54}(2,3,16,20))$ = $Ad_{54}(C_{54}(2,3,16,20, 34,38,51,52))$ 

\hfill = $\{C_{54}(x(2,3,16,20, 34,38,51,52)): x = 1,5,7,11,13,17,19,23,25\}$. 

      = $\{C_{54}(2,3,16,20)$, $C_{54}(8,10,15,26)$, $C_{54}(4,14,21,22)\}$. 

$\Rightarrow$ $C_{54}(3,4,14,22),C_{54}(3,8,10,26)\notin Ad_{54}(C_{54}(2,3,16,20))$.

$\Rightarrow$ $(C_{54}(2,3,16,20)$, $C_{54}(3,4,14,22)$ $\&$  $C_{54}(3,8,10,26)$ are Type-2 isomorphic w.r.t. $m$ = 3 = $\gcd(54,3)$.  \hfill $\Box$
\end{enumerate} 

It is easy to obtain the following.

\begin{enumerate}
	\item [\rm (a1)]  $Ad_{54}(C_{54}(1,3,17,19))$ = $\{C_{54}(1,3,17,19)$, $C_{54}(5,13,15,23)$, $C_{54}(7,11,21,25)\}$

\hfill = $Ad_{54}(C_{54}(5,13,15,23))$ = $Ad_{54}(C_{54}(7,11,21,25))$. 
	
	\item [\rm (a2)]  $Ad_{54}(C_{54}(1,6,17,19))$ = $\{C_{54}(1,6,17,19)$, $C_{54}(5,13,23,24)$, $C_{54}(7,11,12,25)\}$

\hfill = $Ad_{54}(C_{54}(5,13,23,24))$ = $Ad_{54}(C_{54}(7,11,12,25))$. 
	
	\item [\rm (a3)]  $Ad_{54}(C_{54}(1,12,17,19))$ = $\{C_{54}(1,12,17,19)$, $C_{54}(5,6,13,23)$, $C_{54}(7,11,24,25)\}$

\hfill = $Ad_{54}(C_{54}(5,6,13,23))$ = $Ad_{54}(C_{54}(7,11,24,25))$.

	\item [\rm (a4)]  $Ad_{54}(C_{54}(1,15,17,19))$ = $\{C_{54}(1,15,17,19)$, $C_{54}(5,13,21,23)$, $C_{54}(3,7,11,25)\}$

\hfill = $Ad_{54}(C_{54}(5,13,21,23))$ = $Ad_{54}(C_{54}(3,7,11,25))$. 
	
	\item [\rm (a5)]  $Ad_{54}(C_{54}(1,17,19,21))$ = $\{C_{54}(1,17,19,21)$, $C_{54}(3,5,13,23)$, $C_{54}(7,11,15,25)\}$

\hfill = $Ad_{54}(C_{54}(3,5,13,23))$ = $Ad_{54}(C_{54}(7,11,15,25))$. 
	
	\item [\rm (a6)]  $Ad_{54}(C_{54}(1,17,19,24))$ = $\{C_{54}(1,17,19,24)$, $C_{54}(5,12,13,23)$, $C_{54}(6,7,11,25)\}$

\hfill = $Ad_{54}(C_{54}(5,12,13,23))$ = $Ad_{54}(C_{54}(6,7,11,25))$. 
	
\item [\rm (b1)]  $Ad_{54}(C_{54}(2,3,16,20))$ = $\{C_{54}(2,3,16,20)$, $C_{54}(8,10,15,26)$, $C_{54}(4,14,21,22)\}$

\hfill = $Ad_{54}(C_{54}(8,10,15,26))$ = $Ad_{54}(C_{54}(4,14,21,22))$. 

\item [\rm (b2)]  $Ad_{54}(C_{54}(2,6,16,20))$ = $\{C_{54}(2,6,16,20)$, $C_{54}(8,10,24,26)$, $C_{54}(4,12,14,22)\}$

\hfill = $Ad_{54}(C_{54}(8,10,24,26))$ = $Ad_{54}(C_{54}(4,12,14,22))$. 

\item [\rm (b3)]  $Ad_{54}(C_{54}(2,12,16,20))$ = $\{C_{54}(2,12,16,20)$, $C_{54}(6,8,10,26)$, $C_{54}(4,14,22,24)\}$

\hfill = $Ad_{54}(C_{54}(6,8,10,26))$ = $Ad_{54}(C_{54}(4,14,22,24))$. 

\item [\rm (b4)]  $Ad_{54}(C_{54}(2,15,16,20))$ = $\{C_{54}(2,15,16,20)$, $C_{54}(8,10,21,26)$, $C_{54}(3,4,14,22)\}$

\hfill = $Ad_{54}(C_{54}(8,10,21,26))$ = $Ad_{54}(C_{54}(3,4,14,22))$. 

\item [\rm (b5)]  $Ad_{54}(C_{54}(2,16,20,21))$ = $\{C_{54}(2,16,20,21)$, $C_{54}(3,8,10,26)$, $C_{54}(4,14,15,22)\}$

\hfill = $Ad_{54}(C_{54}(3,8,10,26))$ = $Ad_{54}(C_{54}(4,14,15,22))$. 

\item [\rm (b6)]  $Ad_{54}(C_{54}(2,16,20,24))$ = $\{C_{54}(2,16,20,24)$, $C_{54}(8,10,12,26)$, $C_{54}(4,6,14,22)\}$

\hfill = $Ad_{54}(C_{54}(8,10,12,26))$ = $Ad_{54}(C_{54}(4,6,14,22))$. 
\end{enumerate} 

In problem \ref{p2}, we proved that circulant graphs 

(a) $C_{54}(1,3,17,19)$, $C_{54}(3,7,11,25)$ and $C_{54}(3,5,13,23)$; and

(b) $C_{54}(2,3,16,20)$, $C_{54}(3,4,14,22)$ and $C_{54}(3,8,10,26)$ are isomorphic of Type-2 w.r.t. $m$ = 3 = $\gcd(54, 3)$. In the next problem, we use remark \ref{r5} in the above triples of Type-2 isomorphic circulant graphs and also use the following different values of $s$ to obtain more triples of isomorphic circulant graphs and all Type-2 isomorphic circulant graphs of order 54.

\noindent
Different values of $s$ that are used to obtain more triples of isomorphic circulant graphs of order 54.
\begin{enumerate}
	\item [\rm (i)]  $\gcd(54, s)$ = 2 for $s$ = 2,4,8,10,14,16,20,22,26; 
	
	\item [\rm (ii)]  $\gcd(54, s)$ = 3 for $s$ = 3,15,21; 
	
	\item [\rm (iii)]  $\gcd(54, s)$ = 6 for $s$ = 6,12,24; 
	
	\item [\rm (iv)]  $\gcd(54, s)$ = 9 for $s$ = 9; 
	
	\item [\rm (v)]  $\gcd(54, s)$ = 18 for $s$ = 18; 
	
	\item [\rm (vi)]  $\gcd(54, s)$ = 27 for $s$ = 27. 
\end{enumerate}

 \begin{prm}\quad \label{p3} {\rm Show that each triple of circulant graphs of order 54 given below are isomorphic; classify their type of isomorphism and show that among them 960 triples are Type-2 isomorphic circulant graphs w.r.t. $m$ = 3. 

\noindent
{\bf {\footnotesize  (a) Triples of isomorphic circulant graphs of order 54 containing jump sizes 1,17,19, 7,11,25, 5,13,23.}}
\begin{enumerate} 	
\item [\rm (1)]  $C_{54}(1,3,17,19)$, $C_{54}(3,7,11,25)$, $C_{54}(3,5,13,23)$; 	

\item [\rm (2)]  $C_{54}(1,6,17,19)$, $C_{54}(6,7,11,25)$, $C_{54}(5,6,13,23)$; 	

\item [\rm (3)]  $C_{54}(1,9,17,19)$, $C_{54}(7,9,11,25)$, $C_{54}(5,9,13,23)$; 	

\item [\rm (4)]  $C_{54}(1,12,17,19)$, $C_{54}(7,11,12,25)$, $C_{54}(5,12,13,23)$; 	

\item [\rm (5)]  $C_{54}(1,15,17,19)$, $C_{54}(7,11,15,25)$, $C_{54}(5,13,15,23)$; 	

\item [\rm (6)]  $C_{54}(1,17,18,19)$, $C_{54}(7,11,18,25)$, $C_{54}(5,13,18,23)$; 	

\item [\rm (7)]  $C_{54}(1,17,19,21)$, $C_{54}(7,11,21,25)$, $C_{54}(5,13,21,23)$; 	

\item [\rm (8)]  $C_{54}(1,17,19,24)$, $C_{54}(7,11,24,25)$, $C_{54}(5,13,23,24)$; 

\item [\rm (9)]  $C_{54}(1,17,19,27)$, $C_{54}(7,11,25,27)$, $C_{54}(5,13,23,27)$; 

\item [\rm (10)]  $C_{54}(1,3,6,17,19)$, $C_{54}(3,6,7,11,25)$, $C_{54}(3,5,6,13,23)$; 	

\item [\rm (11)]  $C_{54}(1,3,9,17,19)$, $C_{54}(3,7,9,11,25)$, $C_{54}(3,5,9,13,23)$; 	

\item [\rm (12)]  $C_{54}(1,3,12,17,19)$, $C_{54}(3,7,11,12,25)$, $C_{54}(3,5,12,13,23)$; 	

\item [\rm (13)]  $C_{54}(1,3,15,17,19)$, $C_{54}(3,7,11,15,25)$, $C_{54}(3,5,13,15,23)$; 	

\item [\rm (14)]  $C_{54}(1,3,17,18,19)$, $C_{54}(3,7,11,18,25)$, $C_{54}(3,5,13,18,23)$; 	

\item [\rm (15)]  $C_{54}(1,3,17,19,21)$, $C_{54}(3,7,11,21,25)$, $C_{54}(3,5,13,21,23)$; 	

\item [\rm (16)]  $C_{54}(1,3,17,19,24)$, $C_{54}(3,7,11,24,25)$, $C_{54}(3,5,13,23,24)$; 	

\item [\rm (17)]  $C_{54}(1,3,17,19,27)$, $C_{54}(3,7,11,25,27)$, $C_{54}(3,5,13,23,27)$; 	

\item [\rm (18)]  $C_{54}(1,6,9,17,19)$, $C_{54}(6,7,9,11,25)$, $C_{54}(5,6,9,13,23)$; 	

\item [\rm (19)]  $C_{54}(1,6,12,17,19)$, $C_{54}(6,7,11,12,25)$, $C_{54}(5,6,12,13,23)$; 	

\item [\rm (20)]  $C_{54}(1,6,15,17,19)$, $C_{54}(6,7,11,15,25)$, $C_{54}(5,6,13,15,23)$; 	

\item [\rm (21)]  $C_{54}(1,6,17,18,19)$, $C_{54}(6,7,11,18,25)$, $C_{54}(5,6,13,18,23)$; 	

\item [\rm (22)]  $C_{54}(1,6,17,19,21)$, $C_{54}(6,7,11,21,25)$, $C_{54}(5,6,13,21,23)$; 	

\item [\rm (23)]  $C_{54}(1,6,17,19,24)$, $C_{54}(6,7,11,24,25)$, $C_{54}(5,6,13,23,24)$; 	

\item [\rm (24)]  $C_{54}(1,6,17,19,27)$, $C_{54}(6,7,11,25,27)$, $C_{54}(5,6,13,23,27)$; 	

\item [\rm (25)]  $C_{54}(1,9,12,17,19)$, $C_{54}(7,9,11,12,25)$, $C_{54}(5,9,12,13,23)$; 	

\item [\rm (26)]  $C_{54}(1,9,15,17,19)$, $C_{54}(7,9,11,15,25)$, $C_{54}(5,9,13,15,23)$; 	

\item [\rm (27)]  $C_{54}(1,9,17,18,19)$, $C_{54}(7,9,11,18,25)$, $C_{54}(5,9,13,18,23)$; 	

\item [\rm (28)]  $C_{54}(1,9,17,19,21)$, $C_{54}(7,9,11,21,25)$, $C_{54}(5,9,13,21,23)$; 	

\item [\rm (29)]  $C_{54}(1,9,17,19,24)$, $C_{54}(7,9,11,24,25)$, $C_{54}(5,9,13,23,24)$; 	

\item [\rm (30)]  $C_{54}(1,9,17,19,27)$, $C_{54}(7,9,11,25,27)$, $C_{54}(5,9,13,23,27)$; 	

\item [\rm (31)]  $C_{54}(1,12,15,17,19)$, $C_{54}(7,11,12,15,25)$, $C_{54}(5,12,13,15,23)$; 	

\item [\rm (32)]  $C_{54}(1,12,17,18,19)$, $C_{54}(7,11,12,18,25)$, $C_{54}(5,12,13,18,23)$; 	

\item [\rm (33)]  $C_{54}(1,12,17,19,21)$, $C_{54}(7,11,12,21,25)$, $C_{54}(5,12,13,21,23)$; 	

\item [\rm (34)]  $C_{54}(1,12,17,19,24)$, $C_{54}(7,11,12,24,25)$, $C_{54}(5,12,13,23,24)$; 	

\item [\rm (35)]  $C_{54}(1,12,17,19,27)$, $C_{54}(7,11,12,25,27)$, $C_{54}(5,12,13,23,27)$; 	

\item [\rm (36)]  $C_{54}(1,15,17,18,19)$, $C_{54}(7,11,15,18,25)$, $C_{54}(5,13,15,18,23)$; 	

\item [\rm (37)]  $C_{54}(1,15,17,19,21)$, $C_{54}(7,11,15,21,25)$, $C_{54}(5,13,15,21,23)$; 	

\item [\rm (38)]  $C_{54}(1,15,17,19,24)$, $C_{54}(7,11,15,24,25)$, $C_{54}(5,13,15,23,24)$; 	

\item [\rm (39)]  $C_{54}(1,15,17,19,27)$, $C_{54}(7,11,15,25,27)$, $C_{54}(5,13,15,23,27)$; 	

\item [\rm (40)]  $C_{54}(1,18,17,19,21)$, $C_{54}(7,11,18,21,25)$, $C_{54}(5,13,18,21,23)$; 	

\item [\rm (41)]  $C_{54}(1,18,17,19,24)$, $C_{54}(7,11,18,24,25)$, $C_{54}(5,13,18,23,24)$; 	

\item [\rm (42)]  $C_{54}(1,18,17,19,27)$, $C_{54}(7,11,18,25,27)$, $C_{54}(5,13,18,23,27)$; 	

\item [\rm (43)]  $C_{54}(1,17,19,21,24)$, $C_{54}(7,11,21,24,25)$, $C_{54}(5,13,21,23,24)$; 	

\item [\rm (44)]  $C_{54}(1,17,19,21,27)$, $C_{54}(7,11,21,25,27)$, $C_{54}(5,13,21,23,27)$; 	

\item [\rm (45)]  $C_{54}(1,17,19,24,27)$, $C_{54}(7,11,24,25,27)$, $C_{54}(5,13,23,24,27)$; 	

\item [\rm (46)]  $C_{54}(1,3,6,9,17,19)$, $C_{54}(3,6,7,9,11,25)$, $C_{54}(3,5,6,9,13,23)$; 	

\item [\rm (47)]  $C_{54}(1,3,6,12,17,19)$, $C_{54}(3,6,7,11,12,25)$, $C_{54}(3,5,6,12,13,23)$; 	

\item [\rm (48)]  $C_{54}(1,3,6,15,17,19)$, $C_{54}(3,6,7,11,15,25)$, $C_{54}(3,5,6,13,15,23)$; 	

\item [\rm (49)]  $C_{54}(1,3,6,17,18,19)$, $C_{54}(3,6,7,11,18,25)$, $C_{54}(3,5,6,13,18,23)$; 	

\item [\rm (50)]  $C_{54}(1,3,6,17,19,21)$, $C_{54}(3,6,7,11,21,25)$, $C_{54}(3,5,6,13,21,23)$; 	

\item [\rm (51)]  $C_{54}(1,3,6,17,19,24)$, $C_{54}(3,6,7,11,24,25)$, $C_{54}(3,5,6,13,23,24)$; 	

\item [\rm (52)]  $C_{54}(1,3,6,17,19,27)$, $C_{54}(3,6,7,11,25,27)$, $C_{54}(3,5,6,13,23,27)$; 	

\item [\rm (53)]  $C_{54}(1,3,9,12,17,19)$, $C_{54}(3,7,9,11,12,25)$, $C_{54}(3,5,9,12,13,23)$; 	

\item [\rm (54)]  $C_{54}(1,3,9,15,17,19)$, $C_{54}(3,7,9,11,15,25)$, $C_{54}(3,5,9,13,15,23)$; 	

\item [\rm (55)]  $C_{54}(1,3,9,17,18,19)$, $C_{54}(3,7,9,11,18,25)$, $C_{54}(3,5,9,13,18,23)$; 	

\item [\rm (56)]  $C_{54}(1,3,9,17,19,21)$, $C_{54}(3,7,9,11,21,25)$, $C_{54}(3,5,9,13,21,23)$; 	

\item [\rm (57)]  $C_{54}(1,3,9,17,19,24)$, $C_{54}(3,7,9,11,24,25)$, $C_{54}(3,5,9,13,23,24)$; 	

\item [\rm (58)]  $C_{54}(1,3,9,17,19,27)$, $C_{54}(3,7,9,11,25,27)$, $C_{54}(3,5,9,13,23,27)$; 	

\item [\rm (59)]  $C_{54}(1,3,12,15,17,19)$, $C_{54}(3,7,11,12,15,25)$, $C_{54}(3,5,12,13,15,23)$; 	

\item [\rm (60)]  $C_{54}(1,3,12,17,18,19)$, $C_{54}(3,7,11,12,18,25)$, $C_{54}(3,5,12,13,18,23)$; 	

\item [\rm (61)]  $C_{54}(1,3,12,17,19,21)$, $C_{54}(3,7,11,12,21,25)$, $C_{54}(3,5,12,13,21,23)$; 	

\item [\rm (62)]  $C_{54}(1,3,12,17,19,24)$, $C_{54}(3,7,11,12,24,25)$, $C_{54}(3,5,12,13,23,24)$; 	

\item [\rm (63)]  $C_{54}(1,3,12,17,19,27)$, $C_{54}(3,7,11,12,25,27)$, $C_{54}(3,5,12,13,23,27)$; 	

\item [\rm (64)]  $C_{54}(1,3,15,17,18,19)$, $C_{54}(3,7,11,15,18,25)$, $C_{54}(3,5,13,15,18,23)$; 	

\item [\rm (65)]  $C_{54}(1,3,15,17,19,21)$, $C_{54}(3,7,11,15,21,25)$, $C_{54}(3,5,13,15,21,23)$; 	

\item [\rm (66)]  $C_{54}(1,3,15,17,19,24)$, $C_{54}(3,7,11,15,24,25)$, $C_{54}(3,5,13,15,23,24)$; 	

\item [\rm (67)]  $C_{54}(1,3,15,17,19,27)$, $C_{54}(3,7,11,15,25,27)$, $C_{54}(3,5,13,15,23,27)$; 	

\item [\rm (68)]  $C_{54}(1,3,17,18,19,21)$, $C_{54}(3,7,11,18,21,25)$, $C_{54}(3,5,13,18,21,23)$; 	

\item [\rm (69)]  $C_{54}(1,3,17,18,19,24)$, $C_{54}(3,7,11,18,24,25)$, $C_{54}(3,5,13,18,23,24)$; 	

\item [\rm (70)]  $C_{54}(1,3,17,18,19,27)$, $C_{54}(3,7,11,18,25,27)$, $C_{54}(3,5,13,18,23,27)$; 	

\item [\rm (71)]  $C_{54}(1,3,17,19,21,24)$, $C_{54}(3,7,11,21,24,25)$, $C_{54}(3,5,13,21,23,24)$; 	

\item [\rm (72)]  $C_{54}(1,3,17,19,21,27)$, $C_{54}(3,7,11,21,25,27)$, $C_{54}(3,5,13,21,23,27)$; 	

\item [\rm (73)]  $C_{54}(1,3,17,19,24,27)$, $C_{54}(3,7,11,24,25,27)$, $C_{54}(3,5,13,23,24,27)$; 	

\item [\rm (74)]   $C_{54}(1,6,9,12,17,19)$, $C_{54}(6,7,9,11,12,25)$, $C_{54}(5,6,9,12,13,23)$; 	

\item [\rm (75)]   $C_{54}(1,6,9,15,17,19)$, $C_{54}(6,7,9,11,15,25)$, $C_{54}(5,6,9,13,15,23)$; 	

\item [\rm (76)]   $C_{54}(1,6,9,17,18,19)$, $C_{54}(6,7,9,11,18,25)$, $C_{54}(5,6,9,13,18,23)$; 	

\item [\rm (77)]   $C_{54}(1,6,9,17,19,21)$, $C_{54}(6,7,9,11,21,25)$, $C_{54}(5,6,9,13,21,23)$; 	

\item [\rm (78)]   $C_{54}(1,6,9,17,19,24)$, $C_{54}(6,7,9,11,24,25)$, $C_{54}(5,6,9,13,23,27)$; 	

\item [\rm (79)]   $C_{54}(1,6,9,17,19,27)$, $C_{54}(6,7,9,11,25,27)$, $C_{54}(5,6,9,13,23,27)$; 	

\item [\rm (80)]   $C_{54}(1,6,12,15,17,19)$, $C_{54}(6,7,11,12,15,25)$, $C_{54}(5,6,12,13,15,23)$; 	

\item [\rm (81)]   $C_{54}(1,6,12,17,18,19)$, $C_{54}(6,7,11,12,18,25)$, $C_{54}(5,6,12,13,18,23)$; 	

\item [\rm (82)]   $C_{54}(1,6,12,17,19,21)$, $C_{54}(6,7,11,12,21,25)$, $C_{54}(5,6,12,13,21,23)$; 	

\item [\rm (83)]   $C_{54}(1,6,12,17,19,24)$, $C_{54}(6,7,11,12,24,25)$, $C_{54}(5,6,12,13,23,24)$; 	

\item [\rm (84)]   $C_{54}(1,6,12,17,19,27)$, $C_{54}(6,7,11,12,25,27)$, $C_{54}(5,6,12,13,23,27)$; 	

\item [\rm (85)]  $C_{54}(1,6,15,17,18,19)$, $C_{54}(6,7,11,15,18,25)$, $C_{54}(5,6,13,15,18,23)$; 	

\item [\rm (86)]  $C_{54}(1,6,15,17,19,21)$, $C_{54}(6,7,11,15,21,25)$, $C_{54}(5,6,13,15,21,23)$; 	

\item [\rm (87)]  $C_{54}(1,6,15,17,19,24)$, $C_{54}(6,7,11,15,24,25)$, $C_{54}(5,6,13,15,23,24)$; 	

\item [\rm (88)]  $C_{54}(1,6,15,17,19,27)$, $C_{54}(6,7,11,15,25,27)$, $C_{54}(5,6,13,15,23,27)$; 	

\item [\rm (89)]  $C_{54}(1,6,17,18,19,21)$, $C_{54}(6,7,11,18,21,25)$, $C_{54}(5,6,13,18,21,23)$; 	

\item [\rm (90)]  $C_{54}(1,6,17,18,19,24)$, $C_{54}(6,7,11,18,24,25)$, $C_{54}(5,6,13,18,23,24)$; 	

\item [\rm (91)]  $C_{54}(1,6,17,18,19,27)$, $C_{54}(6,7,11,18,25,27)$, $C_{54}(5,6,13,18,23,27)$; 	

\item [\rm (92)]  $C_{54}(1,6,17,19,21,24)$, $C_{54}(6,7,11,21,24,25)$, $C_{54}(5,6,13,21,23,24)$; 	

\item [\rm (93)]  $C_{54}(1,6,17,19,21,27)$, $C_{54}(6,7,11,21,25,27)$, $C_{54}(5,6,13,21,23,27)$; 	

\item [\rm (94)]  $C_{54}(1,6,17,19,24,27)$, $C_{54}(6,7,11,24,25,27)$, $C_{54}(5,6,13,23,24,27)$; 	

\item [\rm (95)]    $C_{54}(1,9,12,15,17,19)$, $C_{54}(7,9,11,12,15,25)$, $C_{54}(5,9,12,13,15,23)$; 	

\item [\rm (96)]    $C_{54}(1,9,12,17,18,19)$, $C_{54}(7,9,11,12,18,25)$, $C_{54}(5,9,12,13,18,23)$; 	

\item [\rm (97)]    $C_{54}(1,9,12,17,19,21)$, $C_{54}(7,9,11,12,21,25)$, $C_{54}(5,9,12,13,21,23)$; 	

\item [\rm (98)]    $C_{54}(1,9,12,17,19,24)$, $C_{54}(7,9,11,12,24,25)$, $C_{54}(5,9,12,13,23,24)$; 	

\item [\rm (99)]    $C_{54}(1,9,12,17,19,27)$, $C_{54}(7,9,11,12,25,27)$, $C_{54}(5,9,12,13,23,27)$; 	

\item [\rm (100)]    $C_{54}(1,9,15,17,18,19)$, $C_{54}(7,9,11,15,18,25)$, $C_{54}(5,9,13,15,18,23)$; 	

\item [\rm (101)]    $C_{54}(1,9,15,17,19,21)$, $C_{54}(7,9,11,15,21,25)$, $C_{54}(5,9,13,15,21,23)$; 	

\item [\rm (102)]    $C_{54}(1,9,15,17,19,24)$, $C_{54}(7,9,11,15,24,25)$, $C_{54}(5,9,13,15,23,24)$; 	

\item [\rm (103)]    $C_{54}(1,9,15,17,19,27)$, $C_{54}(7,9,11,15,25,27)$, $C_{54}(5,9,13,15,23,27)$; 	

\item [\rm (104)]    $C_{54}(1,9,17,18,19,21)$, $C_{54}(7,9,11,18,21,25)$, $C_{54}(5,9,13,18,21,23)$; 	

\item [\rm (105)]    $C_{54}(1,9,17,18,19,24)$, $C_{54}(7,9,11,18,24,25)$, $C_{54}(5,9,13,18,23,24)$; 	

\item [\rm (106)]    $C_{54}(1,9,17,18,19,27)$, $C_{54}(7,9,11,18,25,27)$, $C_{54}(5,9,13,18,23,27)$; 	

\item [\rm (107)]    $C_{54}(1,9,17,19,21,24)$, $C_{54}(7,9,11,21,24,25)$, $C_{54}(5,9,13,21,23,24)$; 	

\item [\rm (108)]    $C_{54}(1,9,17,19,21,27)$, $C_{54}(7,9,11,21,25,27)$, $C_{54}(5,9,13,21,23,27)$; 	

\item [\rm (109)]    $C_{54}(1,9,17,19,24,27)$, $C_{54}(7,9,11,24,25,27)$, $C_{54}(5,9,13,23,24,27)$; 	

\item [\rm (110)]  $C_{54}(1,12,15,17,18,19)$, $C_{54}(7,11,12,15,18,25)$, $C_{54}(5,12,13,15,18,23)$; 	

\item [\rm (111)]  $C_{54}(1,12,15,17,19,21)$, $C_{54}(7,11,12,15,21,25)$, $C_{54}(5,12,13,15,21,23)$; 	

\item [\rm (112)]  $C_{54}(1,12,15,17,19,24)$, $C_{54}(7,11,12,15,24,25)$, $C_{54}(5,12,13,15,23,24)$; 	

\item [\rm (113)]  $C_{54}(1,12,15,17,19,27)$, $C_{54}(7,11,12,15,25,27)$, $C_{54}(5,12,13,15,23,27)$; 	

\item [\rm (114)]  $C_{54}(1,12,17,18,19,21)$, $C_{54}(7,11,12,18,21,25)$, $C_{54}(5,12,13,18,21,23)$; 	

\item [\rm (115)]  $C_{54}(1,12,17,18,19,24)$, $C_{54}(7,11,12,18,24,25)$, $C_{54}(5,12,13,18,23,24)$; 	

\item [\rm (116)]  $C_{54}(1,12,17,18,19,27)$, $C_{54}(7,11,12,18,25,27)$, $C_{54}(5,12,13,18,23,27)$; 	

\item [\rm (117)]  $C_{54}(1,12,17,19,21,24)$, $C_{54}(7,11,12,21,24,25)$, $C_{54}(5,12,13,21,23,24)$; 	

\item [\rm (118)]  $C_{54}(1,12,17,19,21,27)$, $C_{54}(7,11,12,21,25,27)$, $C_{54}(5,12,13,21,23,27)$; 	

\item [\rm (119)]  $C_{54}(1,12,17,19,24,27)$, $C_{54}(7,11,12,24,25,27)$, $C_{54}(5,12,13,23,24,27)$; 	

\item [\rm (120)]  $C_{54}(1,15,17,18,19,21)$, $C_{54}(7,11,15,18,21,25)$, $C_{54}(5,13,15,18,21,23)$; 	

\item [\rm (121)]  $C_{54}(1,15,17,18,19,24)$, $C_{54}(7,11,15,18,24,25)$, $C_{54}(5,13,15,18,23,24)$; 	

\item [\rm (122)]  $C_{54}(1,15,17,18,19,27)$, $C_{54}(7,11,15,18,25,27)$, $C_{54}(5,13,15,18,23,27)$; 	

\item [\rm (123)]  $C_{54}(1,15,17,19,21,24)$, $C_{54}(7,11,15,21,24,25)$, $C_{54}(5,13,15,21,23,24)$; 	

\item [\rm (124)]  $C_{54}(1,15,17,19,21,27)$, $C_{54}(7,11,15,21,25,27)$, $C_{54}(5,13,15,21,23,27)$; 	

\item [\rm (125)]  $C_{54}(1,15,17,19,24,27)$, $C_{54}(7,11,15,24,25,27)$, $C_{54}(5,13,15,23,24,27)$; 	

\item [\rm (126)]  $C_{54}(1,17,18,19,21,24)$, $C_{54}(7,11,18,21,24,25)$, $C_{54}(5,13,18,21,23,24)$; 	

\item [\rm (127)]  $C_{54}(1,17,18,19,21,27)$, $C_{54}(7,11,18,21,25,27)$, $C_{54}(5,13,18,21,23,27)$; 	

\item [\rm (128)]  $C_{54}(1,17,18,19,24,27)$, $C_{54}(7,11,18,24,25,27)$, $C_{54}(5,13,18,23,24,27)$; 	

\item [\rm (129)]  $C_{54}(1,17,19,21,24,27)$, $C_{54}(7,11,21,24,25,27)$, $C_{54}(5,13,21,23,24,27)$; 	

\item [\rm (130)]  $C_{54}(1,3,6,9,12,17,19)$, $C_{54}(3,6,7,9,11,12,25)$, $C_{54}(3,5,6,9,12,13,23)$; 	

\item [\rm (131)]  $C_{54}(1,3,6,9,15,17,19)$, $C_{54}(3,6,7,9,11,15,25)$, $C_{54}(3,5,6,9,13,15,23)$; 	

\item [\rm (132)]  $C_{54}(1,3,6,9,17,18,19)$, $C_{54}(3,6,7,9,11,18,25)$, $C_{54}(3,5,6,9,13,18,23)$; 	

\item [\rm (133)]  $C_{54}(1,3,6,9,17,19,21)$, $C_{54}(3,6,7,9,11,21,25)$, $C_{54}(3,5,6,9,13,21,23)$; 	

\item [\rm (134)]  $C_{54}(1,3,6,9,17,19,24)$, $C_{54}(3,6,7,9,11,24,25)$, $C_{54}(3,5,6,9,13,23,24)$; 	

\item [\rm (135)]  $C_{54}(1,3,6,9,17,19,27)$, $C_{54}(3,6,7,9,11,25,27)$, $C_{54}(3,5,6,9,13,23,27)$; 	

\item [\rm (136)]  $C_{54}(1,3,6,12,15,17,19)$, $C_{54}(3,6,7,11,12,15,25)$, $C_{54}(3,5,6,12,13,15,23)$; 	

\item [\rm (137)]  $C_{54}(1,3,6,12,17,18,19)$, $C_{54}(3,6,7,11,12,18,25)$, $C_{54}(3,5,6,12,13,18,23)$; 	

\item [\rm (138)]  $C_{54}(1,3,6,12,17,19,21)$, $C_{54}(3,6,7,11,12,21,25)$, $C_{54}(3,5,6,12,13,21,23)$; 	

\item [\rm (139)]  $C_{54}(1,3,6,12,17,19,24)$, $C_{54}(3,6,7,11,12,24,25)$, $C_{54}(3,5,6,12,13,23,24)$; 	

\item [\rm (140)]  $C_{54}(1,3,6,12,17,19,27)$, $C_{54}(3,6,7,11,12,25,27)$, $C_{54}(3,5,6,12,13,23,27)$; 	

\item [\rm (141)]  $C_{54}(1,3,6,15,17,18,19)$, $C_{54}(3,6,7,11,15,18,25)$, $C_{54}(3,5,6,13,15,18,23)$; 	

\item [\rm (142)]  $C_{54}(1,3,6,15,17,19,21)$, $C_{54}(3,6,7,11,15,21,25)$, $C_{54}(3,5,6,13,15,21,23)$; 	

\item [\rm (143)]  $C_{54}(1,3,6,15,17,19,24)$, $C_{54}(3,6,7,11,15,24,25)$, $C_{54}(3,5,6,13,15,23,24)$; 	

\item [\rm (144)]  $C_{54}(1,3,6,15,17,19,27)$, $C_{54}(3,6,7,11,15,25,27)$, $C_{54}(3,5,6,13,15,23,27)$; 	

\item [\rm (145)]  $C_{54}(1,3,6,17,18,19,21)$, $C_{54}(3,6,7,11,18,21,25)$, $C_{54}(3,5,6,13,18,21,23)$; 	

\item [\rm (146)]  $C_{54}(1,3,6,17,18,19,24)$, $C_{54}(3,6,7,11,18,24,25)$, $C_{54}(3,5,6,13,18,23,24)$; 	

\item [\rm (147)]  $C_{54}(1,3,6,17,18,19,27)$, $C_{54}(3,6,7,11,18,25,27)$, $C_{54}(3,5,6,13,18,23,27)$; 	

\item [\rm (148)]  $C_{54}(1,3,6,17,19,21,24)$, $C_{54}(3,6,7,11,21,24,25)$, $C_{54}(3,5,6,13,21,23,24)$; 	

\item [\rm (149)]  $C_{54}(1,3,6,17,19,21,27)$, $C_{54}(3,6,7,11,21,25,27)$, $C_{54}(3,5,6,13,21,23,27)$; 	

\item [\rm (150)]  $C_{54}(1,3,6,17,19,24,27)$, $C_{54}(3,6,7,11,24,25,27)$, $C_{54}(3,5,6,13,23,24,27)$; 	

\item [\rm (151)]  $C_{54}(1,3,9,12,15,17,19)$, $C_{54}(3,7,9,11,12,15,25)$, $C_{54}(3,5,9,12,13,15,23)$; 	

\item [\rm (152)]  $C_{54}(1,3,9,12,17,18,19)$, $C_{54}(3,7,9,11,12,18,25)$, $C_{54}(3,5,9,12,13,18,23)$; 	

\item [\rm (153)]  $C_{54}(1,3,9,12,17,19,21)$, $C_{54}(3,7,9,11,12,21,25)$, $C_{54}(3,5,9,12,13,21,23)$; 	

\item [\rm (154)]  $C_{54}(1,3,9,12,17,19,24)$, $C_{54}(3,7,9,11,12,24,25)$, $C_{54}(3,5,9,12,13,23,24)$; 	

\item [\rm (155)]  $C_{54}(1,3,9,12,17,19,27)$, $C_{54}(3,7,9,11,12,25,27)$, $C_{54}(3,5,9,12,13,23,27)$; 	

\item [\rm (156)]  $C_{54}(1,3,9,15,17,18,19)$, $C_{54}(3,7,9,11,15,18,25)$, $C_{54}(3,5,9,13,15,18,23)$; 	

\item [\rm (157)]  $C_{54}(1,3,9,15,17,19,21)$, $C_{54}(3,7,9,11,15,21,25)$, $C_{54}(3,5,9,13,15,21,23)$; 	

\item [\rm (158)]  $C_{54}(1,3,9,15,17,19,24)$, $C_{54}(3,7,9,11,15,24,25)$, $C_{54}(3,5,9,13,15,23,24)$; 	

\item [\rm (159)]  $C_{54}(1,3,9,15,17,19,27)$, $C_{54}(3,7,9,11,15,25,27)$, $C_{54}(3,5,9,13,15,23,27)$; 	

\item [\rm (160)]  $C_{54}(1,3,9,17,18,19,21)$, $C_{54}(3,7,9,11,18,21,25)$, $C_{54}(3,5,9,13,18,21,23)$; 	

\item [\rm (161)]  $C_{54}(1,3,9,17,18,19,24)$, $C_{54}(3,7,9,11,18,24,25)$, $C_{54}(3,5,9,13,18,23,24)$; 	

\item [\rm (162)]  $C_{54}(1,3,9,17,18,19,27)$, $C_{54}(3,7,9,11,18,25,27)$, $C_{54}(3,5,9,13,18,23,27)$; 	

\item [\rm (163)]  $C_{54}(1,3,9,17,19,21,24)$, $C_{54}(3,7,9,11,21,24,25)$, $C_{54}(3,5,9,13,21,23,24)$; 	

\item [\rm (164)]  $C_{54}(1,3,9,17,19,21,27)$, $C_{54}(3,7,9,11,21,25,27)$, $C_{54}(3,5,9,13,21,23,27)$; 	

\item [\rm (165)]  $C_{54}(1,3,9,17,19,24,27)$, $C_{54}(3,7,9,11,24,25,27)$, $C_{54}(3,5,9,13,23,24,27)$; 	

\item [\rm (166)]  $C_{54}(1,3,12,15,17,18,19)$, $C_{54}(3,7,11,12,15,18,25)$, $C_{54}(3,5,12,13,15,18,23)$; 	

\item [\rm (167)]  $C_{54}(1,3,12,15,17,19,21)$, $C_{54}(3,7,11,12,15,21,25)$, $C_{54}(3,5,12,13,15,21,23)$; 	

\item [\rm (168)]  $C_{54}(1,3,12,15,17,19,24)$, $C_{54}(3,7,11,12,15,24,25)$, $C_{54}(3,5,12,13,15,23,24)$; 	

\item [\rm (169)]  $C_{54}(1,3,12,15,17,19,27)$, $C_{54}(3,7,11,12,15,25,27)$, $C_{54}(3,5,12,13,15,23,27)$; 	

\item [\rm (170)]  $C_{54}(1,3,12,17,18,19,21)$, $C_{54}(3,7,11,12,18,21,25)$, $C_{54}(3,5,12,13,18,21,23)$; 	

\item [\rm (171)]  $C_{54}(1,3,12,17,18,19,24)$, $C_{54}(3,7,11,12,18,24,25)$, $C_{54}(3,5,12,13,18,23,24)$; 	

\item [\rm (172)]  $C_{54}(1,3,12,17,18,19,27)$, $C_{54}(3,7,11,12,18,25,27)$, $C_{54}(3,5,12,13,18,23,27)$; 	

\item [\rm (173)]  $C_{54}(1,3,12,17,19,21,24)$, $C_{54}(3,7,11,12,21,24,25)$, $C_{54}(3,5,12,13,21,23,24)$; 	

\item [\rm (174)]  $C_{54}(1,3,12,17,19,21,27)$, $C_{54}(3,7,11,12,21,25,27)$, $C_{54}(3,5,12,13,21,23,27)$; 	

\item [\rm (175)]  $C_{54}(1,3,12,17,19,24,27)$, $C_{54}(3,7,11,12,24,25,27)$, $C_{54}(3,5,12,13,23,24,27)$; 	

\item [\rm (176)]  $C_{54}(1,3,15,17,18,19,21)$, $C_{54}(3,7,11,15,18,21,25)$, $C_{54}(3,5,13,15,18,21,23)$; 	

\item [\rm (177)]  $C_{54}(1,3,15,17,18,19,24)$, $C_{54}(3,7,11,15,18,24,25)$, $C_{54}(3,5,13,15,18,23,24)$; 	

\item [\rm (178)]  $C_{54}(1,3,15,17,18,19,27)$, $C_{54}(3,7,11,15,18,25,27)$, $C_{54}(3,5,13,15,18,23,27)$; 	

\item [\rm (179)]  $C_{54}(1,3,15,17,19,21,24)$, $C_{54}(3,7,11,15,21,24,25)$, $C_{54}(3,5,13,15,21,23,24)$; 	

\item [\rm (180)]  $C_{54}(1,3,15,17,19,21,27)$, $C_{54}(3,7,11,15,21,25,27)$, $C_{54}(3,5,13,15,21,23,27)$; 	

\item [\rm (181)]  $C_{54}(1,3,15,17,19,24,27)$, $C_{54}(3,7,11,15,24,25,27)$, $C_{54}(3,5,13,15,23,24,27)$; 	

\item [\rm (182)]  $C_{54}(1,3,17,18,19,21,24)$, $C_{54}(3,7,11,18,21,24,25)$, $C_{54}(3,5,13,18,21,23,24)$; 	

\item [\rm (183)]  $C_{54}(1,3,17,18,19,21,27)$, $C_{54}(3,7,11,18,21,25,27)$, $C_{54}(3,5,13,18,21,23,27)$; 	

\item [\rm (184)]  $C_{54}(1,3,17,18,19,24,27)$, $C_{54}(3,7,11,18,24,25,27)$, $C_{54}(3,5,13,18,23,24,27)$; 	

\item [\rm (185)]  $C_{54}(1,3,17,19,21,24,27)$, $C_{54}(3,7,11,21,24,25,27)$, $C_{54}(3,5,13,21,23,24,27)$; 	

\item [\rm (186)]  $C_{54}(1,6,9,12,15,17,19)$, $C_{54}(6,7,9,11,12,15,25)$, $C_{54}(5,6,9,12,13,15,23)$; 	

\item [\rm (187)]  $C_{54}(1,6,9,12,17,18,19)$, $C_{54}(6,7,9,11,12,18,25)$, $C_{54}(5,6,9,12,13,18,23)$; 	

\item [\rm (188)]  $C_{54}(1,6,9,12,17,19,21)$, $C_{54}(6,7,9,11,12,21,25)$, $C_{54}(5,6,9,12,13,21,23)$; 	

\item [\rm (189)]  $C_{54}(1,6,9,12,17,19,24)$, $C_{54}(6,7,9,11,12,24,25)$, $C_{54}(5,6,9,12,13,23,24)$; 	

\item [\rm (190)]  $C_{54}(1,6,9,12,17,19,27)$, $C_{54}(6,7,9,11,12,25,27)$, $C_{54}(5,6,9,12,13,23,27)$; 	

\item [\rm (191)]  $C_{54}(1,6,9,15,17,18,19)$, $C_{54}(6,7,9,11,15,18,25)$, $C_{54}(5,6,9,13,15,18,23)$; 	

\item [\rm (192)]  $C_{54}(1,6,9,15,17,19,21)$, $C_{54}(6,7,9,11,15,21,25)$, $C_{54}(5,6,9,13,15,21,23)$; 	

\item [\rm (193)]  $C_{54}(1,6,9,15,17,19,24)$, $C_{54}(6,7,9,11,15,24,25)$, $C_{54}(5,6,9,13,15,23,24)$; 	

\item [\rm (194)]  $C_{54}(1,6,9,15,17,19,27)$, $C_{54}(6,7,9,11,15,25,27)$, $C_{54}(5,6,9,13,15,23,27)$; 	

\item [\rm (195)]  $C_{54}(1,6,9,17,18,19,21)$, $C_{54}(6,7,9,11,18,21,25)$, $C_{54}(5,6,9,13,18,21,23)$; 	

\item [\rm (196)]  $C_{54}(1,6,9,17,18,19,24)$, $C_{54}(6,7,9,11,18,24,25)$, $C_{54}(5,6,9,13,18,23,24)$; 	

\item [\rm (197)]  $C_{54}(1,6,9,17,18,19,27)$, $C_{54}(6,7,9,11,18,25,27)$, $C_{54}(5,6,9,13,18,23,27)$; 	

\item [\rm (198)]  $C_{54}(1,6,9,17,19,21,24)$, $C_{54}(6,7,9,11,21,24,25)$, $C_{54}(5,6,9,13,21,23,24)$; 	

\item [\rm (199)]  $C_{54}(1,6,9,17,19,21,27)$, $C_{54}(6,7,9,11,21,25,27)$, $C_{54}(5,6,9,13,21,23,27)$; 	

\item [\rm (200)]  $C_{54}(1,6,9,17,19,24,27)$, $C_{54}(6,7,9,11,24,25,27)$, $C_{54}(5,6,9,13,23,24,27)$; 	

\item [\rm (201)]  $C_{54}(1,6,12,15,17,18,19)$, $C_{54}(6,7,11,12,15,18,25)$, $C_{54}(5,6,12,13,15,18,23)$; 	

\item [\rm (202)]  $C_{54}(1,6,12,15,17,19,21)$, $C_{54}(6,7,11,12,15,21,25)$, $C_{54}(5,6,12,13,15,21,23)$; 	

\item [\rm (203)]  $C_{54}(1,6,12,15,17,19,24)$, $C_{54}(6,7,11,12,15,24,25)$, $C_{54}(5,6,12,13,15,23,24)$; 	

\item [\rm (204)]  $C_{54}(1,6,12,15,17,19,27)$, $C_{54}(6,7,11,12,15,25,27)$, $C_{54}(5,6,12,13,15,23,27)$; 	

\item [\rm (205)]  $C_{54}(1,6,12,17,18,19,21)$, $C_{54}(6,7,11,12,18,21,25)$, $C_{54}(5,6,12,13,18,21,23)$; 	

\item [\rm (206)]  $C_{54}(1,6,12,17,18,19,24)$, $C_{54}(6,7,11,12,18,24,25)$, $C_{54}(5,6,12,13,18,23,24)$; 	

\item [\rm (207)]  $C_{54}(1,6,12,17,18,19,27)$, $C_{54}(6,7,11,12,18,25,27)$, $C_{54}(5,6,12,13,18,23,27)$; 	

\item [\rm (208)]  $C_{54}(1,6,12,17,19,21,24)$, $C_{54}(6,7,11,12,21,24,25)$, $C_{54}(5,6,12,13,21,23,24)$; 	

\item [\rm ( 209)]  $C_{54}(1,6,12,17,19,21,27)$, $C_{54}(6,7,11,12,21,25,27)$, $C_{54}(5,6,12,13,21,23,27)$; 	

\item [\rm ( 210)]  $C_{54}(1,6,12,17,19,24,27)$, $C_{54}(6,7,11,12,24,25,27)$, $C_{54}(5,6,12,13,23,24,27)$; 	

\item [\rm ( 211)]  $C_{54}(1,6,15,17,18,19,21)$, $C_{54}(6,7,11,15,18,21,25)$, $C_{54}(5,6,13,15,18,21,23)$; 	

\item [\rm ( 212)]  $C_{54}(1,6,15,17,18,19,24)$, $C_{54}(6,7,11,15,18,24,25)$, $C_{54}(5,6,13,15,18,23,24)$; 	

\item [\rm ( 213)]  $C_{54}(1,6,15,17,18,19,27)$, $C_{54}(6,7,11,15,18,25,27)$, $C_{54}(5,6,13,15,18,23,27)$; 	

\item [\rm ( 214)]  $C_{54}(1,6,15,17,19,21,24)$, $C_{54}(6,7,11,15,21,24,25)$, $C_{54}(5,6,13,15,21,23,24)$; 	

\item [\rm ( 215)]  $C_{54}(1,6,15,17,19,21,27)$, $C_{54}(6,7,11,15,21,25,27)$, $C_{54}(5,6,13,15,21,23,27)$; 	

\item [\rm ( 216)]  $C_{54}(1,6,15,17,19,24,27)$, $C_{54}(6,7,11,15,24,25,27)$, $C_{54}(5,6,13,15,23,24,27)$; 	

\item [\rm ( 217)]  $C_{54}(1,6,17,18,19,21,24)$, $C_{54}(6,7,11,18,21,24,25)$, $C_{54}(5,6,13,18,21,23,24)$; 	

\item [\rm ( 218)]  $C_{54}(1,6,17,18,19,21,27)$, $C_{54}(6,7,11,18,21,25,27)$, $C_{54}(5,6,13,18,21,23,27)$; 	

\item [\rm ( 219)]  $C_{54}(1,6,17,18,19,24,27)$, $C_{54}(6,7,11,18,24,25,27)$, $C_{54}(5,6,13,18,23,24,27)$; 	

\item [\rm ( 220)]  $C_{54}(1,6,17,19,21,24,27)$, $C_{54}(6,7,11,21,24,25,27)$, $C_{54}(5,6,13,21,23,24,27)$; 	

\item [\rm ( 221)]  $C_{54}(1,9,12,15,17,18,19)$, $C_{54}(7,9,11,12,15,18,25)$, $C_{54}(5,9,12,13,15,18,23)$; 	

\item [\rm ( 222)]  $C_{54}(1,9,12,15,17,19,21)$, $C_{54}(7,9,11,12,15,21,25)$, $C_{54}(5,9,12,13,15,21,23)$; 	

\item [\rm ( 223)]  $C_{54}(1,9,12,15,17,19,24)$, $C_{54}(7,9,11,12,15,24,25)$, $C_{54}(5,9,12,13,15,23,24)$; 	

\item [\rm ( 224)]  $C_{54}(1,9,12,15,17,19,27)$, $C_{54}(7,9,11,12,15,25,27)$, $C_{54}(5,9,12,13,15,23,27)$; 	

\item [\rm ( 225)]  $C_{54}(1,9,12,17,18,19,21)$, $C_{54}(7,9,11,12,18,21,25)$, $C_{54}(5,9,12,13,18,21,23)$; 	

\item [\rm ( 226)]  $C_{54}(1,9,12,17,18,19,24)$, $C_{54}(7,9,11,12,18,24,25)$, $C_{54}(5,9,12,13,18,23,24)$; 	

\item [\rm ( 227)]  $C_{54}(1,9,12,17,18,19,27)$, $C_{54}(7,9,11,12,18,25,27)$, $C_{54}(5,9,12,13,18,23,27)$; 	

\item [\rm ( 228)]  $C_{54}(1,9,12,17,19,21,24)$, $C_{54}(7,9,11,12,21,24,25)$, $C_{54}(5,9,12,13,21,23,24)$; 	

\item [\rm ( 229)]  $C_{54}(1,9,12,17,19,21,27)$, $C_{54}(7,9,11,12,21,25,27)$, $C_{54}(5,9,12,13,21,23,27)$; 	

\item [\rm ( 230)]  $C_{54}(1,9,12,17,19,24,27)$, $C_{54}(7,9,11,12,24,25,27)$, $C_{54}(5,9,12,13,23,24,27)$; 	

\item [\rm ( 231)]  $C_{54}(1,9,15,17,18,19,21)$, $C_{54}(7,9,11,15,18,21,25)$, $C_{54}(5,9,13,15,18,21,23)$; 	

\item [\rm ( 232)]  $C_{54}(1,9,15,17,18,19,24)$, $C_{54}(7,9,11,15,18,24,25)$, $C_{54}(5,9,13,15,18,23,24)$; 	

\item [\rm ( 233)]  $C_{54}(1,9,15,17,18,19,27)$, $C_{54}(7,9,11,15,18,25,27)$, $C_{54}(5,9,13,15,18,23,27)$; 	

\item [\rm ( 234)]  $C_{54}(1,9,15,17,19,21,24)$, $C_{54}(7,9,11,15,21,24,25)$, $C_{54}(5,9,13,15,21,23,24)$; 	

\item [\rm ( 235)]  $C_{54}(1,9,15,17,19,21,27)$, $C_{54}(7,9,11,15,21,25,27)$, $C_{54}(5,9,13,15,21,23,27)$; 	

\item [\rm ( 236)]  $C_{54}(1,9,15,17,19,24,27)$, $C_{54}(7,9,11,15,24,25,27)$, $C_{54}(5,9,13,15,23,24,27)$; 	

\item [\rm ( 237)]  $C_{54}(1,9,17,18,19,21,24)$, $C_{54}(7,9,11,18,21,24,25)$, $C_{54}(5,9,13,18,21,23,24)$; 	

\item [\rm ( 238)]  $C_{54}(1,9,17,18,19,21,27)$, $C_{54}(7,9,11,18,21,25,27)$, $C_{54}(5,9,13,18,21,23,27)$; 	

\item [\rm ( 239)]  $C_{54}(1,9,17,18,19,24,27)$, $C_{54}(7,9,11,18,24,25,27)$, $C_{54}(5,9,13,18,23,24,27)$; 	

\item [\rm ( 240)]  $C_{54}(1,9,17,19,21,24,27)$, $C_{54}(7,9,11,21,24,25,27)$, $C_{54}(5,9,13,21,23,24,27)$; 	

\item [\rm ( 241)]  $C_{54}(1,12,15,17,18,19,21)$, $C_{54}(7,11,12,15,18,21,25)$, $C_{54}(5,12,13,15,18,21,23)$; 	

\item [\rm ( 242)]  $C_{54}(1,12,15,17,18,19,24)$, $C_{54}(7,11,12,15,18,24,25)$, $C_{54}(5,12,13,15,18,23,24)$; 	

\item [\rm ( 243)]  $C_{54}(1,12,15,17,18,19,27)$, $C_{54}(7,11,12,15,18,25,27)$, $C_{54}(5,12,13,15,18,23,27)$; 	

\item [\rm ( 244)]  $C_{54}(1,12,15,17,19,21,24)$, $C_{54}(7,11,12,15,21,24,25)$, $C_{54}(5,12,13,15,21,23,24)$; 	

\item [\rm ( 245)]  $C_{54}(1,12,15,17,19,21,27)$, $C_{54}(7,11,12,15,21,25,27)$, $C_{54}(5,12,13,15,21,23,27)$; 	

\item [\rm ( 246)]  $C_{54}(1,12,15,17,19,24,27)$, $C_{54}(7,11,12,15,24,25,27)$, $C_{54}(5,12,13,15,23,24,27)$; 	

\item [\rm ( 247)]  $C_{54}(1,12,17,18,19,21,24)$, $C_{54}(7,11,12,18,21,24,25)$, $C_{54}(5,12,13,18,21,23,24)$; 	

\item [\rm ( 248)]  $C_{54}(1,12,17,18,19,21,27)$, $C_{54}(7,11,12,18,21,25,27)$, $C_{54}(5,12,13,18,21,23,27)$; 	

\item [\rm ( 249)]  $C_{54}(1,12,17,18,19,24,27)$, $C_{54}(7,11,12,18,24,25,27)$, $C_{54}(5,12,13,18,23,24,27)$; 	

\item [\rm ( 250)]  $C_{54}(1,12,17,19,21,24,27)$, $C_{54}(7,11,12,21,24,25,27)$, $C_{54}(5,12,13,21,23,24,27)$; 	

\item [\rm (251)]   $C_{54}(1,15,17,18,19,21,24)$, $C_{54}(7,11,15,18,21,24,25)$, $C_{54}(5,13,15,18,21,23,24)$; 	

\item [\rm (252)]   $C_{54}(1,15,17,18,19,21,27)$, $C_{54}(7,11,15,18,21,25,27)$, $C_{54}(5,13,15,18,21,23,27)$; 	

\item [\rm (253)]  $C_{54}(1,15,17,18,19,24,27)$, $C_{54}(7,11,15,18,24,25,27)$, $C_{54}(5,13,15,18,23,24,27)$; 	

\item [\rm (254)]  $C_{54}(1,15,17,19,21,24,27)$, $C_{54}(7,11,15,21,24,25,27)$, $C_{54}(5,13,15,21,23,24,27)$; 	

\item [\rm (255)]  $C_{54}(1,17,18,19,21,24,27)$, $C_{54}(7,11,18,21,24,25,27)$, $C_{54}(5,13,18,21,23,24,27)$; 	

\item [\rm (256)]  $C_{54}(1,3,6,9,12,15,17,19)$, $C_{54}(3,6,7,9,11,12,15,25)$, $C_{54}(3,5,6,9,12,13,15,23)$; 	

\item [\rm (257)]  $C_{54}(1,3,6,9,12,17,18,19)$, $C_{54}(3,6,7,9,11,12,18,25)$, $C_{54}(3,5,6,9,12,13,18,23)$; 	

\item [\rm (258)]  $C_{54}(1,3,6,9,12,17,19,21)$, $C_{54}(3,6,7,9,11,12,21,25)$, $C_{54}(3,5,6,9,12,13,21,23)$; 	

\item [\rm (259)]  $C_{54}(1,3,6,9,12,17,19,24)$, $C_{54}(3,6,7,9,11,12,24,25)$, $C_{54}(3,5,6,9,12,13,23,24)$; 	

\item [\rm (260)]  $C_{54}(1,3,6,9,12,17,19,27)$, $C_{54}(3,6,7,9,11,12,25,27)$, $C_{54}(3,5,6,9,12,13,23,27)$; 	

\item [\rm (261)]  $C_{54}(1,3,6,9,15,17,18,19)$, $C_{54}(3,6,7,9,11,15,18,25)$, $C_{54}(3,5,6,9,13,15,18,23)$; 	

\item [\rm (262)]  $C_{54}(1,3,6,9,15,17,19,21)$, $C_{54}(3,6,7,9,11,15,21,25)$, $C_{54}(3,5,6,9,13,15,21,23)$; 	

\item [\rm (263)]  $C_{54}(1,3,6,9,15,17,19,24)$, $C_{54}(3,6,7,9,11,15,24,25)$, $C_{54}(3,5,6,9,13,15,23,24)$; 	

\item [\rm (264)]  $C_{54}(1,3,6,9,15,17,19,27)$, $C_{54}(3,6,7,9,11,15,25,27)$, $C_{54}(3,5,6,9,13,15,23,27)$; 	

\item [\rm (265)]  $C_{54}(1,3,6,9,17,18,19,21)$, $C_{54}(3,6,7,9,11,18,21,25)$, $C_{54}(3,5,6,9,13,18,21,23)$; 	

\item [\rm (266)]  $C_{54}(1,3,6,9,17,18,19,24)$, $C_{54}(3,6,7,9,11,18,24,25)$, $C_{54}(3,5,6,9,13,18,23,24)$; 	

\item [\rm (267)]  $C_{54}(1,3,6,9,17,18,19,27)$, $C_{54}(3,6,7,9,11,18,25,27)$, $C_{54}(3,5,6,9,13,18,23,27)$; 	

\item [\rm (268)]  $C_{54}(1,3,6,9,17,19,21,24)$, $C_{54}(3,6,7,9,11,21,24,25)$, $C_{54}(3,5,6,9,13,21,23,24)$; 	

\item [\rm (269)]  $C_{54}(1,3,6,9,17,19,21,27)$, $C_{54}(3,6,7,9,11,21,25,27)$, $C_{54}(3,5,6,9,13,21,23,27)$; 	

\item [\rm (270)]  $C_{54}(1,3,6,9,17,19,24,27)$, $C_{54}(3,6,7,9,11,24,25,27)$, $C_{54}(3,5,6,9,13,23,24,27)$; 	

\item [\rm (271)]  $C_{54}(1,3,6,12,15,17,18,19)$, $C_{54}(3,6,7,11,12,15,18,25)$, $C_{54}(3,5,6,12,13,15,18,23)$; 	

\item [\rm (272)]  $C_{54}(1,3,6,12,15,17,19,21)$, $C_{54}(3,6,7,11,12,15,21,25)$, $C_{54}(3,5,6,12,13,15,21,23)$; 	

\item [\rm (273)]  $C_{54}(1,3,6,12,15,17,19,24)$, $C_{54}(3,6,7,11,12,15,24,25)$, $C_{54}(3,5,6,12,13,15,23,24)$; 	

\item [\rm (274)]  $C_{54}(1,3,6,12,15,17,19,27)$, $C_{54}(3,6,7,11,12,15,25,27)$, $C_{54}(3,5,6,12,13,15,23,27)$; 	

\item [\rm (275)]  $C_{54}(1,3,6,12,17,18,19,21)$, $C_{54}(3,6,7,11,12,18,21,25)$, $C_{54}(3,5,6,12,13,18,21,23)$; 	

\item [\rm (276)]  $C_{54}(1,3,6,12,17,18,19,24)$, $C_{54}(3,6,7,11,12,18,24,25)$, $C_{54}(3,5,6,12,13,18,23,24)$; 	

\item [\rm (277)]  $C_{54}(1,3,6,12,17,18,19,27)$, $C_{54}(3,6,7,11,12,18,25,27)$, $C_{54}(3,5,6,12,13,18,23,27)$; 	

\item [\rm (278)]  $C_{54}(1,3,6,12,17,19,21,24)$, $C_{54}(3,6,7,11,12,21,24,25)$, $C_{54}(3,5,6,12,13,21,23,24)$; 	

\item [\rm (279)]  $C_{54}(1,3,6,12,17,19,21,27)$, $C_{54}(3,6,7,11,12,21,25,27)$, $C_{54}(3,5,6,12,13,21,23,27)$; 	

\item [\rm (280)]  $C_{54}(1,3,6,12,17,19,24,27)$, $C_{54}(3,6,7,11,12,24,25,27)$, $C_{54}(3,5,6,12,13,23,24,27)$; 	

\item [\rm (281)]  $C_{54}(1,3,6,15,17,18,19,21)$, $C_{54}(3,6,7,11,15,18,21,25)$, $C_{54}(3,5,6,13,15,18,21,23)$; 	

\item [\rm (282)]  $C_{54}(1,3,6,15,17,18,19,24)$, $C_{54}(3,6,7,11,15,18,24,25)$, $C_{54}(3,5,6,13,15,18,23,24)$; 	

\item [\rm (283)]  $C_{54}(1,3,6,15,17,18,19,27)$, $C_{54}(3,6,7,11,15,18,25,27)$, $C_{54}(3,5,6,13,15,18,23,27)$; 	

\item [\rm (284)]  $C_{54}(1,3,6,15,17,19,21,24)$, $C_{54}(3,6,7,11,15,21,24,25)$, $C_{54}(3,5,6,13,15,21,23,24)$; 	

\item [\rm (285)]  $C_{54}(1,3,6,15,17,19,21,27)$, $C_{54}(3,6,7,11,15,21,25,27)$, $C_{54}(3,5,6,13,15,21,23,27)$; 	

\item [\rm (286)]  $C_{54}(1,3,6,15,17,19,24,27)$, $C_{54}(3,6,7,11,15,24,25,27)$, $C_{54}(3,5,6,13,15,23,24,27)$; 	

\item [\rm (287)]  $C_{54}(1,3,6,17,18,19,21,24)$, $C_{54}(3,6,7,11,18,21,24,25)$, $C_{54}(3,5,6,13,18,21,23,24)$; 	

\item [\rm (288)]  $C_{54}(1,3,6,17,18,19,21,27)$, $C_{54}(3,6,7,11,18,21,25,27)$, $C_{54}(3,5,6,13,18,21,23,27)$; 	

\item [\rm (289)]  $C_{54}(1,3,6,17,18,19,24,27)$, $C_{54}(3,6,7,11,18,24,25,27)$, $C_{54}(3,5,6,13,18,23,24,27)$; 	

\item [\rm (290)]  $C_{54}(1,3,6,17,19,21,24,27)$, $C_{54}(3,6,7,11,21,24,25,27)$, $C_{54}(3,5,6,13,21,23,24,27)$; 	

\item [\rm (291)]  $C_{54}(1,3,9,12,15,17,18,19)$, $C_{54}(3,7,9,11,12,15,18,25)$, $C_{54}(3,5,9,12,13,15,18,23)$; 	

\item [\rm (292)]  $C_{54}(1,3,9,12,15,17,19,21)$, $C_{54}(3,7,9,11,12,15,21,25)$, $C_{54}(3,5,9,12,13,15,21,23)$; 	

\item [\rm (293)]  $C_{54}(1,3,9,12,15,17,19,24)$, $C_{54}(3,7,9,11,12,15,24,25)$, $C_{54}(3,5,9,12,13,15,23,24)$; 	

\item [\rm (294)]  $C_{54}(1,3,9,12,15,17,19,27)$, $C_{54}(3,7,9,11,12,15,25,27)$, $C_{54}(3,5,9,12,13,15,23,27)$; 	

\item [\rm (295)]  $C_{54}(1,3,9,12,17,18,19,21)$, $C_{54}(3,7,9,11,12,18,21,25)$, $C_{54}(3,5,9,12,13,18,21,23)$; 	

\item [\rm (296)]  $C_{54}(1,3,9,12,17,18,19,24)$, $C_{54}(3,7,9,11,12,18,24,25)$, $C_{54}(3,5,9,12,13,18,23,24)$; 	

\item [\rm (297)]  $C_{54}(1,3,9,12,17,18,19,27)$, $C_{54}(3,7,9,11,12,18,25,27)$, $C_{54}(3,5,9,12,13,18,23,27)$; 	

\item [\rm (298)]  $C_{54}(1,3,9,12,17,19,21,24)$, $C_{54}(3,7,9,11,12,21,24,25)$, $C_{54}(3,5,9,12,13,21,23,24)$; 	

\item [\rm (299)]  $C_{54}(1,3,9,12,17,19,21,27)$, $C_{54}(3,7,9,11,12,21,25,27)$, $C_{54}(3,5,9,12,13,21,23,27)$; 	

\item [\rm (300)]  $C_{54}(1,3,9,12,17,19,24,27)$, $C_{54}(3,7,9,11,12,24,25,27)$, $C_{54}(3,5,9,12,13,23,24,27)$; 	

\item [\rm (301)]  $C_{54}(1,3,9,15,17,18,19,21)$, $C_{54}(3,7,9,11,15,18,21,25)$, $C_{54}(3,5,9,13,15,18,21,23)$; 	

\item [\rm (302)]  $C_{54}(1,3,9,15,17,18,19,24)$, $C_{54}(3,7,9,11,15,18,24,25)$, $C_{54}(3,5,9,13,15,18,23,24)$; 	

\item [\rm (303)]  $C_{54}(1,3,9,15,17,18,19,27)$, $C_{54}(3,7,9,11,15,18,25,27)$, $C_{54}(3,5,9,13,15,18,23,27)$; 	

\item [\rm (304)]  $C_{54}(1,3,9,15,17,19,21,24)$, $C_{54}(3,7,9,11,15,21,24,25)$, $C_{54}(3,5,9,13,15,21,23,24)$; 	

\item [\rm (305)]  $C_{54}(1,3,9,15,17,19,21,27)$, $C_{54}(3,7,9,11,15,21,25,27)$, $C_{54}(3,5,9,13,15,21,23,27)$; 	

\item [\rm (306)]  $C_{54}(1,3,9,15,17,19,24,27)$, $C_{54}(3,7,9,11,15,24,25,27)$, $C_{54}(3,5,9,13,15,23,24,27)$; 	

\item [\rm (307)]  $C_{54}(1,3,9,17,18,19,21,24)$, $C_{54}(3,7,9,11,18,21,24,25)$, $C_{54}(3,5,9,13,18,21,23,24)$; 	

\item [\rm (308)]  $C_{54}(1,3,9,17,18,19,21,27)$, $C_{54}(3,7,9,11,18,21,24,25)$, $C_{54}(3,5,9,13,18,21,23,24)$; 	

\item [\rm (309)]  $C_{54}(1,3,9,17,18,19,24,27)$, $C_{54}(3,7,9,11,18,24,25,27)$, $C_{54}(3,5,9,13,18,23,24,27)$; 	

\item [\rm (310)]  $C_{54}(1,3,9,17,19,21,24,27)$, $C_{54}(3,7,9,11,21,24,25,27)$, $C_{54}(3,5,9,13,21,23,24,27)$; 	

\item [\rm (311)]  $C_{54}(1,3,12,15,17,18,19,21)$, $C_{54}(3,7,11,12,15,18,21,25)$, $C_{54}(3,5,12,13,15,18,21,23)$; 	

\item [\rm (312)]  $C_{54}(1,3,12,15,17,18,19,24)$, $C_{54}(3,7,11,12,15,18,24,25)$, $C_{54}(3,5,12,13,15,18,23,24)$; 	

\item [\rm (313)]  $C_{54}(1,3,12,15,17,18,19,27)$, $C_{54}(3,7,11,12,15,18,25,27)$, $C_{54}(3,5,12,13,15,18,23,27)$; 	

\item [\rm (314)]  $C_{54}(1,3,12,15,17,19,21,24)$, $C_{54}(3,7,11,12,15,21,24,25)$, $C_{54}(3,5,12,13,15,21,23,24)$; 	

\item [\rm (315)]  $C_{54}(1,3,12,15,17,19,21,27)$, $C_{54}(3,7,11,12,15,21,25,27)$, $C_{54}(3,5,12,13,15,21,23,27)$; 	

\item [\rm (316)]  $C_{54}(1,3,12,15,17,19,24,27)$, $C_{54}(3,7,11,12,15,24,25,27)$, $C_{54}(3,5,12,13,15,23,24,27)$; 	

\item [\rm (317)]  $C_{54}(1,3,12,17,18,19,21,24)$, $C_{54}(3,7,11,12,18,21,24,25)$, $C_{54}(3,5,12,13,18,21,23,24)$; 	

\item [\rm (318)]  $C_{54}(1,3,12,17,18,19,21,27)$, $C_{54}(3,7,11,12,18,21,25,27)$, $C_{54}(3,5,12,13,18,21,23,27)$; 	

\item [\rm (319)]  $C_{54}(1,3,12,17,18,19,24,27)$, $C_{54}(3,7,11,12,18,24,25,27)$, $C_{54}(3,5,12,13,18,23,24,27)$; 	

\item [\rm (320)]  $C_{54}(1,3,12,17,19,21,24,27)$, $C_{54}(3,7,11,12,21,24,25,27)$, $C_{54}(3,5,12,13,21,23,24,27)$; 	

\item [\rm (321)]  $C_{54}(1,3,15,17,18,19,21,24)$, $C_{54}(3,7,11,15,18,21,24,25)$, $C_{54}(3,5,13,15,18,21,23,24)$; 	

\item [\rm (322)]  $C_{54}(1,3,15,17,18,19,21,27)$, $C_{54}(3,7,11,15,18,21,25,27)$, $C_{54}(3,5,13,15,18,21,23,27)$; 	

\item [\rm (323)]  $C_{54}(1,3,15,17,18,19,24,27)$, $C_{54}(3,7,11,15,18,24,25,27)$, $C_{54}(3,5,13,15,18,23,24,27)$; 	

\item [\rm (324)]  $C_{54}(1,3,15,17,19,21,24,27)$, $C_{54}(3,7,11,15,21,24,25,27)$, $C_{54}(3,5,13,15,21,23,24,27)$; 	

\item [\rm (325)]  $C_{54}(1,3,17,18,19,21,24,27)$, $C_{54}(3,7,11,18,21,24,25,27)$, $C_{54}(3,5,13,18,21,23,24,27)$; 	

\item [\rm (326)]  $C_{54}(1,6,9,12,15,17,18,19)$, $C_{54}(6,7,9,11,12,15,18,25)$, $C_{54}(5,6,9,12,13,15,18,23)$; 	

\item [\rm (327)]  $C_{54}(1,6,9,12,15,17,19,21)$, $C_{54}(6,7,9,11,12,15,21,25)$, $C_{54}(5,6,9,12,13,15,21,23)$; 	

\item [\rm (328)]  $C_{54}(1,6,9,12,15,17,19,24)$, $C_{54}(6,7,9,11,12,15,24,25)$, $C_{54}(5,6,9,12,13,15,23,24)$; 	

\item [\rm (329)]  $C_{54}(1,6,9,12,15,17,19,27)$, $C_{54}(6,7,9,11,12,15,25,27)$, $C_{54}(5,6,9,12,13,15,23,27)$; 	

\item [\rm (330)]  $C_{54}(1,6,9,12,18,17,19,21)$, $C_{54}(6,7,9,11,12,18,21,25)$, $C_{54}(5,6,9,12,13,18,21,23)$; 	

\item [\rm (331)]  $C_{54}(1,6,9,12,18,17,19,24)$, $C_{54}(6,7,9,11,12,18,24,25)$, $C_{54}(5,6,9,12,13,18,23,24)$; 	

\item [\rm (332)]  $C_{54}(1,6,9,12,18,17,19,27)$, $C_{54}(6,7,9,11,12,18,25,27)$, $C_{54}(5,6,9,12,13,18,23,27)$; 	

\item [\rm (333)]  $C_{54}(1,6,9,12,17,19,21,24)$, $C_{54}(6,7,9,11,12,21,24,25)$, $C_{54}(5,6,9,12,13,21,23,24)$; 	

\item [\rm (334)]  $C_{54}(1,6,9,12,17,19,21,27)$, $C_{54}(6,7,9,11,12,21,25,27)$, $C_{54}(5,6,9,12,13,21,23,27)$; 	

\item [\rm (335)]  $C_{54}(1,6,9,12,17,19,24,27)$, $C_{54}(6,7,9,11,12,24,25,27)$, $C_{54}(5,6,9,12,13,23,24,27)$; 	

\item [\rm (336)]  $C_{54}(1,6,9,15,17,18,19,21)$, $C_{54}(6,7,9,11,15,18,21,25)$, $C_{54}(5,6,9,13,15,18,21,23)$; 	

\item [\rm (337)]  $C_{54}(1,6,9,15,17,18,19,24)$, $C_{54}(6,7,9,11,15,18,24,25)$, $C_{54}(5,6,9,13,15,18,23,24)$; 	

\item [\rm (338)]  $C_{54}(1,6,9,15,17,18,19,27)$, $C_{54}(6,7,9,11,15,18,25,27)$, $C_{54}(5,6,9,13,15,18,23,27)$; 	

\item [\rm (339)]  $C_{54}(1,6,9,15,17,19,21,24)$, $C_{54}(6,7,9,11,15,21,24,25)$, $C_{54}(5,6,9,13,15,21,23,24)$; 	

\item [\rm (340)]  $C_{54}(1,6,9,15,17,19,21,27)$, $C_{54}(6,7,9,11,15,21,25,27)$, $C_{54}(5,6,9,13,15,21,23,27)$; 	

\item [\rm (341)]  $C_{54}(1,6,9,15,17,19,24,27)$, $C_{54}(6,7,9,11,15,24,25,27)$, $C_{54}(5,6,9,13,15,23,24,27)$; 	

\item [\rm (342)]  $C_{54}(1,6,9,17,18,19,21,24)$, $C_{54}(6,7,9,11,18,21,24,25)$, $C_{54}(5,6,9,13,18,21,23,24)$; 	

\item [\rm (343)]  $C_{54}(1,6,9,17,18,19,21,27)$, $C_{54}(6,7,9,11,18,21,25,27)$, $C_{54}(5,6,9,13,18,21,23,27)$; 	

\item [\rm (344)]  $C_{54}(1,6,9,17,18,19,24,27)$, $C_{54}(6,7,9,11,18,24,25,27)$, $C_{54}(5,6,9,13,18,23,24,27)$; 	

\item [\rm (345)]  $C_{54}(1,6,9,17,19,21,24,27)$, $C_{54}(6,7,9,11,21,24,25,27)$, $C_{54}(5,6,9,13,21,23,24,27)$; 	

\item [\rm (346)]  $C_{54}(1,6,12,15,17,18,19,21)$, $C_{54}(6,7,11,12,15,18,21,25)$, $C_{54}(5,6,12,13,15,18,21,23)$; 	

\item [\rm (347)]  $C_{54}(1,6,12,15,17,18,19,24)$, $C_{54}(6,7,11,12,15,18,24,25)$, $C_{54}(5,6,12,13,15,18,23,24)$; 	

\item [\rm (348)]  $C_{54}(1,6,12,15,17,18,19,27)$, $C_{54}(6,7,11,12,15,18,25,27)$, $C_{54}(5,6,12,13,15,18,23,27)$; 	

\item [\rm (349)]  $C_{54}(1,6,12,15,17,19,21,24)$, $C_{54}(6,7,11,12,15,21,24,25)$, $C_{54}(5,6,12,13,15,21,23,24)$; 	

\item [\rm (350)]  $C_{54}(1,6,12,15,17,19,21,27)$, $C_{54}(6,7,11,12,15,21,25,27)$, $C_{54}(5,6,12,13,15,21,23,27)$; 	

\item [\rm (351)]  $C_{54}(1,6,12,15,17,19,24,27)$, $C_{54}(6,7,11,12,15,24,25,27)$, $C_{54}(5,6,12,13,15,23,24,27)$; 	

\item [\rm (352)]  $C_{54}(1,6,12,17,18,19,21,24)$, $C_{54}(6,7,11,12,18,21,24,25)$, $C_{54}(5,6,12,13,18,21,23,24)$; 	

\item [\rm (353)]  $C_{54}(1,6,12,17,18,19,21,27)$, $C_{54}(6,7,11,12,18,21,25,27)$, $C_{54}(5,6,12,13,18,21,23,27)$; 	

\item [\rm (354)]  $C_{54}(1,6,12,17,18,19,24,27)$, $C_{54}(6,7,11,12,18,24,25,27)$, $C_{54}(5,6,12,13,18,23,24,27)$; 	

\item [\rm (355)]  $C_{54}(1,6,12,17,19,21,24,27)$, $C_{54}(6,7,11,12,21,24,25,27)$, $C_{54}(5,6,12,13,21,23,24,27)$; 	

\item [\rm (356)]  $C_{54}(1,6,15,17,18,19,21,24)$, $C_{54}(6,7,11,15,18,21,24,25)$, $C_{54}(5,6,13,15,18,21,23,24)$; 	

\item [\rm (357)]  $C_{54}(1,6,15,17,18,19,21,27)$, $C_{54}(6,7,11,15,18,21,25,27)$, $C_{54}(5,6,13,15,18,21,23,27)$; 	

\item [\rm (358)]  $C_{54}(1,6,15,17,18,19,24,27)$, $C_{54}(6,7,11,15,18,24,25,27)$, $C_{54}(5,6,13,15,18,23,24,27)$; 	

\item [\rm (359)]  $C_{54}(1,6,15,17,19,21,24,27)$, $C_{54}(6,7,11,15,21,24,25,27)$, $C_{54}(5,6,13,15,21,23,24,27)$; 	

\item [\rm (360)]  $C_{54}(1,6,17,18,19,21,24,27)$, $C_{54}(6,7,11,18,21,24,25,27)$, $C_{54}(5,6,13,18,21,23,24,27)$; 	

\item [\rm (361)]  $C_{54}(1,9,12,15,17,18,19,21)$, $C_{54}(7,9,11,12,15,18,21,25)$, $C_{54}(5,9,12,13,15,18,21,23)$; 	

\item [\rm (362)]  $C_{54}(1,9,12,15,17,18,19,24)$, $C_{54}(7,9,11,12,15,18,24,25)$, $C_{54}(5,9,12,13,15,18,23,24)$; 	

\item [\rm (363)]  $C_{54}(1,9,12,15,17,18,19,27)$, $C_{54}(7,9,11,12,15,18,25,27)$, $C_{54}(5,9,12,13,15,18,23,27)$; 	

\item [\rm (364)]  $C_{54}(1,9,12,15,17,19,21,24)$, $C_{54}(7,9,11,12,15,21,24,25)$, $C_{54}(5,9,12,13,15,21,23,24)$; 	

\item [\rm (365)]  $C_{54}(1,9,12,15,17,19,21,27)$, $C_{54}(7,9,11,12,15,21,25,27)$, $C_{54}(5,9,12,13,15,21,23,27)$; 	

\item [\rm (366)]  $C_{54}(1,9,12,15,17,19,24,27)$, $C_{54}(7,9,11,12,15,24,25,27)$, $C_{54}(5,9,12,13,15,23,24,27)$; 	

\item [\rm (367)]  $C_{54}(1,9,12,17,18,19,21,24)$, $C_{54}(7,9,11,12,18,21,24,25)$, $C_{54}(5,9,12,13,18,21,23,24)$; 	

\item [\rm (368)]  $C_{54}(1,9,12,17,18,19,21,27)$, $C_{54}(7,9,11,12,18,21,25,27)$, $C_{54}(5,9,12,13,18,21,23,27)$; 	

\item [\rm (369)]  $C_{54}(1,9,12,17,18,19,24,27)$, $C_{54}(7,9,11,12,18,24,25,27)$, $C_{54}(5,9,12,13,18,23,24,27)$; 	

\item [\rm (370)]  $C_{54}(1,9,12,17,19,21,24,27)$, $C_{54}(7,9,11,12,21,24,25,27)$, $C_{54}(5,9,12,13,21,23,24,27)$; 	

\item [\rm (371)]  $C_{54}(1,9,15,17,18,19,21,24)$, $C_{54}(7,9,11,15,18,21,24,25)$, $C_{54}(5,9,13,15,18,21,23,24)$; 	

\item [\rm (372)]  $C_{54}(1,9,15,17,18,19,21,27)$, $C_{54}(7,9,11,15,18,21,25,27)$, $C_{54}(5,9,13,15,18,21,23,27)$; 	

\item [\rm (373)]  $C_{54}(1,9,15,17,18,19,24,27)$, $C_{54}(7,9,11,15,18,24,25,27)$, $C_{54}(5,9,13,15,18,23,24,27)$; 	

\item [\rm (374)]  $C_{54}(1,9,15,17,19,21,24,27)$, $C_{54}(7,9,11,15,21,24,25,27)$, $C_{54}(5,9,13,15,21,23,24,27)$; 	

\item [\rm (375)]  $C_{54}(1,9,17,18,19,21,24,27)$, $C_{54}(7,9,11,18,21,24,25,27)$, $C_{54}(5,9,13,18,21,23,24,27)$; 	

\item [\rm (376)]  $C_{54}(1,12,15,17,18,19,21,24)$, $C_{54}(7,11,12,15,18,21,24,25)$, $C_{54}(5,12,13,15,18,21,23,24)$; 

\item [\rm (377)]  $C_{54}(1,12,15,17,18,19,21,27)$, $C_{54}(7,11,12,15,18,21,25,27)$, $C_{54}(5,12,13,15,18,21,23,27)$; 

\item [\rm (378)] $C_{54}(1,12,15,17,18,19,24,27)$, $C_{54}(7,11,12,15,18,24,25,27)$, $C_{54}(5,12,13,15,18,23,24,27)$; 

\item [\rm (379)]  $C_{54}(1,12,15,17,19,21,24,27)$, $C_{54}(7,11,12,15,21,24,25,27)$, $C_{54}(5,12,13,15,21,23,24,27)$; 

\item [\rm (380)]  $C_{54}(1,12,17,18,19,21,24,27)$, $C_{54}(7,11,12,18,21,24,25,27)$, $C_{54}(5,12,13,18,21,23,24,27)$; 

\item [\rm (381)] $C_{54}(1,15,17,18,19,21,24,27)$, $C_{54}(7,11,15,18,21,24,25,27)$, $C_{54}(5,13,15,18,21,23,24,27)$; 	

\item [\rm (382)]   $C_{54}(1,3,6,9,12,15,17,18,19)$, $C_{54}(3,6,7,9,11,12,15,18,25)$, $C_{54}(3,5,6,9,12,13,15,18,23)$; 	

\item [\rm (383)]   $C_{54}(1,3,6,9,12,15,17,19,21)$, $C_{54}(3,6,7,9,11,12,15,21,25)$, $C_{54}(3,5,6,9,12,13,15,21,23)$; 	

\item [\rm (384)]   $C_{54}(1,3,6,9,12,15,17,19,24)$, $C_{54}(3,6,7,9,11,12,15,24,25)$, $C_{54}(3,5,6,9,12,13,15,23,24)$; 	

\item [\rm (385)]   $C_{54}(1,3,6,9,12,15,17,19,27)$, $C_{54}(3,6,7,9,11,12,15,25,27)$, $C_{54}(3,5,6,9,12,13,15,23,27)$; 	

\item [\rm (386)]  $C_{54}(1,3,6,9,12,17,18,19,21)$, $C_{54}(3,6,7,9,11,12,18,21,25)$, $C_{54}(3,5,6,9,12,13,18,21,23)$; 	

\item [\rm (387)]  $C_{54}(1,3,6,9,12,17,18,19,24)$, $C_{54}(3,6,7,9,11,12,18,24,25)$, $C_{54}(3,5,6,9,12,13,18,23,24)$; 	

\item [\rm (388)]  $C_{54}(1,3,6,9,12,17,18,19,27)$, $C_{54}(3,6,7,9,11,12,18,25,27)$, $C_{54}(3,5,6,9,12,13,18,23,27)$; 	

\item [\rm (389)]  $C_{54}(1,3,6,9,12,17,19,21,24)$, $C_{54}(3,6,7,9,11,12,21,24,25)$, $C_{54}(3,5,6,9,12,13,21,23,24)$; 	

\item [\rm (390)]  $C_{54}(1,3,6,9,12,17,19,21,27)$, $C_{54}(3,6,7,9,11,12,21,25,27)$, $C_{54}(3,5,6,9,12,13,21,23,27)$; 	

\item [\rm (391)]  $C_{54}(1,3,6,9,12,17,19,24,27)$, $C_{54}(3,6,7,9,11,12,24,25,27)$, $C_{54}(3,5,6,9,12,13,23,24,27)$; 	

\item [\rm (392)]  $C_{54}(1,3,6,9,15,17,18,19,21)$, $C_{54}(3,6,7,9,11,15,18,21,25)$, $C_{54}(3,5,6,9,13,15,18,21,23)$; 	

\item [\rm (393)]  $C_{54}(1,3,6,9,15,17,18,19,24)$, $C_{54}(3,6,7,9,11,15,18,24,25)$, $C_{54}(3,5,6,9,13,15,18,23,24)$; 	

\item [\rm (394)]  $C_{54}(1,3,6,9,15,17,18,19,27)$, $C_{54}(3,6,7,9,11,15,18,25,27)$, $C_{54}(3,5,6,9,13,15,18,23,27)$; 	

\item [\rm (395)]  $C_{54}(1,3,6,9,15,17,19,21,24)$, $C_{54}(3,6,7,9,11,15,21,24,25)$, $C_{54}(3,5,6,9,13,15,21,23,24)$; 	

\item [\rm (396)]  $C_{54}(1,3,6,9,15,17,19,21,27)$, $C_{54}(3,6,7,9,11,15,21,25,27)$, $C_{54}(3,5,6,9,13,15,21,23,27)$; 	

\item [\rm (397)]  $C_{54}(1,3,6,9,15,17,19,24,27)$, $C_{54}(3,6,7,9,11,15,24,25,27)$, $C_{54}(3,5,6,9,13,15,23,24,27)$; 	

\item [\rm (398)]  $C_{54}(1,3,6,9,17,18,19,21,24)$, $C_{54}(3,6,7,9,11,18,21,24,25)$, $C_{54}(3,5,6,9,13,18,21,23,24)$; 	

\item [\rm (399)]  $C_{54}(1,3,6,9,17,18,19,21,27)$, $C_{54}(3,6,7,9,11,18,21,25,27)$, $C_{54}(3,5,6,9,13,18,21,23,27)$; 	

\item [\rm (400)]  $C_{54}(1,3,6,9,17,18,19,24,27)$, $C_{54}(3,6,7,9,11,18,24,25,27)$, $C_{54}(3,5,6,9,13,18,23,24,27)$; 	

\item [\rm (401)]  $C_{54}(1,3,6,9,17,19,21,24,27)$, $C_{54}(3,6,7,9,11,21,24,25,27)$, $C_{54}(3,5,6,9,13,21,23,24,27)$; 	

\item [\rm (402)]  $C_{54}(1,3,6,12,15,17,18,19,21)$, $C_{54}(3,6,7,11,12,15,18,21,25)$, $C_{54}(3,5,6,12,13,15,18,21,23)$; 	

\item [\rm (403)] $C_{54}(1,3,6,12,15,17,18,19,24)$, $C_{54}(3,6,7,11,12,15,18,24,25)$, $C_{54}(3,5,6,12,13,15,18,23,24)$; 	

\item [\rm (404)]  $C_{54}(1,3,6,12,15,17,18,19,27)$, $C_{54}(3,6,7,11,12,15,18,25,27)$, $C_{54}(3,5,6,12,13,15,18,23,27)$; 	

\item [\rm (405)]  $C_{54}(1,3,6,12,15,17,19,21,24)$, $C_{54}(3,6,7,11,12,15,21,24,25)$, $C_{54}(3,5,6,12,13,15,21,23,24)$; 	

\item [\rm (406)]  $C_{54}(1,3,6,12,15,17,19,21,27)$, $C_{54}(3,6,7,11,12,15,21,25,27)$, $C_{54}(3,5,6,12,13,15,21,23,27)$; 	

\item [\rm (407)]  $C_{54}(1,3,6,12,15,17,19,24,27)$, $C_{54}(3,6,7,11,12,15,24,25,27)$, $C_{54}(3,5,6,12,13,15,23,24,27)$; 	

\item [\rm (408)]  $C_{54}(1,3,6,12,17,18,19,21,24)$, $C_{54}(3,6,7,11,12,18,21,24,25)$, $C_{54}(3,5,6,12,13,18,21,23,24)$; 	

\item [\rm (409)]  $C_{54}(1,3,6,12,17,18,19,21,27)$, $C_{54}(3,6,7,11,12,18,21,25,27)$, $C_{54}(3,5,6,12,13,18,21,23,27)$; 	

\item [\rm (410)]  $C_{54}(1,3,6,12,17,18,19,24,27)$, $C_{54}(3,6,7,11,12,18,24,25,27)$, $C_{54}(3,5,6,12,13,18,23,24,27)$; 	

\item [\rm (411)]  $C_{54}(1,3,6,12,17,19,21,24,27)$, $C_{54}(3,6,7,11,12,21,24,25,27)$, $C_{54}(3,5,6,12,13,21,23,24,27)$; 	

\item [\rm (412)]  $C_{54}(1,3,6,15,17,18,19,21,24)$, $C_{54}(3,6,7,11,15,18,21,24,25)$, $C_{54}(3,5,6,13,15,18,21,23,24)$; 	

\item [\rm (413)]  $C_{54}(1,3,6,15,17,18,19,21,27)$, $C_{54}(3,6,7,11,15,18,21,25,27)$, $C_{54}(3,5,6,13,15,18,21,23,27)$; 	

\item [\rm (414)]  $C_{54}(1,3,6,15,17,18,19,24,27)$, $C_{54}(3,6,7,11,15,18,24,25,27)$, $C_{54}(3,5,6,13,15,18,23,24,27)$; 	

\item [\rm (415)]  $C_{54}(1,3,6,15,17,19,21,24,27)$, $C_{54}(3,6,7,11,15,21,24,25,27)$, $C_{54}(3,5,6,13,15,21,23,24,27)$; 	

\item [\rm (416)]  $C_{54}(1,3,6,17,18,19,21,24,27)$, $C_{54}(3,6,7,11,18,21,24,25,27)$, $C_{54}(3,5,6,13,18,21,23,24,27)$; 	

\item [\rm (417)]  $C_{54}(1,3,9,12,15,17,18,19,21)$, $C_{54}(3,7,9,11,12,15,18,21,25)$, $C_{54}(3,5,9,12,13,15,18,21,23)$; 	

\item [\rm (418)]  $C_{54}(1,3,9,12,15,17,18,19,24)$, $C_{54}(3,7,9,11,12,15,18,24,25)$, $C_{54}(3,5,9,12,13,15,18,23,24)$; 	

\item [\rm (419)]  $C_{54}(1,3,9,12,15,17,18,19,27)$, $C_{54}(3,7,9,11,12,15,18,25,27)$, $C_{54}(3,5,9,12,13,15,18,23,27)$; 	

\item [\rm (420)]  $C_{54}(1,3,9,12,15,17,19,21,24)$, $C_{54}(3,7,9,11,12,15,21,24,25)$, $C_{54}(3,5,9,12,13,15,21,23,24)$; 	

\item [\rm (421)]  $C_{54}(1,3,9,12,15,17,19,21,27)$, $C_{54}(3,7,9,11,12,15,21,25,27)$, $C_{54}(3,5,9,12,13,15,21,23,27)$; 	

\item [\rm (422)]  $C_{54}(1,3,9,12,15,17,19,24,27)$, $C_{54}(3,7,9,11,12,15,24,25,27)$, $C_{54}(3,5,9,12,13,15,23,24,27)$; 	

\item [\rm (423)]  $C_{54}(1,3,9,12,17,18,19,21,24)$, $C_{54}(3,7,9,11,12,18,21,24,25)$, $C_{54}(3,5,9,12,13,18,21,23,24)$; 	

\item [\rm (424)]  $C_{54}(1,3,9,12,17,18,19,21,27)$, $C_{54}(3,7,9,11,12,18,21,25,27)$, $C_{54}(3,5,9,12,13,18,21,23,27)$; 	

\item [\rm (425)]  $C_{54}(1,3,9,12,17,18,19,24,27)$, $C_{54}(3,7,9,11,12,18,24,25,27)$, $C_{54}(3,5,9,12,13,18,23,24,27)$; 	

\item [\rm (426)]  $C_{54}(1,3,9,12,17,19,21,24,27)$, $C_{54}(3,7,9,11,12,21,24,25,27)$, $C_{54}(3,5,9,12,13,21,23,24,27)$; 	

\item [\rm (427)]  $C_{54}(1,3,9,15,17,18,19,21,24)$, $C_{54}(3,7,9,11,15,18,21,24,25)$, $C_{54}(3,5,9,13,15,18,21,23,24)$; 	

\item [\rm (428)]  $C_{54}(1,3,9,15,17,18,19,21,27)$, $C_{54}(3,7,9,11,15,18,21,25,27)$, $C_{54}(3,5,9,13,15,18,21,23,27)$; 	

\item [\rm (429)]  $C_{54}(1,3,9,15,17,18,19,24,27)$, $C_{54}(3,7,9,11,15,18,24,25,27)$, $C_{54}(3,5,9,13,15,18,23,24,27)$; 	

\item [\rm (430)]  $C_{54}(1,3,9,15,17,19,21,24,27)$, $C_{54}(3,7,9,11,15,21,24,25,27)$, $C_{54}(3,5,9,13,15,21,23,24,27)$; 	

\item [\rm (431)]  $C_{54}(1,3,9,17,18,19,21,24,27)$, $C_{54}(3,7,9,11,18,21,24,25,27)$, $C_{54}(3,5,9,13,18,21,23,24,27)$; 	

\item [\rm (432)]  $C_{54}(1,3,12,15,17,18,19,21,24)$, $C_{54}(3,7,11,12,15,18,21,24,25)$, 

\hfill $C_{54}(3,5,12,13,15,18,21,23,24)$; 	

\item [\rm (433)]  $C_{54}(1,3,12,15,17,18,19,21,27)$, $C_{54}(3,7,11,12,15,18,21,25,27)$, 

\hfill $C_{54}(3,5,12,13,15,18,21,23,27)$; 	

\item [\rm (434)]  $C_{54}(1,3,12,15,17,18,19,24,27)$, $C_{54}(3,7,11,12,15,18,24,25,27)$, 

\hfill $C_{54}(3,5,12,13,15,18,23,27,27)$; 	

\item [\rm (435)]  $C_{54}(1,3,12,15,17,19,21,24,27)$, $C_{54}(3,7,11,12,15,21,24,25,27)$, 

\hfill $C_{54}(3,5,12,13,15,21,23,27,27)$; 	

\item [\rm (436)]  $C_{54}(1,3,12,17,18,19,21,24,27)$, $C_{54}(3,7,11,12,18,21,24,25,27)$, 

\hfill $C_{54}(3,5,12,13,18,21,23,27,27)$; 	

\item [\rm (437)]  $C_{54}(1,3,15,17,18,19,21,24,27)$, $C_{54}(3,7,11,15,18,21,24,25,27)$, 

\hfill $C_{54}(3,5,13,15,18,21,23,27,27)$; 	

\item [\rm (438)]  $C_{54}(1,6,9,12,15,17,18,19,21)$, $C_{54}(6,7,9,11,12,15,18,21,25)$, 

\hfill $C_{54}(5,6,9,12,13,15,18,21,23)$; 	

\item [\rm (439)]  $C_{54}(1,6,9,12,15,17,18,19,24)$, $C_{54}(6,7,9,11,12,15,18,24,25)$, 

\hfill $C_{54}(5,6,9,12,13,15,18,23,24)$; 	

\item [\rm (440)]  $C_{54}(1,6,9,12,15,17,18,19,27)$, $C_{54}(6,7,9,11,12,15,18,25,27)$, 

\hfill $C_{54}(5,6,9,12,13,15,18,23,27)$; 	

\item [\rm (441)]  $C_{54}(1,6,9,12,15,17,19,21,24)$, $C_{54}(6,7,9,11,12,15,21,24,25)$, 

\hfill $C_{54}(5,6,9,12,13,15,21,23,24)$; 	

\item [\rm (442)]  $C_{54}(1,6,9,12,15,17,19,21,27)$, $C_{54}(6,7,9,11,12,15,21,25,27)$, 

\hfill $C_{54}(5,6,9,12,13,15,21,23,27)$; 	

\item [\rm (443)]  $C_{54}(1,6,9,12,15,17,19,24,27)$, $C_{54}(6,7,9,11,12,15,24,25,27)$, 

\hfill $C_{54}(5,6,9,12,13,15,23,24,27)$; 	

\item [\rm (444)]  $C_{54}(1,6,9,12,17,18,19,21,24)$, $C_{54}(6,7,9,11,12,18,21,24,25)$, 

\hfill $C_{54}(5,6,9,12,13,18,21,23,24)$; 	

\item [\rm (445)]  $C_{54}(1,6,9,12,17,18,19,21,27)$, $C_{54}(6,7,9,11,12,18,21,25,27)$, 

\hfill $C_{54}(5,6,9,12,13,18,21,23,27)$; 	

\item [\rm (446)]  $C_{54}(1,6,9,12,17,18,19,24,27)$, $C_{54}(6,7,9,11,12,18,24,25,27)$, 

\hfill $C_{54}(5,6,9,12,13,18,23,24,27)$; 	

\item [\rm (447)]  $C_{54}(1,6,9,12,17,19,21,24,27)$, $C_{54}(6,7,9,11,12,21,24,25,27)$, 

\hfill $C_{54}(5,6,9,12,13,21,23,24,27)$; 	

\item [\rm (448)]  $C_{54}(1,6,9,15,17,18,19,21,24)$, $C_{54}(6,7,9,11,15,18,21,24,25)$, 

\hfill $C_{54}(5,6,9,13,15,18,21,23,24)$; 	

\item [\rm (449)]  $C_{54}(1,6,9,15,17,18,19,21,27)$, $C_{54}(6,7,9,11,15,18,21,25,27)$, 

\hfill $C_{54}(5,6,9,13,15,18,21,23,27)$; 	

\item [\rm (450)]  $C_{54}(1,6,9,15,17,18,19,24,27)$, $C_{54}(6,7,9,11,15,18,24,25,27)$, 

\hfill $C_{54}(5,6,9,13,15,18,23,24,27)$; 	

\item [\rm (451)]  $C_{54}(1,6,9,15,17,19,21,24,27)$, $C_{54}(6,7,9,11,15,21,24,25,27)$, 

\hfill $C_{54}(5,6,9,13,15,21,23,24,27)$; 	

\item [\rm (452)]  $C_{54}(1,6,9,17,18,19,21,24,27)$, $C_{54}(6,7,9,11,18,21,24,25,27)$, 

\hfill $C_{54}(5,6,9,13,18,21,23,24,27)$; 	

\item [\rm (453)]  $C_{54}(1,6,12,15,17,18,19,21,24)$, $C_{54}(6,7,11,12,15,18,21,24,25)$, 

\hfill $C_{54}(5,6,12,13,15,18,21,23,24)$; 	

\item [\rm (454)]  $C_{54}(1,6,12,15,17,18,19,21,27)$, $C_{54}(6,7,11,12,15,18,21,25,27)$, 

\hfill $C_{54}(5,6,12,13,15,18,21,23,27)$; 	

\item [\rm (455)]  $C_{54}(1,6,12,15,17,18,19,24,27)$, $C_{54}(6,7,11,12,15,18,24,25,27)$, 

\hfill $C_{54}(5,6,12,13,15,18,23,24,27)$; 	

\item [\rm (456)]  $C_{54}(1,6,12,15,17,19,21,24,27)$, $C_{54}(6,7,11,12,15,21,24,25,27)$, 

\hfill $C_{54}(5,6,12,13,15,21,23,24,27)$; 	

\item [\rm (457)]  $C_{54}(1,6,12,17,18,19,21,24,27)$, $C_{54}(6,7,11,12,18,21,24,25,27)$, 

\hfill $C_{54}(5,6,12,13,18,21,23,24,27)$; 	

\item [\rm (458)]  $C_{54}(1,6,15,17,18,19,21,24,27)$, $C_{54}(6,7,11,15,18,21,24,25,27)$, 

\hfill $C_{54}(5,6,13,15,18,21,23,24,27)$; 	

\item [\rm (459)]  $C_{54}(1,9,12,15,17,18,19,21,24)$, $C_{54}(7,9,11,12,15,18,21,24,25)$, 

\hfill $C_{54}(5,9,12,13,15,18,21,23,24)$; 	

\item [\rm (460)]  $C_{54}(1,9,12,15,17,18,19,21,27)$, $C_{54}(7,9,11,12,15,18,21,25,27)$, 

\hfill $C_{54}(5,9,12,13,15,18,21,23,27)$; 	

\item [\rm (461)]  $C_{54}(1,9,12,15,17,18,19,24,27)$, $C_{54}(7,9,11,12,15,18,24,25,27)$, 

\hfill $C_{54}(5,9,12,13,15,18,23,24,27)$; 	

\item [\rm (462)]  $C_{54}(1,9,12,15,17,19,21,24,27)$, $C_{54}(7,9,11,12,15,21,24,25,27)$, 

\hfill $C_{54}(5,9,12,13,15,21,23,24,27)$; 	

\item [\rm (463)]  $C_{54}(1,9,12,17,18,19,21,24,27)$, $C_{54}(7,9,11,12,18,21,24,25,27)$, 

\hfill $C_{54}(5,9,12,13,18,21,23,24,27)$; 	

\item [\rm (464)]  $C_{54}(1,9,15,17,18,19,21,24,27)$, $C_{54}(7,9,11,15,18,21,24,25,27)$, 

\hfill $C_{54}(5,9,13,15,18,21,23,24,27)$; 	

\item [\rm (465)]  $C_{54}(1,12,15,17,18,19,21,24,27)$, $C_{54}(7,11,12,15,18,21,24,25,27)$, 

\hfill $C_{54}(5,12,13,15,18,21,23,24,27)$; 	

\item [\rm (466)]  $C_{54}(1,3,6,9,12,15,17,18,19,21)$, $C_{54}(3,6,7,9,11,12,15,18,21,25)$, 

\hfill $C_{54}(3,5,6,9,12,13,15,18,21,23)$; 	

\item [\rm (467)]  $C_{54}(1,3,6,9,12,15,17,18,19,24)$, $C_{54}(3,6,7,9,11,12,15,18,24,25)$, 

\hfill $C_{54}(3,5,6,9,12,13,15,18,23,24)$; 	

\item [\rm (468)]  $C_{54}(1,3,6,9,12,15,17,18,19,27)$, $C_{54}(3,6,7,9,11,12,15,18,25,27)$, 

\hfill $C_{54}(3,5,6,9,12,13,15,18,23,27)$; 	

\item [\rm (469)]   $C_{54}(1,3,6,9,12,15,17,19,21,24)$, $C_{54}(3,6,7,9,11,12,15,21,24,25)$, 

\hfill $C_{54}(3,5,6,9,12,13,15,21,23,24)$; 	

\item [\rm (470)]   $C_{54}(1,3,6,9,12,15,17,19,21,27)$, $C_{54}(3,6,7,9,11,12,15,21,25,27)$, 

\hfill $C_{54}(3,5,6,9,12,13,15,21,23,27)$; 	

\item [\rm (471)]   $C_{54}(1,3,6,9,12,15,17,19,24,27)$, $C_{54}(3,6,7,9,11,12,15,24,25,27)$, 

\hfill $C_{54}(3,5,6,9,12,13,15,23,24,27)$; 	

\item [\rm (472)]   $C_{54}(1,3,6,9,12,17,18,19,21,24)$, $C_{54}(3,6,7,9,11,12,18,21,24,25)$, 

\hfill $C_{54}(3,5,6,9,12,13,18,21,23,24)$; 	

\item [\rm (473)]   $C_{54}(1,3,6,9,12,17,18,19,21,27)$, $C_{54}(3,6,7,9,11,12,18,21,25,27)$, 

\hfill $C_{54}(3,5,6,9,12,13,18,21,23,27)$; 	

\item [\rm (474)]  $C_{54}(1,3,6,9,12,17,18,19,24,27)$, $C_{54}(3,6,7,9,11,12,18,24,25,27)$, 

\hfill $C_{54}(3,5,6,9,12,13,18,23,24,27)$; 	

\item [\rm (475)]  $C_{54}(1,3,6,9,12,17,19,21,24,27)$, $C_{54}(3,6,7,9,11,12,21,24,25,27)$, 

\hfill $C_{54}(3,5,6,9,12,13,21,23,24,27)$; 	

\item [\rm (476)]  $C_{54}(1,3,6,9,15,17,18,19,21,24)$, $C_{54}(3,6,7,9,11,15,18,21,24,25)$, 

\hfill $C_{54}(3,5,6,9,13,15,18,21,23,24)$; 	

\item [\rm (477)]  $C_{54}(1,3,6,9,15,17,18,19,21,27)$, $C_{54}(3,6,7,9,11,15,18,21,25,27)$, 

\hfill $C_{54}(3,5,6,9,13,15,18,21,23,27)$; 	

\item [\rm (478)]  $C_{54}(1,3,6,9,15,17,18,19,24,27)$, $C_{54}(3,6,7,9,11,15,18,24,25,27)$, 

\hfill $C_{54}(3,5,6,9,13,15,18,23,24,27)$; 	

\item [\rm (479)]  $C_{54}(1,3,6,9,15,17,19,21,24,27)$, $C_{54}(3,6,7,9,11,15,21,24,25,27)$, 

\hfill $C_{54}(3,5,6,9,13,15,21,23,24,27)$; 	

\item [\rm (480)]  $C_{54}(1,3,6,9,17,18,19,21,24,27)$, $C_{54}(3,6,7,9,11,18,21,24,25,27)$, 

\hfill $C_{54}(3,5,6,9,13,18,21,23,24,27)$; 	

\item [\rm (481)]  $C_{54}(1,3,6,12,15,17,18,19,21,24)$, $C_{54}(3,6,7,11,12,15,18,21,24,25)$, 

\hfill $C_{54}(3,5,6,12,13,15,18,21,23,24)$; 	

\item [\rm (482)]  $C_{54}(1,3,6,12,15,17,18,19,21,27)$, $C_{54}(3,6,7,11,12,15,18,21,25,27)$, 

\hfill $C_{54}(3,5,6,12,13,15,18,21,23,27)$; 	

\item [\rm (483)] $C_{54}(1,3,6,12,15,17,18,19,24,27)$, $C_{54}(3,6,7,11,12,15,18,24,25,27)$, 

\hfill $C_{54}(3,5,6,12,13,15,18,23,24,27)$; 	

\item [\rm (484)] $C_{54}(1,3,6,12,15,17,19,21,24,27)$, $C_{54}(3,6,7,11,12,15,21,24,25,27)$, 

\hfill $C_{54}(3,5,6,12,13,15,21,23,24,27)$; 	

\item [\rm (485)]  $C_{54}(1,3,6,12,17,18,19,21,24,27)$, $C_{54}(3,6,7,11,12,18,21,24,25,27)$, 

\hfill $C_{54}(3,5,6,12,13,18,21,23,24,27)$; 	

\item [\rm (486)]  $C_{54}(1,3,6,15,17,18,19,21,24,27)$, $C_{54}(3,6,7,11,15,18,21,24,25,27)$, 

\hfill $C_{54}(3,5,6,13,15,18,21,23,24,27)$; 	

\item [\rm (487)]  $C_{54}(1,3,9,12,15,17,18,19,21,24)$, $C_{54}(3,7,9,11,12,15,18,21,24,25)$, 

\hfill $C_{54}(3,5,9,12,13,15,18,21,23,24)$; 	

\item [\rm (488)]  $C_{54}(1,3,9,12,15,17,18,19,21,27)$, $C_{54}(3,7,9,11,12,15,18,21,25,27)$, 

\hfill $C_{54}(3,5,9,12,13,15,18,21,23,27)$; 	

\item [\rm (489)]  $C_{54}(1,3,9,12,15,17,18,19,24,27)$, $C_{54}(3,7,9,11,12,15,18,24,25,27)$, 

\hfill $C_{54}(3,5,9,12,13,15,18,23,24,27)$; 	

\item [\rm (490)]  $C_{54}(1,3,9,12,15,17,19,21,24,27)$, $C_{54}(3,7,9,11,12,15,21,24,25,27)$, 

\hfill $C_{54}(3,5,9,12,13,15,21,23,24,27)$; 	

\item [\rm (491)]  $C_{54}(1,3,9,12,17,18,19,21,24,27)$, $C_{54}(3,7,9,11,12,18,21,24,25,27)$, 

\hfill $C_{54}(3,5,9,12,13,18,21,23,24,27)$; 	

\item [\rm (492)]  $C_{54}(1,3,9,15,17,18,19,21,24,27)$, $C_{54}(3,7,9,11,15,18,21,24,25,27)$, 

\hfill $C_{54}(3,5,9,13,15,18,21,23,24,27)$; 	

\item [\rm (493)]  $C_{54}(1,3,12,15,17,18,19,21,24,27)$, $C_{54}(3,7,11,12,15,18,21,24,25,27)$, 

\hfill $C_{54}(3,5,12,13,15,18,21,23,24,27)$; 	

\item [\rm (494)]  $C_{54}(1,6,9,12,15,17,18,19,21,24)$, $C_{54}(6,7,9,11,12,15,18,21,24,25)$, 

\hfill $C_{54}(5,6,9,12,13,15,18,21,23,24)$; 	

\item [\rm (495)]  $C_{54}(1,6,9,12,15,17,18,19,21,27)$, $C_{54}(6,7,9,11,12,15,18,21,25,27)$, 

\hfill $C_{54}(5,6,9,12,13,15,18,21,23,27)$; 	

\item [\rm (496)]  $C_{54}(1,6,9,12,15,17,18,19,24,27)$, $C_{54}(6,7,9,11,12,15,18,24,25,27)$, 

\hfill $C_{54}(5,6,9,12,13,15,18,23,24,27)$; 	

\item [\rm (497)]  $C_{54}(1,6,9,12,15,17,19,21,24,27)$, $C_{54}(6,7,9,11,12,15,21,24,25,27)$, 

\hfill $C_{54}(5,6,9,12,13,15,21,23,24,27)$; 	

\item [\rm (498)]  $C_{54}(1,6,9,12,17,18,19,21,24,27)$, $C_{54}(6,7,9,11,12,18,21,24,25,27)$, 

\hfill $C_{54}(5,6,9,12,13,18,21,23,24,27)$; 	

\item [\rm (499)]  $C_{54}(1,6,9,15,17,18,19,21,24,27)$, $C_{54}(6,7,9,11,15,18,21,24,25,27)$, 

\hfill $C_{54}(5,6,9,13,15,18,21,23,24,27)$; 	

\item [\rm (500)]  $C_{54}(1,6,12,15,17,18,19,21,24,27)$, $C_{54}(6,7,11,12,15,18,21,24,25,27)$, 

\hfill $C_{54}(5,6,12,13,15,18,21,23,24,27)$; 	

\item [\rm (501)]  $C_{54}(1,9,12,15,17,18,19,21,24,27)$, $C_{54}(7,9,11,12,15,18,21,24,25,27)$, 

\hfill $C_{54}(5,9,12,13,15,18,21,23,24,27)$; 	

\item [\rm (502)]  $C_{54}(1,3,6,9,12,15,17,18,19,21,24)$, $C_{54}(3,6,7,9,11,12,15,18,21,24,25)$, 

\hfill $C_{54}(3,5,6,9,12,13,15,18,21,23,24)$; 	

\item [\rm (503)]  $C_{54}(1,3,6,9,12,15,17,18,19,21,27)$, $C_{54}(3,6,7,9,11,12,15,18,21,25,27)$, 

\hfill $C_{54}(3,5,6,9,12,13,15,18,21,23,27)$; 	

\item [\rm (504)]  $C_{54}(1,3,6,9,12,15,17,18,19,24,27)$, $C_{54}(3,6,7,9,11,12,15,18,24,25,27)$, 

\hfill $C_{54}(3,5,6,9,12,13,15,18,23,24,27)$; 	

\item [\rm (505)]  $C_{54}(1,3,6,9,12,15,17,19,21,24,27)$, $C_{54}(3,6,7,9,11,12,15,21,24,25,27)$, 

\hfill $C_{54}(3,5,6,9,12,13,15,21,23,24,27)$; 	

\item [\rm (506)]   $C_{54}(1,3,6,9,12,17,18,19,21,24,27)$, $C_{54}(3,6,7,9,11,12,18,21,24,25,27)$, 

\hfill $C_{54}(3,5,6,9,12,13,18,21,23,24,27)$; 	

\item [\rm (507)]   $C_{54}(1,3,6,9,15,17,18,19,21,24,27)$, $C_{54}(3,6,7,9,11,15,18,21,24,25,27)$, 

\hfill $C_{54}(3,5,6,9,13,15,18,21,23,24,27)$; 	

\item [\rm (508)]  $C_{54}(1,3,6,12,15,17,18,19,21,24,27)$, $C_{54}(3,6,7,11,12,15,18,21,24,25,27)$, 

\hfill $C_{54}(3,5,6,12,13,15,18,21,23,24,27)$; 	

\item [\rm (509)]  $C_{54}(1,3,9,12,15,17,18,19,21,24,27)$, $C_{54}(3,7,9,11,12,15,18,21,24,25,27)$, 

\hfill $C_{54}(3,5,9,12,13,15,18,21,23,24,27)$; 	

\item [\rm (510)]  $C_{54}(1,6,9,12,15,17,18,19,21,24,27)$, $C_{54}(6,7,9,11,12,15,18,21,24,25,27)$, 

\hfill $C_{54}(5,6,9,12,13,15,18,21,23,24,27)$; 	

\item [\rm (511)]  $C_{54}(1,3,6,9,12,15,17,18,19,21,24,27)$, $C_{54}(3,6,7,9,11,12,15,18,21,24,25,27)$, 

\hfill $C_{54}(3,5,6,9,12,13,15,18,21,23,24,27)$.
\end{enumerate}

\noindent
{\bf   {\footnotesize (b) Triples of isomorphic circulant graphs of order 54 with jump sizes 2,16,20, 4,14,22, 8,10,26.}}
\begin{enumerate} 	
\item [\rm (1)]  $C_{54}(2,3,16,20)$, $C_{54}(3,4,14,22)$, $C_{54}(3,8,10,26)$; 	

\item [\rm (2)]  $C_{54}(2,6,16,20)$, $C_{54}(4,6,14,22)$, $C_{54}(6,8,10,26)$; 	

\item [\rm (3)]  $C_{54}(2,9,16,20)$, $C_{54}(4,9,14,22)$, $C_{54}(8,9,10,26)$; 	

\item [\rm (4)]  $C_{54}(2,12,16,20)$, $C_{54}(4,12,14,22)$, $C_{54}(8,10,12,26)$; 	

\item [\rm (5)]  $C_{54}(2,15,16,20)$, $C_{54}(4,14,15,22)$, $C_{54}(8,10,15,26)$; 	

\item [\rm (6)]  $C_{54}(2,16,18,20)$, $C_{54}(4,14,18,22)$, $C_{54}(8,10,18,26)$; 	

\item [\rm (7)]  $C_{54}(2,16,20,21)$, $C_{54}(4,14,21,22)$, $C_{54}(8,10,21,26)$; 	

\item [\rm (8)]  $C_{54}(2,16,20,24)$, $C_{54}(4,14,22,24)$, $C_{54}(8,10,24,26)$; 

\item [\rm (9)]  $C_{54}(2,16,20,27)$, $C_{54}(4,14,22,27)$, $C_{54}(8,10,26,27)$; 

\item [\rm (10)]  $C_{54}(2,3,6,16,20)$, $C_{54}(3,4,6,14,22)$, $C_{54}(3,6,8,10,26)$; 	

\item [\rm (11)]  $C_{54}(2,3,9,16,20)$, $C_{54}(3,4,9,14,22)$, $C_{54}(3,8,9,10,26)$; 	

\item [\rm (12)]  $C_{54}(2,3,12,16,20)$, $C_{54}(3,4,12,14,22)$, $C_{54}(3,8,10,12,26)$; 	

\item [\rm (13)]  $C_{54}(2,3,15,16,20)$, $C_{54}(3,4,14,15,22)$, $C_{54}(3,8,10,15,26)$; 	

\item [\rm (14)]  $C_{54}(2,3,16,18,20)$, $C_{54}(3,4,14,18,22)$, $C_{54}(3,8,10,18,26)$; 	

\item [\rm (15)]  $C_{54}(2,3,16,20,21)$, $C_{54}(3,4,14,21,22)$, $C_{54}(3,8,10,21,26)$; 	

\item [\rm (16)]  $C_{54}(2,3,16,20,24)$, $C_{54}(3,4,14,22,24)$, $C_{54}(3,8,10,24,26)$; 	

\item [\rm (17)]  $C_{54}(2,3,16,20,27)$, $C_{54}(3,4,14,22,27)$, $C_{54}(3,8,10,26,27)$; 	

\item [\rm (18)]  $C_{54}(2,6,9,16,20)$, $C_{54}(4,6,9,14,22)$, $C_{54}(6,8,9,10,26)$; 	

\item [\rm (19)]  $C_{54}(2,6,12,16,20)$, $C_{54}(4,6,12,14,22)$, $C_{54}(6,8,10,12,26)$; 	

\item [\rm (20)]  $C_{54}(2,6,15,16,20)$, $C_{54}(4,6,14,15,22)$, $C_{54}(6,8,10,15,26)$; 	

\item [\rm (21)]  $C_{54}(2,6,16,18,20)$, $C_{54}(4,6,14,18,22)$, $C_{54}(6,8,10,18,26)$; 	

\item [\rm (22)]  $C_{54}(2,6,16,20,21)$, $C_{54}(4,6,14,21,22)$, $C_{54}(6,8,10,21,26)$; 	

\item [\rm (23)]  $C_{54}(2,6,16,20,24)$, $C_{54}(4,6,14,22,24)$, $C_{54}(6,8,10,24,26)$; 	

\item [\rm (24)]  $C_{54}(2,6,16,20,27)$, $C_{54}(4,6,14,22,27)$, $C_{54}(6,8,10,26,27)$; 	

\item [\rm (25)]  $C_{54}(2,9,12,16,20)$, $C_{54}(4,9,12,14,22)$, $C_{54}(8,9,10,12,26)$; 	

\item [\rm (26)]  $C_{54}(2,9,15,16,20)$, $C_{54}(4,9,14,15,22)$, $C_{54}(8,9,10,15,26)$; 	

\item [\rm (27)]  $C_{54}(2,9,16,18,20)$, $C_{54}(4,9,14,18,22)$, $C_{54}(8,9,10,18,26)$; 	

\item [\rm (28)]  $C_{54}(2,9,16,20,21)$, $C_{54}(4,9,14,21,22)$, $C_{54}(8,9,10,21,26)$; 	

\item [\rm (29)]  $C_{54}(2,9,16,20,24)$, $C_{54}(4,9,14,22,24)$, $C_{54}(8,9,10,24,26)$; 	

\item [\rm (30)]  $C_{54}(2,9,16,20,27)$, $C_{54}(4,9,14,22,27)$, $C_{54}(8,9,10,26,27)$; 	

\item [\rm (31)]  $C_{54}(2,12,15,16,20)$, $C_{54}(4,12,14,15,22)$, $C_{54}(8,10,12,15,26)$; 	

\item [\rm (32)]  $C_{54}(2,12,16,18,20)$, $C_{54}(4,12,14,18,22)$, $C_{54}(8,10,12,18,26)$; 	

\item [\rm (33)]  $C_{54}(2,12,16,20,21)$, $C_{54}(4,12,14,21,22)$, $C_{54}(8,10,12,21,26)$; 	

\item [\rm (34)]  $C_{54}(2,12,16,20,24)$, $C_{54}(4,12,14,22,24)$, $C_{54}(8,10,12,24,26)$; 	

\item [\rm (35)]  $C_{54}(2,12,16,20,27)$, $C_{54}(4,12,14,22,27)$, $C_{54}(8,10,12,26,27)$; 	

\item [\rm (36)]  $C_{54}(2,15,16,18,20)$, $C_{54}(4,14,15,18,22)$, $C_{54}(8,10,15,18,26)$; 	

\item [\rm (37)]  $C_{54}(2,15,16,20,21)$, $C_{54}(4,14,15,21,22)$, $C_{54}(8,10,15,21,26)$; 	

\item [\rm (38)]  $C_{54}(2,15,16,20,24)$, $C_{54}(4,14,15,22,24)$, $C_{54}(8,10,15,24,26)$; 	

\item [\rm (39)]  $C_{54}(2,15,16,20,27)$, $C_{54}(4,14,15,22,27)$, $C_{54}(8,10,15,26,27)$; 	

\item [\rm (40)]  $C_{54}(2,18,16,20,21)$, $C_{54}(4,14,18,21,22)$, $C_{54}(8,10,18,21,26)$; 	

\item [\rm (41)]  $C_{54}(2,18,16,20,24)$, $C_{54}(4,14,18,22,24)$, $C_{54}(8,10,18,24,26)$; 	

\item [\rm (42)]  $C_{54}(2,18,16,20,27)$, $C_{54}(4,14,18,22,27)$, $C_{54}(8,10,18,26,27)$; 	

\item [\rm (43)]  $C_{54}(2,16,20,21,24)$, $C_{54}(4,14,21,22,24)$, $C_{54}(8,10,21,24,26)$; 	

\item [\rm (44)]  $C_{54}(2,16,20,21,27)$, $C_{54}(4,14,21,22,27)$, $C_{54}(8,10,21,26,27)$; 	

\item [\rm (45)]  $C_{54}(2,16,20,24,27)$, $C_{54}(4,14,22,24,27)$, $C_{54}(8,10,24,26,27)$; 	

\item [\rm (46)]  $C_{54}(2,3,6,9,16,20)$, $C_{54}(3,4,6,9,14,22)$, $C_{54}(3,6,8,9,10,26)$; 	

\item [\rm (47)]  $C_{54}(2,3,6,12,16,20)$, $C_{54}(3,4,6,12,14,22)$, $C_{54}(3,6,8,10,12,26)$; 	

\item [\rm (48)]  $C_{54}(2,3,6,15,16,20)$, $C_{54}(3,4,6,14,15,22)$, $C_{54}(3,6,8,10,15,26)$; 	

\item [\rm (49)]  $C_{54}(2,3,6,16,18,20)$, $C_{54}(3,4,6,14,18,22)$, $C_{54}(3,6,8,10,18,26)$; 	

\item [\rm (50)]  $C_{54}(2,3,6,16,20,21)$, $C_{54}(3,4,6,14,21,22)$, $C_{54}(3,6,8,10,21,26)$; 	

\item [\rm (51)]  $C_{54}(2,3,6,16,20,24)$, $C_{54}(3,4,6,14,22,24)$, $C_{54}(3,6,8,10,24,26)$; 	

\item [\rm (52)]  $C_{54}(2,3,6,16,20,27)$, $C_{54}(3,4,6,14,22,27)$, $C_{54}(3,6,8,10,26,27)$; 	

\item [\rm (53)]  $C_{54}(2,3,9,12,16,20)$, $C_{54}(3,4,9,12,14,22)$, $C_{54}(3,8,9,10,12,26)$; 	

\item [\rm (54)]  $C_{54}(2,3,9,15,16,20)$, $C_{54}(3,4,9,14,15,22)$, $C_{54}(3,8,9,10,15,26)$; 	

\item [\rm (55)]  $C_{54}(2,3,9,16,18,20)$, $C_{54}(3,4,9,14,18,22)$, $C_{54}(3,8,9,10,18,26)$; 	

\item [\rm (56)]  $C_{54}(2,3,9,16,20,21)$, $C_{54}(3,4,9,14,21,22)$, $C_{54}(3,8,9,10,21,26)$; 	

\item [\rm (57)]  $C_{54}(2,3,9,16,20,24)$, $C_{54}(3,4,9,14,22,24)$, $C_{54}(3,8,9,10,24,26)$; 	

\item [\rm (58)]  $C_{54}(2,3,9,16,20,27)$, $C_{54}(3,4,9,14,22,27)$, $C_{54}(3,8,9,10,26,27)$; 	

\item [\rm (59)]  $C_{54}(2,3,12,15,16,20)$, $C_{54}(3,4,12,14,15,22)$, $C_{54}(3,8,10,12,15,26)$; 	

\item [\rm (60)]  $C_{54}(2,3,12,16,18,20)$, $C_{54}(3,4,12,14,18,22)$, $C_{54}(3,8,10,12,18,26)$; 	

\item [\rm (61)]  $C_{54}(2,3,12,16,20,21)$, $C_{54}(3,4,12,14,21,22)$, $C_{54}(3,8,10,12,21,26)$; 	

\item [\rm (62)]  $C_{54}(2,3,12,16,20,24)$, $C_{54}(3,4,12,14,22,24)$, $C_{54}(3,8,10,12,24,26)$; 	

\item [\rm (63)]  $C_{54}(2,3,12,16,20,27)$, $C_{54}(3,4,12,14,22,27)$, $C_{54}(3,8,10,12,26,27)$; 	

\item [\rm (64)]  $C_{54}(2,3,15,16,18,20)$, $C_{54}(3,4,14,15,18,22)$, $C_{54}(3,8,10,15,18,26)$; 	

\item [\rm (65)]  $C_{54}(2,3,15,16,20,21)$, $C_{54}(3,4,14,15,21,22)$, $C_{54}(3,8,10,15,21,26)$; 	

\item [\rm (66)]  $C_{54}(2,3,15,16,20,24)$, $C_{54}(3,4,14,15,22,24)$, $C_{54}(3,8,10,15,24,26)$; 	

\item [\rm (67)]  $C_{54}(2,3,15,16,20,27)$, $C_{54}(3,4,14,15,22,27)$, $C_{54}(3,8,10,15,26,27)$; 	

\item [\rm (68)]  $C_{54}(2,3,16,18,20,21)$, $C_{54}(3,4,14,18,21,22)$, $C_{54}(3,8,10,18,21,26)$; 	

\item [\rm (69)]  $C_{54}(2,3,16,18,20,24)$, $C_{54}(3,4,14,18,22,24)$, $C_{54}(3,8,10,18,24,26)$; 	

\item [\rm (70)]  $C_{54}(2,3,16,18,20,27)$, $C_{54}(3,4,14,18,22,27)$, $C_{54}(3,8,10,18,26,27)$; 	

\item [\rm (71)]  $C_{54}(2,3,16,20,21,24)$, $C_{54}(3,4,14,21,22,24)$, $C_{54}(3,8,10,21,24,26)$; 	

\item [\rm (72)]  $C_{54}(2,3,16,20,21,27)$, $C_{54}(3,4,14,21,22,27)$, $C_{54}(3,8,10,21,26,27)$; 	

\item [\rm (73)]  $C_{54}(2,3,16,20,24,27)$, $C_{54}(3,4,14,22,24,27)$, $C_{54}(3,8,10,24,26,27)$; 	

\item [\rm (74)]   $C_{54}(2,6,9,12,16,20)$, $C_{54}(4,6,9,12,14,22)$, $C_{54}(6,8,9,10,12,26)$; 	

\item [\rm (75)]   $C_{54}(2,6,9,15,16,20)$, $C_{54}(4,6,9,14,15,22)$, $C_{54}(6,8,9,10,15,26)$; 	

\item [\rm (76)]   $C_{54}(2,6,9,16,18,20)$, $C_{54}(4,6,9,14,18,22)$, $C_{54}(6,8,9,10,18,26)$; 	

\item [\rm (77)]   $C_{54}(2,6,9,16,20,21)$, $C_{54}(4,6,9,14,21,22)$, $C_{54}(6,8,9,10,21,26)$; 	

\item [\rm (78)]   $C_{54}(2,6,9,16,20,24)$, $C_{54}(4,6,9,14,22,24)$, $C_{54}(6,8,9,10,24,26)$; 	

\item [\rm (79)]   $C_{54}(2,6,9,16,20,27)$, $C_{54}(4,6,9,14,22,27)$, $C_{54}(6,8,9,10,26,27)$; 	

\item [\rm (80)]   $C_{54}(2,6,12,15,16,20)$, $C_{54}(4,6,12,14,15,22)$, $C_{54}(6,8,10,12,15,26)$; 	

\item [\rm (81)]   $C_{54}(2,6,12,16,18,20)$, $C_{54}(4,6,12,14,18,22)$, $C_{54}(6,8,10,12,18,26)$; 	

\item [\rm (82)]   $C_{54}(2,6,12,16,20,21)$, $C_{54}(4,6,12,14,21,22)$, $C_{54}(6,8,10,12,21,26)$; 	

\item [\rm (83)]   $C_{54}(2,6,12,16,20,24)$, $C_{54}(4,6,12,14,22,24)$, $C_{54}(6,8,10,12,24,26)$; 	

\item [\rm (84)]   $C_{54}(2,6,12,16,20,27)$, $C_{54}(4,6,12,14,22,27)$, $C_{54}(6,8,10,12,26,27)$; 	

\item [\rm (85)]  $C_{54}(2,6,15,16,18,20)$, $C_{54}(4,6,14,15,18,22)$, $C_{54}(6,8,10,15,18,26)$; 	

\item [\rm (86)]  $C_{54}(2,6,15,16,20,21)$, $C_{54}(4,6,14,15,21,22)$, $C_{54}(6,8,10,15,21,26)$; 	

\item [\rm (87)]  $C_{54}(2,6,15,16,20,24)$, $C_{54}(4,6,14,15,22,24)$, $C_{54}(6,8,10,15,24,26)$; 	

\item [\rm (88)]  $C_{54}(2,6,15,16,20,27)$, $C_{54}(4,6,14,15,22,27)$, $C_{54}(6,8,10,15,26,27)$; 	

\item [\rm (89)]  $C_{54}(2,6,16,18,20,21)$, $C_{54}(4,6,14,18,21,22)$, $C_{54}(6,8,10,18,21,26)$; 	

\item [\rm (90)]  $C_{54}(2,6,16,18,20,24)$, $C_{54}(4,6,14,18,22,24)$, $C_{54}(6,8,10,18,24,26)$; 	

\item [\rm (91)]  $C_{54}(2,6,16,18,20,27)$, $C_{54}(4,6,14,18,22,27)$, $C_{54}(6,8,10,18,26,27)$; 	

\item [\rm (92)]  $C_{54}(2,6,16,20,21,24)$, $C_{54}(4,6,14,21,22,24)$, $C_{54}(6,8,10,21,24,26)$; 	

\item [\rm (93)]  $C_{54}(2,6,16,20,21,27)$, $C_{54}(4,6,14,21,22,27)$, $C_{54}(6,8,10,21,26,27)$; 	

\item [\rm (94)]  $C_{54}(2,6,16,20,24,27)$, $C_{54}(4,6,14,22,24,27)$, $C_{54}(6,8,10,24,26,27)$; 	

\item [\rm (95)]    $C_{54}(2,9,12,15,16,20)$, $C_{54}(4,9,12,14,15,22)$, $C_{54}(8,9,10,12,15,26)$; 	

\item [\rm (96)]    $C_{54}(2,9,12,16,18,20)$, $C_{54}(4,9,12,14,18,22)$, $C_{54}(8,9,10,12,18,26)$; 	

\item [\rm (97)]    $C_{54}(2,9,12,16,20,21)$, $C_{54}(4,9,12,14,21,22)$, $C_{54}(8,9,10,12,21,26)$; 	

\item [\rm (98)]    $C_{54}(2,9,12,16,20,24)$, $C_{54}(4,9,12,14,22,24)$, $C_{54}(8,9,10,12,24,26)$; 	

\item [\rm (99)]    $C_{54}(2,9,12,16,20,27)$, $C_{54}(4,9,12,14,22,27)$, $C_{54}(8,9,10,12,26,27)$; 	

\item [\rm (100)]    $C_{54}(2,9,15,16,18,20)$, $C_{54}(4,9,14,15,18,22)$, $C_{54}(8,9,10,15,18,26)$; 	

\item [\rm (101)]    $C_{54}(2,9,15,16,20,21)$, $C_{54}(4,9,14,15,21,22)$, $C_{54}(8,9,10,15,21,26)$; 	

\item [\rm (102)]    $C_{54}(2,9,15,16,20,24)$, $C_{54}(4,9,14,15,22,24)$, $C_{54}(8,9,10,15,24,26)$; 	

\item [\rm (103)]    $C_{54}(2,9,15,16,20,27)$, $C_{54}(4,9,14,15,22,27)$, $C_{54}(8,9,10,15,26,27)$; 	

\item [\rm (104)]    $C_{54}(2,9,16,18,20,21)$, $C_{54}(4,9,14,18,21,22)$, $C_{54}(8,9,10,18,21,26)$; 	

\item [\rm (105)]    $C_{54}(2,9,16,18,20,24)$, $C_{54}(4,9,14,18,22,24)$, $C_{54}(8,9,10,18,24,26)$; 	

\item [\rm (106)]    $C_{54}(2,9,16,18,20,27)$, $C_{54}(4,9,14,18,22,27)$, $C_{54}(8,9,10,18,26,27)$; 	

\item [\rm (107)]    $C_{54}(2,9,16,20,21,24)$, $C_{54}(4,9,14,21,22,24)$, $C_{54}(8,9,10,21,24,26)$; 	

\item [\rm (108)]    $C_{54}(2,9,16,20,21,27)$, $C_{54}(4,9,14,21,22,27)$, $C_{54}(8,9,10,21,26,27)$; 	

\item [\rm (109)]    $C_{54}(2,9,16,20,24,27)$, $C_{54}(4,9,14,22,24,27)$, $C_{54}(8,9,10,24,26,27)$; 	

\item [\rm (110)]  $C_{54}(2,12,15,16,18,20)$, $C_{54}(4,12,14,15,18,22)$, $C_{54}(8,10,12,15,18,26)$; 	

\item [\rm (111)]  $C_{54}(2,12,15,16,20,21)$, $C_{54}(4,12,14,15,21,22)$, $C_{54}(8,10,12,15,21,26)$; 	

\item [\rm (112)]  $C_{54}(2,12,15,16,20,24)$, $C_{54}(4,12,14,15,22,24)$, $C_{54}(8,10,12,15,24,26)$; 	

\item [\rm (113)]  $C_{54}(2,12,15,16,20,27)$, $C_{54}(4,12,14,15,22,27)$, $C_{54}(8,10,12,15,26,27)$; 	

\item [\rm (114)]  $C_{54}(2,12,16,18,20,21)$, $C_{54}(4,12,14,18,21,22)$, $C_{54}(8,10,12,18,21,26)$; 	

\item [\rm (115)]  $C_{54}(2,12,16,18,20,24)$, $C_{54}(4,12,14,18,22,24)$, $C_{54}(8,10,12,18,24,26)$; 	

\item [\rm (116)]  $C_{54}(2,12,16,18,20,27)$, $C_{54}(4,12,14,18,22,27)$, $C_{54}(8,10,12,18,26,27)$; 	

\item [\rm (117)]  $C_{54}(2,12,16,20,21,24)$, $C_{54}(4,12,14,21,22,24)$, $C_{54}(8,10,12,21,24,26)$; 	

\item [\rm (118)]  $C_{54}(2,12,16,20,21,27)$, $C_{54}(4,12,14,21,22,27)$, $C_{54}(8,10,12,21,26,27)$; 	

\item [\rm (119)]  $C_{54}(2,12,16,20,24,27)$, $C_{54}(4,12,14,22,24,27)$, $C_{54}(8,10,12,24,26,27)$; 	

\item [\rm (120)]  $C_{54}(2,15,16,18,20,21)$, $C_{54}(4,14,15,18,21,22)$, $C_{54}(8,10,15,18,21,26)$; 	

\item [\rm (121)]  $C_{54}(2,15,16,18,20,24)$, $C_{54}(4,14,15,18,22,24)$, $C_{54}(8,10,15,18,24,26)$; 	

\item [\rm (122)]  $C_{54}(2,15,16,18,20,27)$, $C_{54}(4,14,15,18,22,27)$, $C_{54}(8,10,15,18,26,27)$; 	

\item [\rm (123)]  $C_{54}(2,15,16,20,21,24)$, $C_{54}(4,14,15,21,22,24)$, $C_{54}(8,10,15,21,24,26)$; 	

\item [\rm (124)]  $C_{54}(2,15,16,20,21,27)$, $C_{54}(4,14,15,21,22,27)$, $C_{54}(8,10,15,21,26,27)$; 	

\item [\rm (125)]  $C_{54}(2,15,16,20,24,27)$, $C_{54}(4,14,15,22,24,27)$, $C_{54}(8,10,15,24,26,27)$; 	

\item [\rm (126)]  $C_{54}(2,16,18,20,21,24)$, $C_{54}(4,14,18,21,22,24)$, $C_{54}(8,10,18,21,24,26)$; 	

\item [\rm (127)]  $C_{54}(2,16,18,20,21,27)$, $C_{54}(4,14,18,21,22,27)$, $C_{54}(8,10,18,21,26,27)$; 	

\item [\rm (128)]  $C_{54}(2,16,18,20,24,27)$, $C_{54}(4,14,18,22,24,27)$, $C_{54}(8,10,18,24,26,27)$; 	

\item [\rm (129)]  $C_{54}(2,16,20,21,24,27)$, $C_{54}(4,14,21,22,24,27)$, $C_{54}(8,10,21,24,26,27)$; 	

\item [\rm (130)]  $C_{54}(2,3,6,9,12,16,20)$, $C_{54}(3,4,6,9,12,14,22)$, $C_{54}(3,6,8,9,10,12,26)$; 	

\item [\rm (131)]  $C_{54}(2,3,6,9,15,16,20)$, $C_{54}(3,4,6,9,14,15,22)$, $C_{54}(3,6,8,9,10,15,26)$; 	

\item [\rm (132)]  $C_{54}(2,3,6,9,16,18,20)$, $C_{54}(3,4,6,9,14,18,22)$, $C_{54}(3,6,8,9,10,18,26)$; 	

\item [\rm (133)]  $C_{54}(2,3,6,9,16,20,21)$, $C_{54}(3,4,6,9,14,21,22)$, $C_{54}(3,6,8,9,10,21,26)$; 	

\item [\rm (134)]  $C_{54}(2,3,6,9,16,20,24)$, $C_{54}(3,4,6,9,14,22,24)$, $C_{54}(3,6,8,9,10,24,26)$; 	

\item [\rm (135)]  $C_{54}(2,3,6,9,16,20,27)$, $C_{54}(3,4,6,9,14,22,27)$, $C_{54}(3,6,8,9,10,26,27)$; 	

\item [\rm (136)]  $C_{54}(2,3,6,12,15,16,20)$, $C_{54}(3,4,6,12,14,15,22)$, $C_{54}(3,6,8,10,12,15,26)$; 	

\item [\rm (137)]  $C_{54}(2,3,6,12,16,18,20)$, $C_{54}(3,4,6,12,14,18,22)$, $C_{54}(3,6,8,10,12,18,26)$; 	

\item [\rm (138)]  $C_{54}(2,3,6,12,16,20,21)$, $C_{54}(3,4,6,12,14,21,22)$, $C_{54}(3,6,8,10,12,21,26)$; 	

\item [\rm (139)]  $C_{54}(2,3,6,12,16,20,24)$, $C_{54}(3,4,6,12,14,22,24)$, $C_{54}(3,6,8,10,12,24,26)$; 	

\item [\rm (140)]  $C_{54}(2,3,6,12,16,20,27)$, $C_{54}(3,4,6,12,14,22,27)$, $C_{54}(3,6,8,10,12,26,27)$; 	

\item [\rm (141)]  $C_{54}(2,3,6,15,16,18,20)$, $C_{54}(3,4,6,14,15,18,22)$, $C_{54}(3,6,8,10,15,18,26)$; 	

\item [\rm (142)]  $C_{54}(2,3,6,15,16,20,21)$, $C_{54}(3,4,6,14,15,21,22)$, $C_{54}(3,6,8,10,15,21,26)$; 	

\item [\rm (143)]  $C_{54}(2,3,6,15,16,20,24)$, $C_{54}(3,4,6,14,15,22,24)$, $C_{54}(3,6,8,10,15,24,26)$; 	

\item [\rm (144)]  $C_{54}(2,3,6,15,16,20,27)$, $C_{54}(3,4,6,14,15,22,27)$, $C_{54}(3,6,8,10,15,26,27)$; 	

\item [\rm (145)]  $C_{54}(2,3,6,16,18,20,21)$, $C_{54}(3,4,6,14,18,21,22)$, $C_{54}(3,6,8,10,18,21,26)$; 	

\item [\rm (146)]  $C_{54}(2,3,6,16,18,20,24)$, $C_{54}(3,4,6,14,18,22,24)$, $C_{54}(3,6,8,10,18,24,26)$; 	

\item [\rm (147)]  $C_{54}(2,3,6,16,18,20,27)$, $C_{54}(3,4,6,14,18,22,27)$, $C_{54}(3,6,8,10,18,26,27)$; 	

\item [\rm (148)]  $C_{54}(2,3,6,16,20,21,24)$, $C_{54}(3,4,6,14,21,22,24)$, $C_{54}(3,6,8,10,21,24,26)$; 	

\item [\rm (149)]  $C_{54}(2,3,6,16,20,21,27)$, $C_{54}(3,4,6,14,21,22,27)$, $C_{54}(3,6,8,10,21,26,27)$; 	

\item [\rm (150)]  $C_{54}(2,3,6,16,20,24,27)$, $C_{54}(3,4,6,14,22,24,27)$, $C_{54}(3,6,8,10,24,26,27)$; 	

\item [\rm (151)]  $C_{54}(2,3,9,12,15,16,20)$, $C_{54}(3,4,9,12,14,15,22)$, $C_{54}(3,8,9,10,12,15,26)$; 	

\item [\rm (152)]  $C_{54}(2,3,9,12,16,18,20)$, $C_{54}(3,4,9,12,14,18,22)$, $C_{54}(3,8,9,10,12,18,26)$; 	

\item [\rm (153)]  $C_{54}(2,3,9,12,16,20,21)$, $C_{54}(3,4,9,12,14,21,22)$, $C_{54}(3,8,9,10,12,21,26)$; 	

\item [\rm (154)]  $C_{54}(2,3,9,12,16,20,24)$, $C_{54}(3,4,9,12,14,22,24)$, $C_{54}(3,8,9,10,12,24,26)$; 	

\item [\rm (155)]  $C_{54}(2,3,9,12,16,20,27)$, $C_{54}(3,4,9,12,14,22,27)$, $C_{54}(3,8,9,10,12,26,27)$; 	

\item [\rm (156)]  $C_{54}(2,3,9,15,16,18,20)$, $C_{54}(3,4,9,14,15,18,22)$, $C_{54}(3,8,9,10,15,18,26)$; 	

\item [\rm (157)]  $C_{54}(2,3,9,15,16,20,21)$, $C_{54}(3,4,9,14,15,21,22)$, $C_{54}(3,8,9,10,15,21,26)$; 	

\item [\rm (158)]  $C_{54}(2,3,9,15,16,20,24)$, $C_{54}(3,4,9,14,15,22,24)$, $C_{54}(3,8,9,10,15,24,26)$; 	

\item [\rm (159)]  $C_{54}(2,3,9,15,16,20,27)$, $C_{54}(3,4,9,14,15,22,27)$, $C_{54}(3,8,9,10,15,26,27)$; 	

\item [\rm (160)]  $C_{54}(2,3,9,16,18,20,21)$, $C_{54}(3,4,9,14,18,21,22)$, $C_{54}(3,8,9,10,18,21,26)$; 	

\item [\rm (161)]  $C_{54}(2,3,9,16,18,20,24)$, $C_{54}(3,4,9,14,18,22,24)$, $C_{54}(3,8,9,10,18,24,26)$; 	

\item [\rm (162)]  $C_{54}(2,3,9,16,18,20,27)$, $C_{54}(3,4,9,14,18,22,27)$, $C_{54}(3,8,9,10,18,26,27)$; 	

\item [\rm (163)]  $C_{54}(2,3,9,16,20,21,24)$, $C_{54}(3,4,9,14,21,22,24)$, $C_{54}(3,8,9,10,21,24,26)$; 	

\item [\rm (164)]  $C_{54}(2,3,9,16,20,21,27)$, $C_{54}(3,4,9,14,21,22,27)$, $C_{54}(3,8,9,10,21,26,27)$; 	

\item [\rm (165)]  $C_{54}(2,3,9,16,20,24,27)$, $C_{54}(3,4,9,14,22,24,27)$, $C_{54}(3,8,9,10,24,26,27)$; 	

\item [\rm (166)]  $C_{54}(2,3,12,15,16,18,20)$, $C_{54}(3,4,12,14,15,18,22)$, $C_{54}(3,8,10,12,15,18,26)$; 	

\item [\rm (167)]  $C_{54}(2,3,12,15,16,20,21)$, $C_{54}(3,4,12,14,15,21,22)$, $C_{54}(3,8,10,12,15,21,26)$; 	

\item [\rm (168)]  $C_{54}(2,3,12,15,16,20,24)$, $C_{54}(3,4,12,14,15,22,24)$, $C_{54}(3,8,10,12,15,24,26)$; 	

\item [\rm (169)]  $C_{54}(2,3,12,15,16,20,27)$, $C_{54}(3,4,12,14,15,22,27)$, $C_{54}(3,8,10,12,15,26,27)$; 	

\item [\rm (170)]  $C_{54}(2,3,12,16,18,20,21)$, $C_{54}(3,4,12,14,18,21,22)$, $C_{54}(3,8,10,12,18,21,26)$; 	

\item [\rm (171)]  $C_{54}(2,3,12,16,18,20,24)$, $C_{54}(3,4,12,14,18,22,24)$, $C_{54}(3,8,10,12,18,24,26)$; 	

\item [\rm (172)]  $C_{54}(2,3,12,16,18,20,27)$, $C_{54}(3,4,12,14,18,22,27)$, $C_{54}(3,8,10,12,18,26,27)$; 	

\item [\rm (173)]  $C_{54}(2,3,12,16,20,21,24)$, $C_{54}(3,4,12,14,21,22,24)$, $C_{54}(3,8,10,12,21,24,26)$; 	

\item [\rm (174)]  $C_{54}(2,3,12,16,20,21,27)$, $C_{54}(3,4,12,14,21,22,27)$, $C_{54}(3,8,10,12,21,26,27)$; 	

\item [\rm (175)]  $C_{54}(2,3,12,16,20,24,27)$, $C_{54}(3,4,12,14,22,24,27)$, $C_{54}(3,8,10,12,24,26,27)$; 	

\item [\rm (176)]  $C_{54}(2,3,15,16,18,20,21)$, $C_{54}(3,4,14,15,18,21,22)$, $C_{54}(3,8,10,15,18,21,26)$; 	

\item [\rm (177)]  $C_{54}(2,3,15,16,18,20,24)$, $C_{54}(3,4,14,15,18,22,24)$, $C_{54}(3,8,10,15,18,24,26)$; 	

\item [\rm (178)]  $C_{54}(2,3,15,16,18,20,27)$, $C_{54}(3,4,14,15,18,22,27)$, $C_{54}(3,8,10,15,18,26,27)$; 	

\item [\rm (179)]  $C_{54}(2,3,15,16,20,21,24)$, $C_{54}(3,4,14,15,21,22,24)$, $C_{54}(3,8,10,15,21,24,26)$; 	

\item [\rm (180)]  $C_{54}(2,3,15,16,20,21,27)$, $C_{54}(3,4,14,15,21,22,27)$, $C_{54}(3,8,10,15,21,26,27)$; 	

\item [\rm (181)]  $C_{54}(2,3,15,16,20,24,27)$, $C_{54}(3,4,14,15,22,24,27)$, $C_{54}(3,8,10,15,24,26,27)$; 	

\item [\rm (182)]  $C_{54}(2,3,16,18,20,21,24)$, $C_{54}(3,4,14,18,21,22,24)$, $C_{54}(3,8,10,18,21,24,26)$; 	

\item [\rm (183)]  $C_{54}(2,3,16,18,20,21,27)$, $C_{54}(3,4,14,18,21,22,27)$, $C_{54}(3,8,10,18,21,26,27)$; 	

\item [\rm (184)]  $C_{54}(2,3,16,18,20,24,27)$, $C_{54}(3,4,14,18,22,24,27)$, $C_{54}(3,8,10,18,24,26,27)$; 	

\item [\rm (185)]  $C_{54}(2,3,16,20,21,24,27)$, $C_{54}(3,4,14,21,22,24,27)$, $C_{54}(3,8,10,21,24,26,27)$; 	

\item [\rm (186)]  $C_{54}(2,6,9,12,15,16,20)$, $C_{54}(4,6,9,12,14,15,22)$, $C_{54}(6,8,9,10,12,15,26)$; 	

\item [\rm (187)]  $C_{54}(2,6,9,12,16,18,20)$, $C_{54}(4,6,9,12,14,18,22)$, $C_{54}(6,8,9,10,12,18,26)$; 	

\item [\rm (188)]  $C_{54}(2,6,9,12,16,20,21)$, $C_{54}(4,6,9,12,14,21,22)$, $C_{54}(6,8,9,10,12,21,26)$; 	

\item [\rm (189)]  $C_{54}(2,6,9,12,16,20,24)$, $C_{54}(4,6,9,12,14,22,24)$, $C_{54}(6,8,9,10,12,24,26)$; 	

\item [\rm (190)]  $C_{54}(2,6,9,12,16,20,27)$, $C_{54}(4,6,9,12,14,22,27)$, $C_{54}(6,8,9,10,12,26,27)$; 	

\item [\rm (191)]  $C_{54}(2,6,9,15,16,18,20)$, $C_{54}(4,6,9,14,15,18,22)$, $C_{54}(6,8,9,10,15,18,26)$; 	

\item [\rm (192)]  $C_{54}(2,6,9,15,16,20,21)$, $C_{54}(4,6,9,14,15,21,22)$, $C_{54}(6,8,9,10,15,21,26)$; 	

\item [\rm (193)]  $C_{54}(2,6,9,15,16,20,24)$, $C_{54}(4,6,9,14,15,22,24)$, $C_{54}(6,8,9,10,15,24,26)$; 	

\item [\rm (194)]  $C_{54}(2,6,9,15,16,20,27)$, $C_{54}(4,6,9,14,15,22,27)$, $C_{54}(6,8,9,10,15,26,27)$; 	

\item [\rm (195)]  $C_{54}(2,6,9,16,18,20,21)$, $C_{54}(4,6,9,14,18,21,22)$, $C_{54}(6,8,9,10,18,21,26)$; 	

\item [\rm (196)]  $C_{54}(2,6,9,16,18,20,24)$, $C_{54}(4,6,9,14,18,22,24)$, $C_{54}(6,8,9,10,18,24,26)$; 	

\item [\rm (197)]  $C_{54}(2,6,9,16,18,20,27)$, $C_{54}(4,6,9,14,18,22,27)$, $C_{54}(6,8,9,10,18,26,27)$; 	

\item [\rm (198)]  $C_{54}(2,6,9,16,20,21,24)$, $C_{54}(4,6,9,14,21,22,24)$, $C_{54}(6,8,9,10,21,24,26)$; 	

\item [\rm (199)]  $C_{54}(2,6,9,16,20,21,27)$, $C_{54}(4,6,9,14,21,22,27)$, $C_{54}(6,8,9,10,21,26,27)$; 	

\item [\rm (200)]  $C_{54}(2,6,9,16,20,24,27)$, $C_{54}(4,6,9,14,22,24,27)$, $C_{54}(6,8,9,10,24,26,27)$; 	

\item [\rm (201)]  $C_{54}(2,6,12,15,16,18,20)$, $C_{54}(4,6,12,14,15,18,22)$, $C_{54}(6,8,10,12,15,18,26)$; 	

\item [\rm (202)]  $C_{54}(2,6,12,15,16,20,21)$, $C_{54}(4,6,12,14,15,21,22)$, $C_{54}(6,8,10,12,15,21,26)$; 	

\item [\rm (203)]  $C_{54}(2,6,12,15,16,20,24)$, $C_{54}(4,6,12,14,15,22,24)$, $C_{54}(6,8,10,12,15,24,26)$; 	

\item [\rm (204)]  $C_{54}(2,6,12,15,16,20,27)$, $C_{54}(4,6,12,14,15,22,27)$, $C_{54}(6,8,10,12,15,26,27)$; 	

\item [\rm (205)]  $C_{54}(2,6,12,16,18,20,21)$, $C_{54}(4,6,12,14,18,21,22)$, $C_{54}(6,8,10,12,18,21,26)$; 	

\item [\rm (206)]  $C_{54}(2,6,12,16,18,20,24)$, $C_{54}(4,6,12,14,18,22,24)$, $C_{54}(6,8,10,12,18,24,26)$; 	

\item [\rm (207)]  $C_{54}(2,6,12,16,18,20,27)$, $C_{54}(4,6,12,14,18,22,27)$, $C_{54}(6,8,10,12,18,26,27)$; 	

\item [\rm (208)]  $C_{54}(2,6,12,16,20,21,24)$, $C_{54}(4,6,12,14,21,22,24)$, $C_{54}(6,8,10,12,21,24,26)$; 	

\item [\rm ( 209)]  $C_{54}(2,6,12,16,20,21,27)$, $C_{54}(4,6,12,14,21,22,27)$, $C_{54}(6,8,10,12,21,26,27)$; 	

\item [\rm ( 210)]  $C_{54}(2,6,12,16,20,24,27)$, $C_{54}(4,6,12,14,22,24,27)$, $C_{54}(6,8,10,12,24,26,27)$; 	

\item [\rm ( 211)]  $C_{54}(2,6,15,16,18,20,21)$, $C_{54}(4,6,14,15,18,21,22)$, $C_{54}(6,8,10,15,18,21,26)$; 	

\item [\rm ( 212)]  $C_{54}(2,6,15,16,18,20,24)$, $C_{54}(4,6,14,15,18,22,24)$, $C_{54}(6,8,10,15,18,24,26)$; 	

\item [\rm ( 213)]  $C_{54}(2,6,15,16,18,20,27)$, $C_{54}(4,6,14,15,18,22,27)$, $C_{54}(6,8,10,15,18,26,27)$; 	

\item [\rm ( 214)]  $C_{54}(2,6,15,16,20,21,24)$, $C_{54}(4,6,14,15,21,22,24)$, $C_{54}(6,8,10,15,21,24,26)$; 	

\item [\rm ( 215)]  $C_{54}(2,6,15,16,20,21,27)$, $C_{54}(4,6,14,15,21,22,27)$, $C_{54}(6,8,10,15,21,26,27)$; 	

\item [\rm ( 216)]  $C_{54}(2,6,15,16,20,24,27)$, $C_{54}(4,6,14,15,22,24,27)$, $C_{54}(6,8,10,15,24,26,27)$; 	

\item [\rm ( 217)]  $C_{54}(2,6,16,18,20,21,24)$, $C_{54}(4,6,14,18,21,22,24)$, $C_{54}(6,8,10,18,21,24,26)$; 	

\item [\rm ( 218)]  $C_{54}(2,6,16,18,20,21,27)$, $C_{54}(4,6,14,18,21,22,27)$, $C_{54}(6,8,10,18,21,26,27)$; 	

\item [\rm ( 219)]  $C_{54}(2,6,16,18,20,24,27)$, $C_{54}(4,6,14,18,22,24,27)$, $C_{54}(6,8,10,18,24,26,27)$; 	

\item [\rm ( 220)]  $C_{54}(2,6,16,20,21,24,27)$, $C_{54}(4,6,14,21,22,24,27)$, $C_{54}(6,8,10,21,24,26,27)$; 	

\item [\rm ( 221)]  $C_{54}(2,9,12,15,16,18,20)$, $C_{54}(4,9,12,14,15,18,22)$, $C_{54}(8,9,10,12,15,18,26)$; 	

\item [\rm ( 222)]  $C_{54}(2,9,12,15,16,20,21)$, $C_{54}(4,9,12,14,15,21,22)$, $C_{54}(8,9,10,12,15,21,26)$; 	

\item [\rm ( 223)]  $C_{54}(2,9,12,15,16,20,24)$, $C_{54}(4,9,12,14,15,22,24)$, $C_{54}(8,9,10,12,15,24,26)$; 	

\item [\rm ( 224)]  $C_{54}(2,9,12,15,16,20,27)$, $C_{54}(4,9,12,14,15,22,27)$, $C_{54}(8,9,10,12,15,26,27)$; 	

\item [\rm ( 225)]  $C_{54}(2,9,12,16,18,20,21)$, $C_{54}(4,9,12,14,18,21,22)$, $C_{54}(8,9,10,12,18,21,26)$; 	

\item [\rm ( 226)]  $C_{54}(2,9,12,16,18,20,24)$, $C_{54}(4,9,12,14,18,22,24)$, $C_{54}(8,9,10,12,18,24,26)$; 	

\item [\rm ( 227)]  $C_{54}(2,9,12,16,18,20,27)$, $C_{54}(4,9,12,14,18,22,27)$, $C_{54}(8,9,10,12,18,26,27)$; 	

\item [\rm ( 228)]  $C_{54}(2,9,12,16,20,21,24)$, $C_{54}(4,9,12,14,21,22,24)$, $C_{54}(8,9,10,12,21,24,26)$; 	

\item [\rm ( 229)]  $C_{54}(2,9,12,16,20,21,27)$, $C_{54}(4,9,12,14,21,22,27)$, $C_{54}(8,9,10,12,21,26,27)$; 	

\item [\rm ( 230)]  $C_{54}(2,9,12,16,20,24,27)$, $C_{54}(4,9,12,14,22,24,27)$, $C_{54}(8,9,10,12,24,26,27)$; 	

\item [\rm ( 231)]  $C_{54}(2,9,15,16,18,20,21)$, $C_{54}(4,9,14,15,18,21,22)$, $C_{54}(8,9,10,15,18,21,26)$; 	

\item [\rm ( 232)]  $C_{54}(2,9,15,16,18,20,24)$, $C_{54}(4,9,14,15,18,22,24)$, $C_{54}(8,9,10,15,18,24,26)$; 	

\item [\rm ( 233)]  $C_{54}(2,9,15,16,18,20,27)$, $C_{54}(4,9,14,15,18,22,27)$, $C_{54}(8,9,10,15,18,26,27)$; 	

\item [\rm ( 234)]  $C_{54}(2,9,15,16,20,21,24)$, $C_{54}(4,9,14,15,21,22,24)$, $C_{54}(8,9,10,15,21,24,26)$; 	

\item [\rm ( 235)]  $C_{54}(2,9,15,16,20,21,27)$, $C_{54}(4,9,14,15,21,22,27)$, $C_{54}(8,9,10,15,21,26,27)$; 	

\item [\rm ( 236)]  $C_{54}(2,9,15,16,20,24,27)$, $C_{54}(4,9,14,15,22,24,27)$, $C_{54}(8,9,10,15,24,26,27)$; 	

\item [\rm ( 237)]  $C_{54}(2,9,16,18,20,21,24)$, $C_{54}(4,9,14,18,21,22,24)$, $C_{54}(8,9,10,18,21,24,26)$; 	

\item [\rm ( 238)]  $C_{54}(2,9,16,18,20,21,27)$, $C_{54}(4,9,14,18,21,22,27)$, $C_{54}(8,9,10,18,21,26,27)$; 	

\item [\rm ( 239)]  $C_{54}(2,9,16,18,20,24,27)$, $C_{54}(4,9,14,18,22,24,27)$, $C_{54}(8,9,10,18,24,26,27)$; 	

\item [\rm ( 240)]  $C_{54}(2,9,16,20,21,24,27)$, $C_{54}(4,9,14,21,22,24,27)$, $C_{54}(8,9,10,21,24,26,27)$; 	

\item [\rm ( 241)]  $C_{54}(2,12,15,16,18,20,21)$, $C_{54}(4,12,14,15,18,21,22)$, $C_{54}(8,10,12,15,18,21,26)$; 	

\item [\rm ( 242)]  $C_{54}(2,12,15,16,18,20,24)$, $C_{54}(4,12,14,15,18,22,24)$, $C_{54}(8,10,12,15,18,24,26)$; 	

\item [\rm ( 243)]  $C_{54}(2,12,15,16,18,20,27)$, $C_{54}(4,12,14,15,18,22,27)$, $C_{54}(8,10,12,15,18,26,27)$; 	

\item [\rm ( 244)]  $C_{54}(2,12,15,16,20,21,24)$, $C_{54}(4,12,14,15,21,22,24)$, $C_{54}(8,10,12,15,21,24,26)$; 	

\item [\rm ( 245)]  $C_{54}(2,12,15,16,20,21,27)$, $C_{54}(4,12,14,15,21,22,27)$, $C_{54}(8,10,12,15,21,26,27)$; 	

\item [\rm ( 246)]  $C_{54}(2,12,15,16,20,24,27)$, $C_{54}(4,12,14,15,22,24,27)$, $C_{54}(8,10,12,15,24,26,27)$; 	

\item [\rm ( 247)]  $C_{54}(2,12,16,18,20,21,24)$, $C_{54}(4,12,14,18,21,22,24)$, $C_{54}(8,10,12,18,21,24,26)$; 	

\item [\rm ( 248)]  $C_{54}(2,12,16,18,20,21,27)$, $C_{54}(4,12,14,18,21,22,27)$, $C_{54}(8,10,12,18,21,26,27)$; 	

\item [\rm ( 249)]  $C_{54}(2,12,16,18,20,24,27)$, $C_{54}(4,12,14,18,22,24,27)$, $C_{54}(8,10,12,18,24,26,27)$; 	

\item [\rm ( 250)]  $C_{54}(2,12,16,20,21,24,27)$, $C_{54}(4,12,14,21,22,24,27)$, $C_{54}(8,10,12,21,24,26,27)$; 	

\item [\rm (251)]   $C_{54}(2,15,16,18,20,21,24)$, $C_{54}(4,14,15,18,21,22,24)$, $C_{54}(8,10,15,18,21,24,26)$; 	

\item [\rm (252)]   $C_{54}(2,15,16,18,20,21,27)$, $C_{54}(4,14,15,18,21,22,27)$, $C_{54}(8,10,15,18,21,26,27)$; 	

\item [\rm (253)]  $C_{54}(2,15,16,18,20,24,27)$, $C_{54}(4,14,15,18,22,24,27)$, $C_{54}(8,10,15,18,24,26,27)$; 	

\item [\rm (254)]  $C_{54}(2,15,16,20,21,24,27)$, $C_{54}(4,14,15,21,22,24,27)$, $C_{54}(8,10,15,21,24,26,27)$; 	

\item [\rm (255)]  $C_{54}(2,16,18,20,21,24,27)$, $C_{54}(4,14,18,21,22,24,27)$, $C_{54}(8,10,18,21,24,26,27)$; 	

\item [\rm (256)]  $C_{54}(2,3,6,9,12,15,16,20)$, $C_{54}(3,4,6,9,12,14,15,22)$, $C_{54}(3,6,8,9,10,12,15,26)$; 	

\item [\rm (257)]  $C_{54}(2,3,6,9,12,16,18,20)$, $C_{54}(3,4,6,9,12,14,18,22)$, $C_{54}(3,6,8,9,10,12,18,26)$; 	

\item [\rm (258)]  $C_{54}(2,3,6,9,12,16,20,21)$, $C_{54}(3,4,6,9,12,14,21,22)$, $C_{54}(3,6,8,9,10,12,21,26)$; 	

\item [\rm (259)]  $C_{54}(2,3,6,9,12,16,20,24)$, $C_{54}(3,4,6,9,12,14,22,24)$, $C_{54}(3,6,8,9,10,12,24,26)$; 	

\item [\rm (260)]  $C_{54}(2,3,6,9,12,16,20,27)$, $C_{54}(3,4,6,9,12,14,22,27)$, $C_{54}(3,6,8,9,10,12,26,27)$; 	

\item [\rm (261)]  $C_{54}(2,3,6,9,15,16,18,20)$, $C_{54}(3,4,6,9,14,15,18,22)$, $C_{54}(3,6,8,9,10,15,18,26)$; 	

\item [\rm (262)]  $C_{54}(2,3,6,9,15,16,20,21)$, $C_{54}(3,4,6,9,14,15,21,22)$, $C_{54}(3,6,8,9,10,15,21,26)$; 	

\item [\rm (263)]  $C_{54}(2,3,6,9,15,16,20,24)$, $C_{54}(3,4,6,9,14,15,22,24)$, $C_{54}(3,6,8,9,10,15,24,26)$; 	

\item [\rm (264)]  $C_{54}(2,3,6,9,15,16,20,27)$, $C_{54}(3,4,6,9,14,15,22,27)$, $C_{54}(3,6,8,9,10,15,26,27)$; 	

\item [\rm (265)]  $C_{54}(2,3,6,9,16,18,20,21)$, $C_{54}(3,4,6,9,14,18,21,22)$, $C_{54}(3,6,8,9,10,18,21,26)$; 	

\item [\rm (266)]  $C_{54}(2,3,6,9,16,18,20,24)$, $C_{54}(3,4,6,9,14,18,22,24)$, $C_{54}(3,6,8,9,10,18,24,26)$; 	

\item [\rm (267)]  $C_{54}(2,3,6,9,16,18,20,27)$, $C_{54}(3,4,6,9,14,18,22,27)$, $C_{54}(3,6,8,9,10,18,26,27)$; 	

\item [\rm (268)]  $C_{54}(2,3,6,9,16,20,21,24)$, $C_{54}(3,4,6,9,14,21,22,24)$, $C_{54}(3,6,8,9,10,21,24,26)$; 	

\item [\rm (269)]  $C_{54}(2,3,6,9,16,20,21,27)$, $C_{54}(3,4,6,9,14,21,22,27)$, $C_{54}(3,6,8,9,10,21,26,27)$; 	

\item [\rm (270)]  $C_{54}(2,3,6,9,16,20,24,27)$, $C_{54}(3,4,6,9,14,22,24,27)$, $C_{54}(3,6,8,9,10,24,26,27)$; 	

\item [\rm (271)]  $C_{54}(2,3,6,12,15,16,18,20)$, $C_{54}(3,4,6,12,14,15,18,22)$, $C_{54}(3,6,8,10,12,15,18,26)$; 	

\item [\rm (272)]  $C_{54}(2,3,6,12,15,16,20,21)$, $C_{54}(3,4,6,12,14,15,21,22)$, $C_{54}(3,6,8,10,12,15,21,26)$; 	

\item [\rm (273)]  $C_{54}(2,3,6,12,15,16,20,24)$, $C_{54}(3,4,6,12,14,15,22,24)$, $C_{54}(3,6,8,10,12,15,24,26)$; 	

\item [\rm (274)]  $C_{54}(2,3,6,12,15,16,20,27)$, $C_{54}(3,4,6,12,14,15,22,27)$, $C_{54}(3,6,8,10,12,15,26,27)$; 	

\item [\rm (275)]  $C_{54}(2,3,6,12,16,18,20,21)$, $C_{54}(3,4,6,12,14,18,21,22)$, $C_{54}(3,6,8,10,12,18,21,26)$; 	

\item [\rm (276)]  $C_{54}(2,3,6,12,16,18,20,24)$, $C_{54}(3,4,6,12,14,18,22,24)$, $C_{54}(3,6,8,10,12,18,24,26)$; 	

\item [\rm (277)]  $C_{54}(2,3,6,12,16,18,20,27)$, $C_{54}(3,4,6,12,14,18,22,27)$, $C_{54}(3,6,8,10,12,18,26,27)$; 	

\item [\rm (278)]  $C_{54}(2,3,6,12,16,20,21,24)$, $C_{54}(3,4,6,12,14,21,22,24)$, $C_{54}(3,6,8,10,12,21,24,26)$; 	

\item [\rm (279)]  $C_{54}(2,3,6,12,16,20,21,27)$, $C_{54}(3,4,6,12,14,21,22,27)$, $C_{54}(3,6,8,10,12,21,26,27)$; 	

\item [\rm (280)]  $C_{54}(2,3,6,12,16,20,24,27)$, $C_{54}(3,4,6,12,14,22,24,27)$, $C_{54}(3,6,8,10,12,24,26,27)$; 	

\item [\rm (281)]  $C_{54}(2,3,6,15,16,18,20,21)$, $C_{54}(3,4,6,14,15,18,21,22)$, $C_{54}(3,6,8,10,15,18,21,26)$; 	

\item [\rm (282)]  $C_{54}(2,3,6,15,16,18,20,24)$, $C_{54}(3,4,6,14,15,18,22,24)$, $C_{54}(3,6,8,10,15,18,24,26)$; 	

\item [\rm (283)]  $C_{54}(2,3,6,15,16,18,20,27)$, $C_{54}(3,4,6,14,15,18,22,27)$, $C_{54}(3,6,8,10,15,18,26,27)$; 	

\item [\rm (284)]  $C_{54}(2,3,6,15,16,20,21,24)$, $C_{54}(3,4,6,14,15,21,22,24)$, $C_{54}(3,6,8,10,15,21,24,26)$; 	

\item [\rm (285)]  $C_{54}(2,3,6,15,16,20,21,27)$, $C_{54}(3,4,6,14,15,21,22,27)$, $C_{54}(3,6,8,10,15,21,26,27)$; 	

\item [\rm (286)]  $C_{54}(2,3,6,15,16,20,24,27)$, $C_{54}(3,4,6,14,15,22,24,27)$, $C_{54}(3,6,8,10,15,24,26,27)$; 	

\item [\rm (287)]  $C_{54}(2,3,6,16,18,20,21,24)$, $C_{54}(3,4,6,14,18,21,22,24)$, $C_{54}(3,6,8,10,18,21,24,26)$; 	

\item [\rm (288)]  $C_{54}(2,3,6,16,18,20,21,27)$, $C_{54}(3,4,6,14,18,21,22,27)$, $C_{54}(3,6,8,10,18,21,26,27)$; 	

\item [\rm (289)]  $C_{54}(2,3,6,16,18,20,24,27)$, $C_{54}(3,4,6,14,18,22,24,27)$, $C_{54}(3,6,8,10,18,24,26,27)$; 	

\item [\rm (290)]  $C_{54}(2,3,6,16,20,21,24,27)$, $C_{54}(3,4,6,14,21,22,24,27)$, $C_{54}(3,6,8,10,21,24,26,27)$; 	

\item [\rm (291)]  $C_{54}(2,3,9,12,15,16,18,20)$, $C_{54}(3,4,9,12,14,15,18,22)$, $C_{54}(3,8,9,10,12,15,18,26)$; 	

\item [\rm (292)]  $C_{54}(2,3,9,12,15,16,20,21)$, $C_{54}(3,4,9,12,14,15,21,22)$, $C_{54}(3,8,9,10,12,15,21,26)$; 	

\item [\rm (293)]  $C_{54}(2,3,9,12,15,16,20,24)$, $C_{54}(3,4,9,12,14,15,22,24)$, $C_{54}(3,8,9,10,12,15,24,26)$; 	

\item [\rm (294)]  $C_{54}(2,3,9,12,15,16,20,27)$, $C_{54}(3,4,9,12,14,15,22,27)$, $C_{54}(3,8,9,10,12,15,26,27)$; 	

\item [\rm (295)]  $C_{54}(2,3,9,12,16,18,20,21)$, $C_{54}(3,4,9,12,14,18,21,22)$, $C_{54}(3,8,9,10,12,18,21,26)$; 	

\item [\rm (296)]  $C_{54}(2,3,9,12,16,18,20,24)$, $C_{54}(3,4,9,12,14,18,22,24)$, $C_{54}(3,8,9,10,12,18,24,26)$; 	

\item [\rm (297)]  $C_{54}(2,3,9,12,16,18,20,27)$, $C_{54}(3,4,9,12,14,18,22,27)$, $C_{54}(3,8,9,10,12,18,26,27)$; 	

\item [\rm (298)]  $C_{54}(2,3,9,12,16,20,21,24)$, $C_{54}(3,4,9,12,14,21,22,24)$, $C_{54}(3,8,9,10,12,21,24,26)$; 	

\item [\rm (299)]  $C_{54}(2,3,9,12,16,20,21,27)$, $C_{54}(3,4,9,12,14,21,22,27)$, $C_{54}(3,8,9,10,12,21,26,27)$; 	

\item [\rm (300)]  $C_{54}(2,3,9,12,16,20,24,27)$, $C_{54}(3,4,9,12,14,22,24,27)$, $C_{54}(3,8,9,10,12,24,26,27)$; 	

\item [\rm (301)]  $C_{54}(2,3,9,15,16,18,20,21)$, $C_{54}(3,4,9,14,15,18,21,22)$, $C_{54}(3,8,9,10,15,18,21,26)$; 	

\item [\rm (302)]  $C_{54}(2,3,9,15,16,18,20,24)$, $C_{54}(3,4,9,14,15,18,22,24)$, $C_{54}(3,8,9,10,15,18,24,26)$; 	

\item [\rm (303)]  $C_{54}(2,3,9,15,16,18,20,27)$, $C_{54}(3,4,9,14,15,18,22,27)$, $C_{54}(3,8,9,10,15,18,26,27)$; 	

\item [\rm (304)]  $C_{54}(2,3,9,15,16,20,21,24)$, $C_{54}(3,4,9,14,15,21,22,24)$, $C_{54}(3,8,9,10,15,21,24,26)$; 	

\item [\rm (305)]  $C_{54}(2,3,9,15,16,20,21,27)$, $C_{54}(3,4,9,14,15,21,22,27)$, $C_{54}(3,8,9,10,15,21,26,27)$; 	

\item [\rm (306)]  $C_{54}(2,3,9,15,16,20,24,27)$, $C_{54}(3,4,9,14,15,22,24,27)$, $C_{54}(3,8,9,10,15,24,26,27)$; 	

\item [\rm (307)]  $C_{54}(2,3,9,16,18,20,21,24)$, $C_{54}(3,4,9,14,18,21,22,24)$, $C_{54}(3,8,9,10,18,21,24,26)$; 	

\item [\rm (308)]  $C_{54}(2,3,9,16,18,20,21,27)$, $C_{54}(3,4,9,14,18,21,22,24)$, $C_{54}(3,8,9,10,18,21,26,27)$; 	

\item [\rm (309)]  $C_{54}(2,3,9,16,18,20,24,27)$, $C_{54}(3,4,9,14,18,22,24,27)$, $C_{54}(3,8,9,10,18,24,26,27)$; 	

\item [\rm (310)]  $C_{54}(2,3,9,16,20,21,24,27)$, $C_{54}(3,4,9,14,21,22,24,27)$, $C_{54}(3,8,9,10,21,24,26,27)$; 	

\item [\rm (311)]  $C_{54}(2,3,12,15,16,18,20,21)$, $C_{54}(3,4,12,14,15,18,21,22)$, $C_{54}(3,8,10,12,15,18,21,26)$; 	

\item [\rm (312)]  $C_{54}(2,3,12,15,16,18,20,24)$, $C_{54}(3,4,12,14,15,18,22,24)$, $C_{54}(3,8,10,12,15,18,24,26)$; 	

\item [\rm (313)]  $C_{54}(2,3,12,15,16,18,20,27)$, $C_{54}(3,4,12,14,15,18,22,27)$, $C_{54}(3,8,10,12,15,18,26,27)$; 	

\item [\rm (314)]  $C_{54}(2,3,12,15,16,20,21,24)$, $C_{54}(3,4,12,14,15,21,22,24)$, $C_{54}(3,8,10,12,15,21,24,26)$; 	

\item [\rm (315)]  $C_{54}(2,3,12,15,16,20,21,27)$, $C_{54}(3,4,12,14,15,21,22,27)$, $C_{54}(3,8,10,12,15,21,26,27)$; 	

\item [\rm (316)]  $C_{54}(2,3,12,15,16,20,24,27)$, $C_{54}(3,4,12,14,15,22,24,27)$, $C_{54}(3,8,10,12,15,24,26,27)$; 	

\item [\rm (317)]  $C_{54}(2,3,12,16,18,20,21,24)$, $C_{54}(3,4,12,14,18,21,22,24)$, $C_{54}(3,8,10,12,18,21,24,26)$; 	

\item [\rm (318)]  $C_{54}(2,3,12,16,18,20,21,27)$, $C_{54}(3,4,12,14,18,21,22,27)$, $C_{54}(3,8,10,12,18,21,26,27)$; 	

\item [\rm (319)]  $C_{54}(2,3,12,16,18,20,24,27)$, $C_{54}(3,4,12,14,18,22,24,27)$, $C_{54}(3,8,10,12,18,24,26,27)$; 	

\item [\rm (320)]  $C_{54}(2,3,12,16,20,21,24,27)$, $C_{54}(3,4,12,14,21,22,24,27)$, $C_{54}(3,8,10,12,21,24,26,27)$; 	

\item [\rm (321)]  $C_{54}(2,3,15,16,18,20,21,24)$, $C_{54}(3,4,14,15,18,21,22,24)$, $C_{54}(3,8,10,15,18,21,24,26)$; 	

\item [\rm (322)]  $C_{54}(2,3,15,16,18,20,21,27)$, $C_{54}(3,4,14,15,18,21,22,27)$, $C_{54}(3,8,10,15,18,21,26,27)$; 	

\item [\rm (323)]  $C_{54}(2,3,15,16,18,20,24,27)$, $C_{54}(3,4,14,15,18,22,24,27)$, $C_{54}(3,8,10,15,18,24,26,27)$; 	

\item [\rm (324)]  $C_{54}(2,3,15,16,20,21,24,27)$, $C_{54}(3,4,14,15,21,22,24,27)$, $C_{54}(3,8,10,15,21,24,26,27)$; 	

\item [\rm (325)]  $C_{54}(2,3,16,18,20,21,24,27)$, $C_{54}(3,4,14,18,21,22,24,27)$, $C_{54}(3,8,10,18,21,24,26,27)$; 	

\item [\rm (326)]  $C_{54}(2,6,9,12,15,16,18,20)$, $C_{54}(4,6,9,12,14,15,18,22)$, $C_{54}(6,8,9,10,12,15,18,26)$; 	

\item [\rm (327)]  $C_{54}(2,6,9,12,15,16,20,21)$, $C_{54}(4,6,9,12,14,15,21,22)$, $C_{54}(6,8,9,10,12,15,21,26)$; 	

\item [\rm (328)]  $C_{54}(2,6,9,12,15,16,20,24)$, $C_{54}(4,6,9,12,14,15,22,24)$, $C_{54}(6,8,9,10,12,15,24,26)$; 	

\item [\rm (329)]  $C_{54}(2,6,9,12,15,16,20,27)$, $C_{54}(4,6,9,12,14,15,22,27)$, $C_{54}(6,8,9,10,12,15,26,27)$; 	

\item [\rm (330)]  $C_{54}(2,6,9,12,18,16,20,21)$, $C_{54}(4,6,9,12,14,18,21,22)$, $C_{54}(6,8,9,10,12,18,21,26)$; 	

\item [\rm (331)]  $C_{54}(2,6,9,12,18,16,20,24)$, $C_{54}(4,6,9,12,14,18,22,24)$, $C_{54}(6,8,9,10,12,18,24,26)$; 	

\item [\rm (332)]  $C_{54}(2,6,9,12,18,16,20,27)$, $C_{54}(4,6,9,12,14,18,22,27)$, $C_{54}(6,8,9,10,12,18,26,27)$; 	

\item [\rm (333)]  $C_{54}(2,6,9,12,16,20,21,24)$, $C_{54}(4,6,9,12,14,21,22,24)$, $C_{54}(6,8,9,10,12,21,24,26)$; 	

\item [\rm (334)]  $C_{54}(2,6,9,12,16,20,21,27)$, $C_{54}(4,6,9,12,14,21,22,27)$, $C_{54}(6,8,9,10,12,21,26,27)$; 	

\item [\rm (335)]  $C_{54}(2,6,9,12,16,20,24,27)$, $C_{54}(4,6,9,12,14,22,24,27)$, $C_{54}(6,8,9,10,12,24,26,27)$; 	

\item [\rm (336)]  $C_{54}(2,6,9,15,16,18,20,21)$, $C_{54}(4,6,9,14,15,18,21,22)$, $C_{54}(6,8,9,10,15,18,21,26)$; 	

\item [\rm (337)]  $C_{54}(2,6,9,15,16,18,20,24)$, $C_{54}(4,6,9,14,15,18,22,24)$, $C_{54}(6,8,9,10,15,18,24,26)$; 	

\item [\rm (338)]  $C_{54}(2,6,9,15,16,18,20,27)$, $C_{54}(4,6,9,14,15,18,22,27)$, $C_{54}(6,8,9,10,15,18,26,27)$; 	

\item [\rm (339)]  $C_{54}(2,6,9,15,16,20,21,24)$, $C_{54}(4,6,9,14,15,21,22,24)$, $C_{54}(6,8,9,10,15,21,24,26)$; 	

\item [\rm (340)]  $C_{54}(2,6,9,15,16,20,21,27)$, $C_{54}(4,6,9,14,15,21,22,27)$, $C_{54}(6,8,9,10,15,21,26,27)$; 	

\item [\rm (341)]  $C_{54}(2,6,9,15,16,20,24,27)$, $C_{54}(4,6,9,14,15,22,24,27)$, $C_{54}(6,8,9,10,15,24,26,27)$; 	

\item [\rm (342)]  $C_{54}(2,6,9,16,18,20,21,24)$, $C_{54}(4,6,9,14,18,21,22,24)$, $C_{54}(6,8,9,10,18,21,24,26)$; 	

\item [\rm (343)]  $C_{54}(2,6,9,16,18,20,21,27)$, $C_{54}(4,6,9,14,18,21,22,27)$, $C_{54}(6,8,9,10,18,21,26,27)$; 	

\item [\rm (344)]  $C_{54}(2,6,9,16,18,20,24,27)$, $C_{54}(4,6,9,14,18,22,24,27)$, $C_{54}(6,8,9,10,18,24,26,27)$; 	

\item [\rm (345)]  $C_{54}(2,6,9,16,20,21,24,27)$, $C_{54}(4,6,9,14,21,22,24,27)$, $C_{54}(6,8,9,10,21,24,26,27)$; 	

\item [\rm (346)]  $C_{54}(2,6,12,15,16,18,20,21)$, $C_{54}(4,6,12,14,15,18,21,22)$, $C_{54}(6,8,10,12,15,18,21,26)$; 	

\item [\rm (347)]  $C_{54}(2,6,12,15,16,18,20,24)$, $C_{54}(4,6,12,14,15,18,22,24)$, $C_{54}(6,8,10,12,15,18,24,26)$; 	

\item [\rm (348)]  $C_{54}(2,6,12,15,16,18,20,27)$, $C_{54}(4,6,12,14,15,18,22,27)$, $C_{54}(6,8,10,12,15,18,26,27)$; 	

\item [\rm (349)]  $C_{54}(2,6,12,15,16,20,21,24)$, $C_{54}(4,6,12,14,15,21,22,24)$, $C_{54}(6,8,10,12,15,21,24,26)$; 	

\item [\rm (350)]  $C_{54}(2,6,12,15,16,20,21,27)$, $C_{54}(4,6,12,14,15,21,22,27)$, $C_{54}(6,8,10,12,15,21,26,27)$; 	

\item [\rm (351)]  $C_{54}(2,6,12,15,16,20,24,27)$, $C_{54}(4,6,12,14,15,22,24,27)$, $C_{54}(6,8,10,12,15,24,26,27)$; 	

\item [\rm (352)]  $C_{54}(2,6,12,16,18,20,21,24)$, $C_{54}(4,6,12,14,18,21,22,24)$, $C_{54}(6,8,10,12,18,21,24,26)$; 	

\item [\rm (353)]  $C_{54}(2,6,12,16,18,20,21,27)$, $C_{54}(4,6,12,14,18,21,22,27)$, $C_{54}(6,8,10,12,18,21,26,27)$; 	

\item [\rm (354)]  $C_{54}(2,6,12,16,18,20,24,27)$, $C_{54}(4,6,12,14,18,22,24,27)$, $C_{54}(6,8,10,12,18,24,26,27)$; 	

\item [\rm (355)]  $C_{54}(2,6,12,16,20,21,24,27)$, $C_{54}(4,6,12,14,21,22,24,27)$, $C_{54}(6,8,10,12,21,24,26,27)$; 	

\item [\rm (356)]  $C_{54}(2,6,15,16,18,20,21,24)$, $C_{54}(4,6,14,15,18,21,22,24)$, $C_{54}(6,8,10,15,18,21,24,26)$; 	

\item [\rm (357)]  $C_{54}(2,6,15,16,18,20,21,27)$, $C_{54}(4,6,14,15,18,21,22,27)$, $C_{54}(6,8,10,15,18,21,26,27)$; 	

\item [\rm (358)]  $C_{54}(2,6,15,16,18,20,24,27)$, $C_{54}(4,6,14,15,18,22,24,27)$, $C_{54}(6,8,10,15,18,24,26,27)$; 	

\item [\rm (359)]  $C_{54}(2,6,15,16,20,21,24,27)$, $C_{54}(4,6,14,15,21,22,24,27)$, $C_{54}(6,8,10,15,21,24,26,27)$; 	

\item [\rm (360)]  $C_{54}(2,6,16,18,20,21,24,27)$, $C_{54}(4,6,14,18,21,22,24,27)$, $C_{54}(6,8,10,18,21,24,26,27)$; 	

\item [\rm (361)]  $C_{54}(2,9,12,15,16,18,20,21)$, $C_{54}(4,9,12,14,15,18,21,22)$, $C_{54}(8,9,10,12,15,18,21,26)$; 	

\item [\rm (362)]  $C_{54}(2,9,12,15,16,18,20,24)$, $C_{54}(4,9,12,14,15,18,22,24)$, $C_{54}(8,9,10,12,15,18,24,26)$; 	

\item [\rm (363)]  $C_{54}(2,9,12,15,16,18,20,27)$, $C_{54}(4,9,12,14,15,18,22,27)$, $C_{54}(8,9,10,12,15,18,26,27)$; 	

\item [\rm (364)]  $C_{54}(2,9,12,15,16,20,21,24)$, $C_{54}(4,9,12,14,15,21,22,24)$, $C_{54}(8,9,10,12,15,21,24,26)$; 	

\item [\rm (365)]  $C_{54}(2,9,12,15,16,20,21,27)$, $C_{54}(4,9,12,14,15,21,22,27)$, $C_{54}(8,9,10,12,15,21,26,27)$; 	

\item [\rm (366)]  $C_{54}(2,9,12,15,16,20,24,27)$, $C_{54}(4,9,12,14,15,22,24,27)$, $C_{54}(8,9,10,12,15,24,26,27)$; 	

\item [\rm (367)]  $C_{54}(2,9,12,16,18,20,21,24)$, $C_{54}(4,9,12,14,18,21,22,24)$, $C_{54}(8,9,10,12,18,21,24,26)$; 	

\item [\rm (368)]  $C_{54}(2,9,12,16,18,20,21,27)$, $C_{54}(4,9,12,14,18,21,22,27)$, $C_{54}(8,9,10,12,18,21,26,27)$; 	

\item [\rm (369)]  $C_{54}(2,9,12,16,18,20,24,27)$, $C_{54}(4,9,12,14,18,22,24,27)$, $C_{54}(8,9,10,12,18,24,26,27)$; 	

\item [\rm (370)]  $C_{54}(2,9,12,16,20,21,24,27)$, $C_{54}(4,9,12,14,21,22,24,27)$, $C_{54}(8,9,10,12,21,24,26,27)$; 	

\item [\rm (371)]  $C_{54}(2,9,15,16,18,20,21,24)$, $C_{54}(4,9,14,15,18,21,22,24)$, $C_{54}(8,9,10,15,18,21,24,26)$; 	

\item [\rm (372)]  $C_{54}(2,9,15,16,18,20,21,27)$, $C_{54}(4,9,14,15,18,21,22,27)$, $C_{54}(8,9,10,15,18,21,26,27)$; 	

\item [\rm (373)]  $C_{54}(2,9,15,16,18,20,24,27)$, $C_{54}(4,9,14,15,18,22,24,27)$, $C_{54}(8,9,10,15,18,24,26,27)$; 	

\item [\rm (374)]  $C_{54}(2,9,15,16,20,21,24,27)$, $C_{54}(4,9,14,15,21,22,24,27)$, $C_{54}(8,9,10,15,21,24,26,27)$; 	

\item [\rm (375)]  $C_{54}(2,9,16,18,20,21,24,27)$, $C_{54}(4,9,14,18,21,22,24,27)$, $C_{54}(8,9,10,18,21,24,26,27)$; 	

\item [\rm (376)]  $C_{54}(2,12,15,16,18,20,21,24)$, $C_{54}(4,12,14,15,18,21,22,24)$, $C_{54}(8,10,12,15,18,21,24,26)$; 

\item [\rm (377)]  $C_{54}(2,12,15,16,18,20,21,27)$, $C_{54}(4,12,14,15,18,21,22,27)$, $C_{54}(8,10,12,15,18,21,26,27)$; 

\item [\rm (378)] $C_{54}(2,12,15,16,18,20,24,27)$, $C_{54}(4,12,14,15,18,22,24,27)$, $C_{54}(8,10,12,15,18,24,26,27)$; 

\item [\rm (379)]  $C_{54}(2,12,15,16,20,21,24,27)$, $C_{54}(4,12,14,15,21,22,24,27)$, $C_{54}(8,10,12,15,21,24,26,27)$; 

\item [\rm (380)]  $C_{54}(2,12,16,18,20,21,24,27)$, $C_{54}(4,12,14,18,21,22,24,27)$, $C_{54}(8,10,12,18,21,24,26,27)$; 

\item [\rm (381)] $C_{54}(2,15,16,18,20,21,24,27)$, $C_{54}(4,14,15,18,21,22,24,27)$, $C_{54}(8,10,15,18,21,24,26,27)$; 	

\item [\rm (382)]   $C_{54}(2,3,6,9,12,15,16,18,20)$, $C_{54}(3,4,6,9,12,14,15,18,22)$, $C_{54}(3,6,8,9,10,12,15,18,26)$; 	

\item [\rm (383)]   $C_{54}(2,3,6,9,12,15,16,20,21)$, $C_{54}(3,4,6,9,12,14,15,21,22)$, $C_{54}(3,6,8,9,10,12,15,21,26)$; 	

\item [\rm (384)]   $C_{54}(2,3,6,9,12,15,16,20,24)$, $C_{54}(3,4,6,9,12,14,15,22,24)$, $C_{54}(3,6,8,9,10,12,15,24,26)$; 	

\item [\rm (385)]   $C_{54}(2,3,6,9,12,15,16,20,27)$, $C_{54}(3,4,6,9,12,14,15,22,27)$, $C_{54}(3,6,8,9,10,12,15,26,27)$; 	

\item [\rm (386)]  $C_{54}(2,3,6,9,12,16,18,20,21)$, $C_{54}(3,4,6,9,12,14,18,21,22)$, $C_{54}(3,6,8,9,10,12,18,21,26)$; 	

\item [\rm (387)]  $C_{54}(2,3,6,9,12,16,18,20,24)$, $C_{54}(3,4,6,9,12,14,18,22,24)$, $C_{54}(3,6,8,9,10,12,18,24,26)$; 	

\item [\rm (388)]  $C_{54}(2,3,6,9,12,16,18,20,27)$, $C_{54}(3,4,6,9,12,14,18,22,27)$, $C_{54}(3,6,8,9,10,12,18,26,27)$; 	

\item [\rm (389)]  $C_{54}(2,3,6,9,12,16,20,21,24)$, $C_{54}(3,4,6,9,12,14,21,22,24)$, $C_{54}(3,6,8,9,10,12,21,24,26)$; 	

\item [\rm (390)]  $C_{54}(2,3,6,9,12,16,20,21,27)$, $C_{54}(3,4,6,9,12,14,21,22,27)$, $C_{54}(3,6,8,9,10,12,21,26,27)$; 	

\item [\rm (391)]  $C_{54}(2,3,6,9,12,16,20,24,27)$, $C_{54}(3,4,6,9,12,14,22,24,27)$, $C_{54}(3,6,8,9,10,12,24,26,27)$; 	

\item [\rm (392)]  $C_{54}(2,3,6,9,15,16,18,20,21)$, $C_{54}(3,4,6,9,14,15,18,21,22)$, $C_{54}(3,6,8,9,10,15,18,21,26)$; 	

\item [\rm (393)]  $C_{54}(2,3,6,9,15,16,18,20,24)$, $C_{54}(3,4,6,9,14,15,18,22,24)$, $C_{54}(3,6,8,9,10,15,18,24,26)$; 	

\item [\rm (394)]  $C_{54}(2,3,6,9,15,16,18,20,27)$, $C_{54}(3,4,6,9,14,15,18,22,27)$, $C_{54}(3,6,8,9,10,15,18,26,27)$; 	

\item [\rm (395)]  $C_{54}(2,3,6,9,15,16,20,21,24)$, $C_{54}(3,4,6,9,14,15,21,22,24)$, $C_{54}(3,6,8,9,10,15,21,24,26)$; 	

\item [\rm (396)]  $C_{54}(2,3,6,9,15,16,20,21,27)$, $C_{54}(3,4,6,9,14,15,21,22,27)$, $C_{54}(3,6,8,9,10,15,21,26,27)$; 	

\item [\rm (397)]  $C_{54}(2,3,6,9,15,16,20,24,27)$, $C_{54}(3,4,6,9,14,15,22,24,27)$, $C_{54}(3,6,8,9,10,15,24,26,27)$; 	

\item [\rm (398)]  $C_{54}(2,3,6,9,16,18,20,21,24)$, $C_{54}(3,4,6,9,14,18,21,22,24)$, $C_{54}(3,6,8,9,10,18,21,24,26)$; 	

\item [\rm (399)]  $C_{54}(2,3,6,9,16,18,20,21,27)$, $C_{54}(3,4,6,9,14,18,21,22,27)$, $C_{54}(3,6,8,9,10,18,21,26,27)$; 	

\item [\rm (400)]  $C_{54}(2,3,6,9,16,18,20,24,27)$, $C_{54}(3,4,6,9,14,18,22,24,27)$, $C_{54}(3,6,8,9,10,18,24,26,27)$; 	

\item [\rm (401)]  $C_{54}(2,3,6,9,16,20,21,24,27)$, $C_{54}(3,4,6,9,14,21,22,24,27)$, $C_{54}(3,6,8,9,10,21,24,26,27)$; 	

\item [\rm (402)]  $C_{54}(2,3,6,12,15,16,18,20,21)$, $C_{54}(3,4,6,12,14,15,18,21,22)$, $C_{54}(3,6,8,10,12,15,18,21,26)$; 	

\item [\rm (403)] $C_{54}(2,3,6,12,15,16,18,20,24)$, $C_{54}(3,4,6,12,14,15,18,22,24)$, $C_{54}(3,6,8,10,12,15,18,24,26)$; 	

\item [\rm (404)]  $C_{54}(2,3,6,12,15,16,18,20,27)$, $C_{54}(3,4,6,12,14,15,18,22,27)$, $C_{54}(3,6,8,10,12,15,18,26,27)$; 	

\item [\rm (405)]  $C_{54}(2,3,6,12,15,16,20,21,24)$, $C_{54}(3,4,6,12,14,15,21,22,24)$, $C_{54}(3,6,8,10,12,15,21,24,26)$; 	

\item [\rm (406)]  $C_{54}(2,3,6,12,15,16,20,21,27)$, $C_{54}(3,4,6,12,14,15,21,22,27)$, $C_{54}(3,6,8,10,12,15,21,26,27)$; 	

\item [\rm (407)]  $C_{54}(2,3,6,12,15,16,20,24,27)$, $C_{54}(3,4,6,12,14,15,22,24,27)$, $C_{54}(3,6,8,10,12,15,24,26,27)$; 	

\item [\rm (408)]  $C_{54}(2,3,6,12,16,18,20,21,24)$, $C_{54}(3,4,6,12,14,18,21,22,24)$, $C_{54}(3,6,8,10,12,18,21,24,26)$; 	

\item [\rm (409)]  $C_{54}(2,3,6,12,16,18,20,21,27)$, $C_{54}(3,4,6,12,14,18,21,22,27)$, $C_{54}(3,6,8,10,12,18,21,26,27)$; 	

\item [\rm (410)]  $C_{54}(2,3,6,12,16,18,20,24,27)$, $C_{54}(3,4,6,12,14,18,22,24,27)$, $C_{54}(3,6,8,10,12,18,24,26,27)$; 	

\item [\rm (411)]  $C_{54}(2,3,6,12,16,20,21,24,27)$, $C_{54}(3,4,6,12,14,21,22,24,27)$, $C_{54}(3,6,8,10,12,21,24,26,27)$; 	

\item [\rm (412)]  $C_{54}(2,3,6,15,16,18,20,21,24)$, $C_{54}(3,4,6,14,15,18,21,22,24)$, $C_{54}(3,6,8,10,15,18,21,24,26)$; 	

\item [\rm (413)]  $C_{54}(2,3,6,15,16,18,20,21,27)$, $C_{54}(3,4,6,14,15,18,21,22,27)$, $C_{54}(3,6,8,10,15,18,21,26,27)$; 	

\item [\rm (414)]  $C_{54}(2,3,6,15,16,18,20,24,27)$, $C_{54}(3,4,6,14,15,18,22,24,27)$, $C_{54}(3,6,8,10,15,18,24,26,27)$; 	

\item [\rm (415)]  $C_{54}(2,3,6,15,16,20,21,24,27)$, $C_{54}(3,4,6,14,15,21,22,24,27)$, $C_{54}(3,6,8,10,15,21,24,26,27)$; 	

\item [\rm (416)]  $C_{54}(2,3,6,16,18,20,21,24,27)$, $C_{54}(3,4,6,14,18,21,22,24,27)$, $C_{54}(3,6,8,10,18,21,24,26,27)$; 	

\item [\rm (417)]  $C_{54}(2,3,9,12,15,16,18,20,21)$, $C_{54}(3,4,9,12,14,15,18,21,22)$, $C_{54}(3,8,9,10,12,15,18,21,26)$; 	

\item [\rm (418)]  $C_{54}(2,3,9,12,15,16,18,20,24)$, $C_{54}(3,4,9,12,14,15,18,22,24)$, $C_{54}(3,8,9,10,12,15,18,24,26)$; 	

\item [\rm (419)]  $C_{54}(2,3,9,12,15,16,18,20,27)$, $C_{54}(3,4,9,12,14,15,18,22,27)$, $C_{54}(3,8,9,10,12,15,18,26,27)$; 	

\item [\rm (420)]  $C_{54}(2,3,9,12,15,16,20,21,24)$, $C_{54}(3,4,9,12,14,15,21,22,24)$, $C_{54}(3,8,9,10,12,15,21,24,26)$; 	

\item [\rm (421)]  $C_{54}(2,3,9,12,15,16,20,21,27)$, $C_{54}(3,4,9,12,14,15,21,22,27)$, $C_{54}(3,8,9,10,12,15,21,26,27)$; 	

\item [\rm (422)]  $C_{54}(2,3,9,12,15,16,20,24,27)$, $C_{54}(3,4,9,12,14,15,22,24,27)$, $C_{54}(3,8,9,10,12,15,24,26,27)$; 	

\item [\rm (423)]  $C_{54}(2,3,9,12,16,18,20,21,24)$, $C_{54}(3,4,9,12,14,18,21,22,24)$, $C_{54}(3,8,9,10,12,18,21,24,26)$; 	

\item [\rm (424)]  $C_{54}(2,3,9,12,16,18,20,21,27)$, $C_{54}(3,4,9,12,14,18,21,22,27)$, $C_{54}(3,8,9,10,12,18,21,26,27)$; 	

\item [\rm (425)]  $C_{54}(2,3,9,12,16,18,20,24,27)$, $C_{54}(3,4,9,12,14,18,22,24,27)$, $C_{54}(3,8,9,10,12,18,24,26,27)$; 	

\item [\rm (426)]  $C_{54}(2,3,9,12,16,20,21,24,27)$, $C_{54}(3,4,9,12,14,21,22,24,27)$, $C_{54}(3,8,9,10,12,21,24,26,27)$; 	

\item [\rm (427)]  $C_{54}(2,3,9,15,16,18,20,21,24)$, $C_{54}(3,4,9,14,15,18,21,22,24)$, $C_{54}(3,8,9,10,15,18,21,24,26)$; 	

\item [\rm (428)]  $C_{54}(2,3,9,15,16,18,20,21,27)$, $C_{54}(3,4,9,14,15,18,21,22,27)$, $C_{54}(3,8,9,10,15,18,21,26,27)$; 	

\item [\rm (429)]  $C_{54}(2,3,9,15,16,18,20,24,27)$, $C_{54}(3,4,9,14,15,18,22,24,27)$, $C_{54}(3,8,9,10,15,18,24,26,27)$; 	

\item [\rm (430)]  $C_{54}(2,3,9,15,16,20,21,24,27)$, $C_{54}(3,4,9,14,15,21,22,24,27)$, $C_{54}(3,8,9,10,15,21,24,26,27)$; 	

\item [\rm (431)]  $C_{54}(2,3,9,16,18,20,21,24,27)$, $C_{54}(3,4,9,14,18,21,22,24,27)$, $C_{54}(3,8,9,10,18,21,24,26,27)$; 	

\item [\rm (432)]  $C_{54}(2,3,12,15,16,18,20,21,24)$, $C_{54}(3,4,12,14,15,18,21,22,24)$, 

\hfill $C_{54}(3,8,10,12,15,18,21,24,26)$; 	

\item [\rm (433)]  $C_{54}(2,3,12,15,16,18,20,21,27)$, $C_{54}(3,4,12,14,15,18,21,22,27)$, 

\hfill $C_{54}(3,8,10,12,15,18,21,26,27)$; 	

\item [\rm (434)]  $C_{54}(2,3,12,15,16,18,20,24,27)$, $C_{54}(3,4,12,14,15,18,22,24,27)$, 

\hfill $C_{54}(3,8,10,12,15,18,24,26,27)$; 	

\item [\rm (435)]  $C_{54}(2,3,12,15,16,20,21,24,27)$, $C_{54}(3,4,12,14,15,21,22,24,27)$, 

\hfill $C_{54}(3,8,10,12,15,21,24,26,27)$; 	

\item [\rm (436)]  $C_{54}(2,3,12,16,18,20,21,24,27)$, $C_{54}(3,4,12,14,18,21,22,24,27)$, 

\hfill $C_{54}(3,8,10,12,18,21,24,26,27)$; 	

\item [\rm (437)]  $C_{54}(2,3,15,16,18,20,21,24,27)$, $C_{54}(3,4,14,15,18,21,22,24,27)$, 

\hfill $C_{54}(3,8,10,15,18,21,24,26,27)$; 	

\item [\rm (438)]  $C_{54}(2,6,9,12,15,16,18,20,21)$, $C_{54}(4,6,9,12,14,15,18,21,22)$, 

\hfill $C_{54}(6,8,9,10,12,15,18,21,26)$; 	

\item [\rm (439)]  $C_{54}(2,6,9,12,15,16,18,20,24)$, $C_{54}(4,6,9,12,14,15,18,22,24)$, 

\hfill $C_{54}(6,8,9,10,12,15,18,24,26)$; 	

\item [\rm (440)]  $C_{54}(2,6,9,12,15,16,18,20,27)$, $C_{54}(4,6,9,12,14,15,18,22,27)$, 

\hfill $C_{54}(6,8,9,10,12,15,18,26,27)$; 	

\item [\rm (441)]  $C_{54}(2,6,9,12,15,16,20,21,24)$, $C_{54}(4,6,9,12,14,15,21,22,24)$, 

\hfill $C_{54}(6,8,9,10,12,15,21,24,26)$; 	

\item [\rm (442)]  $C_{54}(2,6,9,12,15,16,20,21,27)$, $C_{54}(4,6,9,12,14,15,21,22,27)$, 

\hfill $C_{54}(6,8,9,10,12,15,21,26,27)$; 	

\item [\rm (443)]  $C_{54}(2,6,9,12,15,16,20,24,27)$, $C_{54}(4,6,9,12,14,15,22,24,27)$, 

\hfill $C_{54}(6,8,9,10,12,15,24,26,27)$; 	

\item [\rm (444)]  $C_{54}(2,6,9,12,16,18,20,21,24)$, $C_{54}(4,6,9,12,14,18,21,22,24)$, 

\hfill $C_{54}(6,8,9,10,12,18,21,24,26)$; 	

\item [\rm (445)]  $C_{54}(2,6,9,12,16,18,20,21,27)$, $C_{54}(4,6,9,12,14,18,21,22,27)$, 

\hfill $C_{54}(6,8,9,10,12,18,21,26,27)$; 	

\item [\rm (446)]  $C_{54}(2,6,9,12,16,18,20,24,27)$, $C_{54}(4,6,9,12,14,18,22,24,27)$, 

\hfill $C_{54}(6,8,9,10,12,18,24,26,27)$; 	

\item [\rm (447)]  $C_{54}(2,6,9,12,16,20,21,24,27)$, $C_{54}(4,6,9,12,14,21,22,24,27)$, 

\hfill $C_{54}(6,8,9,10,12,21,24,26,27)$; 	

\item [\rm (448)]  $C_{54}(2,6,9,15,16,18,20,21,24)$, $C_{54}(4,6,9,14,15,18,21,22,24)$, 

\hfill $C_{54}(6,8,9,10,15,18,21,24,26)$; 	

\item [\rm (449)]  $C_{54}(2,6,9,15,16,18,20,21,27)$, $C_{54}(4,6,9,14,15,18,21,22,27)$, 

\hfill $C_{54}(6,8,9,10,15,18,21,26,27)$; 	

\item [\rm (450)]  $C_{54}(2,6,9,15,16,18,20,24,27)$, $C_{54}(4,6,9,14,15,18,22,24,27)$, 

\hfill $C_{54}(6,8,9,10,15,18,24,26,27)$; 	

\item [\rm (451)]  $C_{54}(2,6,9,15,16,20,21,24,27)$, $C_{54}(4,6,9,14,15,21,22,24,27)$, 

\hfill $C_{54}(6,8,9,10,15,21,24,26,27)$; 	

\item [\rm (452)]  $C_{54}(2,6,9,16,18,20,21,24,27)$, $C_{54}(4,6,9,14,18,21,22,24,27)$, 

\hfill $C_{54}(6,8,9,10,18,21,24,26,27)$; 	

\item [\rm (453)]  $C_{54}(2,6,12,15,16,18,20,21,24)$, $C_{54}(4,6,12,14,15,18,21,22,24)$, 

\hfill $C_{54}(6,8,10,12,15,18,21,24,26)$; 	

\item [\rm (454)]  $C_{54}(2,6,12,15,16,18,20,21,27)$, $C_{54}(4,6,12,14,15,18,21,22,27)$, 

\hfill $C_{54}(6,8,10,12,15,18,21,26,27)$; 	

\item [\rm (455)]  $C_{54}(2,6,12,15,16,18,20,24,27)$, $C_{54}(4,6,12,14,15,18,22,24,27)$, 

\hfill $C_{54}(6,8,10,12,15,18,24,26,27)$; 	

\item [\rm (456)]  $C_{54}(2,6,12,15,16,20,21,24,27)$, $C_{54}(4,6,12,14,15,21,22,24,27)$, 

\hfill $C_{54}(6,8,10,12,15,21,24,26,27)$; 	

\item [\rm (457)]  $C_{54}(2,6,12,16,18,20,21,24,27)$, $C_{54}(4,6,12,14,18,21,22,24,27)$, 

\hfill $C_{54}(6,8,10,12,18,21,24,26,27)$; 	

\item [\rm (458)]  $C_{54}(2,6,15,16,18,20,21,24,27)$, $C_{54}(4,6,14,15,18,21,22,24,27)$, 

\hfill $C_{54}(6,8,10,15,18,21,24,26,27)$; 	

\item [\rm (459)]  $C_{54}(2,9,12,15,16,18,20,21,24)$, $C_{54}(4,9,12,14,15,18,21,22,24)$, 

\hfill $C_{54}(8,9,10,12,15,18,21,24,26)$; 	

\item [\rm (460)]  $C_{54}(2,9,12,15,16,18,20,21,27)$, $C_{54}(4,9,12,14,15,18,21,22,27)$, 

\hfill $C_{54}(8,9,10,12,15,18,21,26,27)$; 	

\item [\rm (461)]  $C_{54}(2,9,12,15,16,18,20,24,27)$, $C_{54}(4,9,12,14,15,18,22,24,27)$, 

\hfill $C_{54}(8,9,10,12,15,18,24,26,27)$; 	

\item [\rm (462)]  $C_{54}(2,9,12,15,16,20,21,24,27)$, $C_{54}(4,9,12,14,15,21,22,24,27)$, 

\hfill $C_{54}(8,9,10,12,15,21,24,26,27)$; 	

\item [\rm (463)]  $C_{54}(2,9,12,16,18,20,21,24,27)$, $C_{54}(4,9,12,14,18,21,22,24,27)$, 

\hfill $C_{54}(8,9,10,12,18,21,24,26,27)$; 	

\item [\rm (464)]  $C_{54}(2,9,15,16,18,20,21,24,27)$, $C_{54}(4,9,14,15,18,21,22,24,27)$, 

\hfill $C_{54}(8,9,10,15,18,21,24,26,27)$; 	

\item [\rm (465)]  $C_{54}(2,12,15,16,18,20,21,24,27)$, $C_{54}(4,12,14,15,18,21,22,24,27)$, 

\hfill $C_{54}(8,10,12,15,18,21,24,26,27)$; 	

\item [\rm (466)]  $C_{54}(2,3,6,9,12,15,16,18,20,21)$, $C_{54}(3,4,6,9,12,14,15,18,21,22)$, 

\hfill $C_{54}(3,6,8,9,10,12,15,18,21,26)$; 	

\item [\rm (467)]  $C_{54}(2,3,6,9,12,15,16,18,20,24)$, $C_{54}(3,4,6,9,12,14,15,18,22,24)$, 

\hfill $C_{54}(3,6,8,9,10,12,15,18,24,26)$; 	

\item [\rm (468)]  $C_{54}(2,3,6,9,12,15,16,18,20,27)$, $C_{54}(3,4,6,9,12,14,15,18,22,27)$, 

\hfill $C_{54}(3,6,8,9,10,12,15,18,26,27)$; 	

\item [\rm (469)]   $C_{54}(2,3,6,9,12,15,16,20,21,24)$, $C_{54}(3,4,6,9,12,14,15,21,22,24)$, 

\hfill $C_{54}(3,6,8,9,10,12,15,21,24,26)$; 	

\item [\rm (470)]   $C_{54}(2,3,6,9,12,15,16,20,21,27)$, $C_{54}(3,4,6,9,12,14,15,21,22,27)$, 

\hfill $C_{54}(3,6,8,9,10,12,15,21,26,27)$; 	

\item [\rm (471)]   $C_{54}(2,3,6,9,12,15,16,20,24,27)$, $C_{54}(3,4,6,9,12,14,15,22,24,27)$, 

\hfill $C_{54}(3,6,8,9,10,12,15,24,26,27)$; 	

\item [\rm (472)]   $C_{54}(2,3,6,9,12,16,18,20,21,24)$, $C_{54}(3,4,6,9,12,14,18,21,22,24)$, 

\hfill $C_{54}(3,6,8,9,10,12,18,21,24,26)$; 	

\item [\rm (473)]   $C_{54}(2,3,6,9,12,16,18,20,21,27)$, $C_{54}(3,4,6,9,12,14,18,21,22,27)$, 

\hfill $C_{54}(3,6,8,9,10,12,18,21,26,27)$; 	

\item [\rm (474)]  $C_{54}(2,3,6,9,12,16,18,20,24,27)$, $C_{54}(3,4,6,9,12,14,18,22,24,27)$, 

\hfill $C_{54}(3,6,8,9,10,12,18,24,26,27)$; 	

\item [\rm (475)]  $C_{54}(2,3,6,9,12,16,20,21,24,27)$, $C_{54}(3,4,6,9,12,14,21,22,24,27)$, 

\hfill $C_{54}(3,6,8,9,10,12,21,24,26,27)$; 	

\item [\rm (476)]  $C_{54}(2,3,6,9,15,16,18,20,21,24)$, $C_{54}(3,4,6,9,14,15,18,21,22,24)$, 

\hfill $C_{54}(3,6,8,9,10,15,18,21,24,26)$; 	

\item [\rm (477)]  $C_{54}(2,3,6,9,15,16,18,20,21,27)$, $C_{54}(3,4,6,9,14,15,18,21,22,27)$, 

\hfill $C_{54}(3,6,8,9,10,15,18,21,26,27)$; 	

\item [\rm (478)]  $C_{54}(2,3,6,9,15,16,18,20,24,27)$, $C_{54}(3,4,6,9,14,15,18,22,24,27)$, 

\hfill $C_{54}(3,6,8,9,10,15,18,24,26,27)$; 	

\item [\rm (479)]  $C_{54}(2,3,6,9,15,16,20,21,24,27)$, $C_{54}(3,4,6,9,14,15,21,22,24,27)$, 

\hfill $C_{54}(3,6,8,9,10,15,21,24,26,27)$; 	

\item [\rm (480)]  $C_{54}(2,3,6,9,16,18,20,21,24,27)$, $C_{54}(3,4,6,9,14,18,21,22,24,27)$, 

\hfill $C_{54}(3,6,8,9,10,18,21,24,26,27)$; 	

\item [\rm (481)]  $C_{54}(2,3,6,12,15,16,18,20,21,24)$, $C_{54}(3,4,6,12,14,15,18,21,22,24)$, 

\hfill $C_{54}(3,6,8,10,12,15,18,21,24,26)$; 	

\item [\rm (482)]  $C_{54}(2,3,6,12,15,16,18,20,21,27)$, $C_{54}(3,4,6,12,14,15,18,21,22,27)$, 

\hfill $C_{54}(3,6,8,10,12,15,18,21,26,27)$; 	

\item [\rm (483)] $C_{54}(2,3,6,12,15,16,18,20,24,27)$, $C_{54}(3,4,6,12,14,15,18,22,24,27)$, 

\hfill $C_{54}(3,6,8,10,12,15,18,24,26,27)$; 	

\item [\rm (484)] $C_{54}(2,3,6,12,15,16,20,21,24,27)$, $C_{54}(3,4,6,12,14,15,21,22,24,27)$, 

\hfill $C_{54}(3,6,8,10,12,15,21,24,26,27)$; 	

\item [\rm (485)]  $C_{54}(2,3,6,12,16,18,20,21,24,27)$, $C_{54}(3,4,6,12,14,18,21,22,24,27)$, 

\hfill $C_{54}(3,6,8,10,12,18,21,24,26,27)$; 	

\item [\rm (486)]  $C_{54}(2,3,6,15,16,18,20,21,24,27)$, $C_{54}(3,4,6,14,15,18,21,22,24,27)$, 

\hfill $C_{54}(3,6,8,10,15,18,21,24,26,27)$; 	

\item [\rm (487)]  $C_{54}(2,3,9,12,15,16,18,20,21,24)$, $C_{54}(3,4,9,12,14,15,18,21,22,24)$, 

\hfill $C_{54}(3,8,9,10,12,15,18,21,24,26)$; 	

\item [\rm (488)]  $C_{54}(2,3,9,12,15,16,18,20,21,27)$, $C_{54}(3,4,9,12,14,15,18,21,22,27)$, 

\hfill $C_{54}(3,8,9,10,12,15,18,21,26,27)$; 	

\item [\rm (489)]  $C_{54}(2,3,9,12,15,16,18,20,24,27)$, $C_{54}(3,4,9,12,14,15,18,22,24,27)$, 

\hfill $C_{54}(3,8,9,10,12,15,18,24,26,27)$; 	

\item [\rm (490)]  $C_{54}(2,3,9,12,15,16,20,21,24,27)$, $C_{54}(3,4,9,12,14,15,21,22,24,27)$, 

\hfill $C_{54}(3,8,9,10,12,15,21,24,26,27)$; 	

\item [\rm (491)]  $C_{54}(2,3,9,12,16,18,20,21,24,27)$, $C_{54}(3,4,9,12,14,18,21,22,24,27)$, 

\hfill $C_{54}(3,8,9,10,12,18,21,24,26,27)$; 	

\item [\rm (492)]  $C_{54}(2,3,9,15,16,18,20,21,24,27)$, $C_{54}(3,4,9,14,15,18,21,22,24,27)$, 

\hfill $C_{54}(3,8,9,10,15,18,21,24,26,27)$; 	

\item [\rm (493)]  $C_{54}(2,3,12,15,16,18,20,21,24,27)$, $C_{54}(3,4,12,14,15,18,21,22,24,27)$, 

\hfill $C_{54}(3,8,10,12,15,18,21,24,26,27)$; 	

\item [\rm (494)]  $C_{54}(2,6,9,12,15,16,18,20,21,24)$, $C_{54}(4,6,9,12,14,15,18,21,22,24)$, 

\hfill $C_{54}(6,8,9,10,12,15,18,21,24,26)$; 	

\item [\rm (495)]  $C_{54}(2,6,9,12,15,16,18,20,21,27)$, $C_{54}(4,6,9,12,14,15,18,21,22,27)$, 

\hfill $C_{54}(6,8,9,10,12,15,18,21,26,27)$; 	

\item [\rm (496)]  $C_{54}(2,6,9,12,15,16,18,20,24,27)$, $C_{54}(4,6,9,12,14,15,18,22,24,27)$, 

\hfill $C_{54}(6,8,9,10,12,15,18,24,26,27)$; 	

\item [\rm (497)]  $C_{54}(2,6,9,12,15,16,20,21,24,27)$, $C_{54}(4,6,9,12,14,15,21,22,24,27)$, 

\hfill $C_{54}(6,8,9,10,12,15,21,24,26,27)$; 	

\item [\rm (498)]  $C_{54}(2,6,9,12,16,18,20,21,24,27)$, $C_{54}(4,6,9,12,14,18,21,22,24,27)$, 

\hfill $C_{54}(6,8,9,10,12,18,21,24,26,27)$; 	

\item [\rm (499)]  $C_{54}(2,6,9,15,16,18,20,21,24,27)$, $C_{54}(4,6,9,14,15,18,21,22,24,27)$, 

\hfill $C_{54}(6,8,9,10,15,18,21,24,26,27)$; 	

\item [\rm (500)]  $C_{54}(2,6,12,15,16,18,20,21,24,27)$, $C_{54}(4,6,12,14,15,18,21,22,24,27)$, 

\hfill $C_{54}(6,8,10,12,15,18,21,24,26,27)$; 	

\item [\rm (501)]  $C_{54}(2,9,12,15,16,18,20,21,24,27)$, $C_{54}(4,9,12,14,15,18,21,22,24,27)$, 

\hfill $C_{54}(8,9,10,12,15,18,21,24,26,27)$; 	

\item [\rm (502)]  $C_{54}(2,3,6,9,12,15,16,18,20,21,24)$, $C_{54}(3,4,6,9,12,14,15,18,21,22,24)$, 

\hfill $C_{54}(3,6,8,9,10,12,15,18,21,24,26)$; 	

\item [\rm (503)]  $C_{54}(2,3,6,9,12,15,16,18,20,21,27)$, $C_{54}(3,4,6,9,12,14,15,18,21,22,27)$, 

\hfill $C_{54}(3,6,8,9,10,12,15,18,21,26,27)$; 	

\item [\rm (504)]  $C_{54}(2,3,6,9,12,15,16,18,20,24,27)$, $C_{54}(3,4,6,9,12,14,15,18,22,24,27)$, 

\hfill $C_{54}(3,6,8,9,10,12,15,18,24,26,27)$; 	

\item [\rm (505)]  $C_{54}(2,3,6,9,12,15,16,20,21,24,27)$, $C_{54}(3,4,6,9,12,14,15,21,22,24,27)$, 

\hfill $C_{54}(3,6,8,9,10,12,15,21,24,26,27)$; 	

\item [\rm (506)]   $C_{54}(2,3,6,9,12,16,18,20,21,24,27)$, $C_{54}(3,4,6,9,12,14,18,21,22,24,27)$, 

\hfill $C_{54}(3,6,8,9,10,12,18,21,24,26,27)$; 	

\item [\rm (507)]   $C_{54}(2,3,6,9,15,16,18,20,21,24,27)$, $C_{54}(3,4,6,9,14,15,18,21,22,24,27)$, 

\hfill $C_{54}(3,6,8,9,10,15,18,21,24,26,27)$; 	

\item [\rm (508)]  $C_{54}(2,3,6,12,15,16,18,20,21,24,27)$, $C_{54}(3,4,6,12,14,15,18,21,22,24,27)$, 

\hfill $C_{54}(3,6,8,10,12,15,18,21,24,26,27)$; 	

\item [\rm (509)]  $C_{54}(2,3,9,12,15,16,18,20,21,24,27)$, $C_{54}(3,4,9,12,14,15,18,21,22,24,27)$, 

\hfill $C_{54}(3,8,9,10,12,15,18,21,24,26,27)$; 	

\item [\rm (510)]  $C_{54}(2,6,9,12,15,16,18,20,21,24,27)$, $C_{54}(6,4,9,12,14,15,18,21,22,24,27)$, 

\hfill $C_{54}(6,8,9,10,12,15,18,21,24,26,27)$; 	

\item [\rm (511)]  $C_{54}(2,3,6,9,12,15,16,18,20,21,24,27)$, $C_{54}(3,4,6,9,12,14,15,18,21,22,24,27)$, 

\hfill $C_{54}(3,6,8,9,10,12,15,18,21,24,26,27)$.

\end{enumerate}
}
 \end{prm}
 \noindent
 {\bf Solution.}\quad Here, we consider Type-1 and Type-2 isomorphisms of circulant graphs of the form $C_{54}(R)$ and so the possible values of $m > 1$ for the existence of Type-2 isomorphism w.r.t. $m$ of $C_{54}(R)$ $\ni$ $m$ is a divisor of $\gcd(n, r)$ = $\gcd(54, r)$, $m^3$ divides $n$ = 54 = $2\times 3^3$ and $r\in R$ is $m$ = 3. Also, we have $m$ = 3 = $\gcd(54, 3)$ = $\gcd(54, 15)$ = $\gcd(54, 21)$, 6 = $\gcd(54, 6)$ = $\gcd(54, 12)$, 9 = $\gcd(54, 9)$, 18 = $\gcd(54, 18)$ and 27 = $\gcd(54, 27)$. 

At first, we consider triples of circulant graphs of order 54 containing jump sizes 1,17,19, 7,11,25, 5,13 and 23 and then we consider triples of circulant graphs of order 54 containing jump sizes 2,16,20, 4,14,22, 8,10 and 26. 

\vspace{.2cm}
\noindent
{\bf Case (a)}\quad  In problem \ref{p2}, we proved that circulant graphs 

\vspace{.1cm}
(1) $C_{54}(1,3,17,19)$, $C_{54}(3,7,11,25)$ and $C_{54}(3,5,13,23)$ are isomorphic of Type-2 w.r.t. $m$ = 3. 
\\
Using remark \ref{r5} in this triple of Type-2 isomorphic circulant graphs, we obtain 511 triples of isomorphic circulant graphs as given in the problem. Also, for a given circulant graph $C_n(R)$, if all $C_n(S)$ $\ni$ $C_n(S)$ = $\theta_{n,m,t}(C_{n}(R))$ for some $t$ and $C_n(S)\in T1_n(C_n(R))$, $1 \leq t \leq \frac{n}{m}-1$, then $C_n(R)$ has no isomorphic circulant graph of Type-2  w.r.t. $m$ where $r\in R,S$ and $m > 1$ and $m^3$ are divisors of $\gcd(n, r)$ and $n$, respectively. 

Solutions to all the cases in this problem are similar, in cases of proving Type-1 isomorphism and also in the cases of Type-2 isomorphism, and to simplify our work, we present important values related to all cases related to Type-2 isomorphism in Table 1 to Table 23. In each table, in the column of `T1 or T2' corresponds to whether the triple of isomorphic circulant graphs are `Type-1 isomorphic or Type-2 isomorphic'. Even though some triples of circulant graphs are already covered in problem \ref{p1}, we consider all the triples of isomorphic circulant graphs. From these 23 tables, we found that there are 480 number of triples of Type-2 isomorphic circulant graphs w.r.t. $m$ = 3.

\medskip
\noindent
{\bf Case (b)}\quad  In problem \ref{p2}, we proved that circulant graphs 

\vspace{.1cm}
(1) $C_{54}(2,3,16,20)$, $C_{54}(3,4,14,22)$ and $C_{54}(3,8,10,26)$ are isomorphic of Type-2 w.r.t. $m$ = 3. 
\\
Using remark \ref{r5} in this triple of Type-2 isomorphic circulant graphs, we obtain 511 triples of isomorphic circulant graphs as given in the problem. Also, for a given circulant graph $C_n(R)$, if all $C_n(S)$ $\ni$ $C_n(S)$ = $\theta_{n,m,t}(C_{n}(R))$ for some $t$ and $C_n(S)\in T1_n(C_n(R))$, $1 \leq t \leq \frac{n}{m}-1$, then $C_n(R)$ has no isomorphic circulant graph of Type-2  w.r.t. $m$ where $r\in R,S$ and $m > 1$ and $m^3$ are divisors of $\gcd(n, r)$ and $n$, respectively. 

Solutions to all the cases in this problem are similar, in cases of proving Type-1 isomorphism and also in the cases of Type-2 isomorphism, and to simplify our work, we present important values related to all cases related to Type-2 isomorphism in Table 24 to Table 46. In each table, in the column of `T1 or T2' corresponds to whether the triple of isomorphic circulant graphs are `Type-1 isomorphic or Type-2 isomorphic'. Even though some triples of circulant graphs are already covered in problem \ref{p1}, we consider all the triples of isomorphic circulant graphs. 

From these 46 tables, we found that there are 960 = $2\times 3\times 4\times 2^3$ number of triples of Type-2 isomorphic circulant graphs w.r.t. $m$ = 3.

\begin{table}
	\caption{ {\footnotesize Finding $\theta_{54,3,2}(C_{54}(R))$, $\theta_{54,3,4}(C_{54}(R))$ and $T1_{54}(C_{54}(R))$ as $\theta_{54,3,6}(C_{54}(R))$ = $C_{54}(R)$.}}
\begin{center}
\scalebox{.78}{
\begin{tabular}{||c||c|c|c|c|c||}\hline \hline 
  & & & & &  \\
S. No. & $C_{54}(R)$ &  $\theta_{54,3,2}(C_{54}(R))$ = &  $\theta_{54,3,4}(C_{54}(R))$ =  & $C_{54}(S)\in T1_{54}(C_{54}(R))$ & T1 or T2  \\

 & & &   & &  \\ \hline \hline

1.  & $C_{54}(1,3,17,19)$ &  $C_{54}(3,7,11,25)$ &  $C_{54}(3,5,13,23)$ & $C_{54}(1,3,17,19)$, $C_{54}(5,13,15,23)$, &   \\ 

     & & & & $C_{54}(7,11,21,25)$.  &  T2 \\

2.   & $C_{54}(1,6,17,19)$ &   $C_{54}(6,7,11,25)$ &   $C_{54}(5,6,13,23)$ & $C_{54}(1,6,17,19)$, $C_{54}(5,13,23,24)$, &  \\ 

     & & & & $C_{54}(7,11,12,25)$.  &  T2 \\

3.   & $C_{54}(1,9,17,19)$ &   $C_{54}(7,9,11,25)$ &   $C_{54}(5,9,13,23)$ & $C_{54}(1,9,17,19)$, $C_{54}(5,9,13,23)$, &  \\
\hline \hline
\end{tabular}}
\end{center}
\end{table} 

\begin{table}
	\caption{ {\footnotesize Finding $\theta_{54,3,2}(C_{54}(R))$, $\theta_{54,3,4}(C_{54}(R))$ and $T1_{54}(C_{54}(R))$ as $\theta_{54,3,6}(C_{54}(R))$ = $C_{54}(R)$.}}
\begin{center}
\scalebox{.78}{
 \begin{tabular}{||c||c|c|c|c|c||}\hline \hline 
  & & & & &  \\
S. No. & $C_{54}(R)$ &  $\theta_{54,3,2}(C_{54}(R))$ = &  $\theta_{54,3,4}(C_{54}(R))$ =  & $C_{54}(S)\in T1_{54}(C_{54}(R))$ & T1 or  \\

  & & & & & T2  \\

 & & &   & &  \\ \hline \hline

32.   & $C_{54}(1,12,17,18,19)$ &   $C_{54}(7,11,12,18,25)$ &   $C_{54}(5,12,13,18,23)$ & $C_{54}(1,12,17,18,19)$, $C_{54}(5,6,13,18,23)$,   &  \\ 

     & & & &  $C_{54}(7,11,18,24,25)$.  &   T2 \\

33.   & $C_{54}(1,12,17,19,21)$ &   $C_{54}(7,11,12,21,25)$ &   $C_{54}(5,12,13,21,23)$ & $C_{54}(1,12,17,19,21)$, $C_{54}(3,5,6,13,23)$,   &  \\ 

     & & & &  $C_{54}(7,11,15,24,25)$.  &   T2 \\
\hline \hline
\end{tabular}}
\end{center}
\end{table} 

\begin{table}
	\caption{ {\footnotesize Finding $\theta_{54,3,2}(C_{54}(R))$, $\theta_{54,3,4}(C_{54}(R))$ and $T1_{54}(C_{54}(R))$ as $\theta_{54,3,6}(C_{54}(R))$ = $C_{54}(R)$.}}
\begin{center}
\scalebox{.8}{
 \begin{tabular}{||c||c|c|c|c|c||}\hline \hline 
  & & & & &  \\
S. No. & $C_{54}(R)$ &  $\theta_{54,3,2}(C_{54}(R))$ = &  $\theta_{54,3,4}(C_{54}(R))$ =  & $C_{54}(S)\in T1_{54}(C_{54}(R))$ & T1 or T2  \\

 & & &   & &  \\ \hline \hline

57.  & $C_{54}(1,3,9,17,19,24)$ &  $C_{54}(3,7,9,11,24,25)$ &  $C_{54}(3,5,9,13,23,24)$ & $C_{54}(1,3,9,17,19,24)$, & \\

& & & & $C_{54}(5,9,12,13,15,23)$, &   \\ 

     & & & & $C_{54}(6,7,9,11,21,25)$.  &  T2 \\

58.  & $C_{54}(1,3,9,17,19,27)$ &  $C_{54}(3,7,9,11,25,27)$ &  $C_{54}(3,5,9,13,23,27)$ & $C_{54}(1,3,9,17,19,27)$, & \\

& & & & $C_{54}(5,9,13,15,23,27)$, &   \\
\hline \hline
\end{tabular}}
\end{center}
\end{table} 

\begin{table}
	\caption{ {\footnotesize Finding $\theta_{54,3,2}(C_{54}(R))$, $\theta_{54,3,4}(C_{54}(R))$ and $T1_{54}(C_{54}(R))$ as $\theta_{54,3,6}(C_{54}(R))$ = $C_{54}(R)$.}}
\begin{center}
\scalebox{.8}{
 \begin{tabular}{||c||c|c|c|c|c||}\hline \hline 
  & & & & &  \\
S. No. & $C_{54}(R)$ &  $\theta_{54,3,2}(C_{54}(R))$ = &  $\theta_{54,3,4}(C_{54}(R))$ =  & $C_{54}(S)\in T1_{54}(C_{54}(R))$ & T1 or T2  \\

 & & &   & &  \\ \hline \hline

77.   & $C_{54}(1,6,9,17,19,21)$ &   $C_{54}(6,7,9,11,21,25)$ &   $C_{54}(5,6,9,13,21,23)$ & $C_{54}(1,6,9,17,19,21)$, & \\

& & & & $C_{54}(3,5,9,13,23,24)$, &  \\ 

     & & & & $C_{54}(7,9,11,12,15,25)$.  &  T2 \\

78.   & $C_{54}(1,6,9,17,19,24)$ &   $C_{54}(6,7,9,11,24,25)$ &   $C_{54}(5,6,9,13,23,24)$ & $C_{54}(1,6,9,17,19,24)$, & \\

& & & & $C_{54}(5,9,12,13,23,24)$, &  \\
\hline \hline
\end{tabular}}
\end{center}
\end{table} 

\begin{table}
	\caption{ {\footnotesize Finding $\theta_{54,3,2}(C_{54}(R))$, $\theta_{54,3,4}(C_{54}(R))$ and $T1_{54}(C_{54}(R))$ as $\theta_{54,3,6}(C_{54}(R))$ = $C_{54}(R)$.}}
\begin{center}
\scalebox{.8}{
 \begin{tabular}{||c||c|c|c|c|c||}\hline \hline 
  & & & & &  \\
S. No. & $C_{54}(R)$ &  $\theta_{54,3,2}(C_{54}(R))$ = &  $\theta_{54,3,4}(C_{54}(R))$ =  & $C_{54}(S)\in T1_{54}(C_{54}(R))$ & T1 or T2  \\

 & & &   & &  \\ \hline \hline

98.   & $C_{54}(1,9,12,17,19,24)$ &   $C_{54}(7,9,11,12,24,25)$ &   $C_{54}(5,9,12,13,23,24)$ & $C_{54}(1,9,12,17,19,24)$, & \\

& & & & $C_{54}(5,6,9,12,13,23)$, &  \\ 

     & & & & $C_{54}(6,7,9,11,24,25)$.  &  T2 \\

99.   & $C_{54}(1,9,12,17,19,27)$ &   $C_{54}(7,9,11,12,25,27)$ &   $C_{54}(5,9,12,13,23,27)$ & $C_{54}(1,9,12,17,19,27)$, & \\

& & & & $C_{54}(5,6,9,13,23,27)$, &  \\
\hline \hline
\end{tabular}}
\end{center}
\end{table} 

\begin{table}
	\caption{ {\footnotesize Finding $\theta_{54,3,2}(C_{54}(R))$, $\theta_{54,3,4}(C_{54}(R))$ and $T1_{54}(C_{54}(R))$ as $\theta_{54,3,6}(C_{54}(R))$ = $C_{54}(R)$.}}
\begin{center}
\scalebox{.8}{
 \begin{tabular}{||c||c|c|c|c|c||}\hline \hline 
  & & & & &  \\
S. No. & $C_{54}(R)$ &  $\theta_{54,3,2}(C_{54}(R))$ = &  $\theta_{54,3,4}(C_{54}(R))$ =  & $C_{54}(S)\in T1_{54}(C_{54}(R))$ & T1 or \\

& & & & & T2  \\

 & & &   & &  \\ \hline \hline

119.   & $C_{54}(1,12,17,19,24,27)$ &   $C_{54}(7,11,12,24,25,27)$ &   $C_{54}(5,12,13,23,24,27)$ & $C_{54}(1,12,17,19,24,27)$, & \\

& & & & $C_{54}(5,6,12,13,23,27)$,   &  \\ 

     & & & &  $C_{54}(6,7,11,24,25,27)$.  &   T2 \\

120.   & $C_{54}(1,15,17,18,19,21)$ &   $C_{54}(7,11,15,18,21,25)$ &   $C_{54}(5,13,15,18,21,23)$ & $C_{54}(1,15,17,18,19,21)$, & \\
\hline \hline
\end{tabular}}
\end{center}
\end{table} 

\begin{table}
	\caption{ {\footnotesize Finding $\theta_{54,3,2}(C_{54}(R))$, $\theta_{54,3,4}(C_{54}(R))$ and $T1_{54}(C_{54}(R))$ as $\theta_{54,3,6}(C_{54}(R))$ = $C_{54}(R)$.}}
\begin{center}
\scalebox{.8}{
 \begin{tabular}{||c||c|c|c|c|c||}\hline \hline 
  & & & & &  \\
S. No. & $C_{54}(R)$ &  $\theta_{54,3,2}(C_{54}(R))$ = &  $\theta_{54,3,4}(C_{54}(R))$ =  & $C_{54}(S)\in T1_{54}(C_{54}(R))$ & T1 or \\

& & & & & T2  \\

 & & &   & &  \\ \hline \hline

139.  & $C_{54}(1,3,6,12,17,19,24)$ &  $C_{54}(3,6,7,11,12,24,25)$ &  $C_{54}(3,5,6,12,13,23,24)$ & $C_{54}(1,3,6,12,17,19,24)$,  & \\   

& & & & $C_{54}(5,6,12,13,15,23,24)$, &   \\ 

     & & & & $C_{54}(6,7,11,12,21,24,25)$.  &  T2 \\

140.  & $C_{54}(1,3,6,12,17,19,27)$ &  $C_{54}(3,6,7,11,12,25,27)$ &  $C_{54}(3,5,6,12,13,23,27)$ & $C_{54}(1,3,6,12,17,19,27)$,  & \\
\hline \hline
\end{tabular}}
\end{center}
\end{table} 

\begin{table}
	\caption{ {\footnotesize Finding $\theta_{54,3,2}(C_{54}(R))$, $\theta_{54,3,4}(C_{54}(R))$ and $T1_{54}(C_{54}(R))$ as $\theta_{54,3,6}(C_{54}(R))$ = $C_{54}(R)$.}}
\begin{center}
\scalebox{.78}{
 \begin{tabular}{||c||c|c|c|c|c||}\hline \hline 
  & & & & &  \\
S. No. & $C_{54}(R)$ &  $\theta_{54,3,2}(C_{54}(R))$ = &  $\theta_{54,3,4}(C_{54}(R))$ =  & $C_{54}(S)\in T1_{54}(C_{54}(R))$ & T1   or \\

& & & & & T2  \\

 & & &   & &  \\ \hline \hline

159.  & $C_{54}(1,3,9,15,17,19,27)$ &  $C_{54}(3,7,9,11,15,25,27)$ &  $C_{54}(3,5,9,13,15,23,27)$ & $C_{54}(1,3,9,15,17,19,27)$, & \\

& & & & $C_{54}(5,9,13,15,21,23,27)$, &   \\ 

     & & & & $C_{54}(3,7,9,11,21,25,27)$.  &  T2 \\

160.  & $C_{54}(1,3,9,17,18,19,21)$ &  $C_{54}(3,7,9,11,18,21,25)$ &  $C_{54}(3,5,9,13,18,21,23)$ & $C_{54}(1,3,9,17,18,19,21)$, & \\
\hline \hline
\end{tabular}}
\end{center}
\end{table} 

\begin{table}
	\caption{ {\footnotesize Finding $\theta_{54,3,2}(C_{54}(R))$, $\theta_{54,3,4}(C_{54}(R))$ and $T1_{54}(C_{54}(R))$ as $\theta_{54,3,6}(C_{54}(R))$ = $C_{54}(R)$.}}
\begin{center}
\scalebox{.77}{
 \begin{tabular}{||c||c|c|c|c|c||}\hline \hline 
  & & & & &  \\
S. No. & $C_{54}(R)$ &  $\theta_{54,3,2}(C_{54}(R))$ = &  $\theta_{54,3,4}(C_{54}(R))$ =  & $C_{54}(S)\in T1_{54}(C_{54}(R))$ & T1   or \\

& & & & & T2  \\

 & & &   & &  \\ \hline \hline

179.  & $C_{54}(1,3,15,17,19,21,24)$ &  $C_{54}(3,7,11,15,21,24,25)$ &  $C_{54}(3,5,13,15,21,23,24)$ & $C_{54}(1,3,15,17,19,21,24)$, & \\

& & & & $C_{54}(3,5,12,13,15,21,23)$, &   \\ 

     & & & & $C_{54}(3,6,7,11,15,21,25)$.  &  T2 \\

180.  & $C_{54}(1,3,15,17,19,21,27)$ &  $C_{54}(3,7,11,15,21,25,27)$ &  $C_{54}(3,5,13,15,21,23,27)$ & $C_{54}(1,3,15,17,19,21,27)$, & \\
\hline \hline
\end{tabular}}
\end{center}
\end{table} 

\begin{table}
	\caption{ {\footnotesize Finding $\theta_{54,3,2}(C_{54}(R))$, $\theta_{54,3,4}(C_{54}(R))$ $\&$ $T1_{54}(C_{54}(R))$ as $\theta_{54,3,6}(C_{54}(R))$ = $C_{54}(R)$.}}
\begin{center}
\scalebox{.75}{
 \begin{tabular}{||c||c|c|c|c|c||}\hline \hline 
  & & & & &  \\
S. No. & $C_{54}(R)$ &  $\theta_{54,3,2}(C_{54}(R))$ = &  $\theta_{54,3,4}(C_{54}(R))$ =  & $C_{54}(S)\in T1_{54}(C_{54}(R))$ & T1   or \\

& & & & & T2  \\

 & & &   & &  \\ \hline \hline

200.   & $C_{54}(1,6,9,17,19,24,27)$ &   $C_{54}(6,7,9,11,24,25,27)$ &   $C_{54}(5,6,9,13,23,24,27)$ & $C_{54}(1,6,9,17,19,24,27)$, & \\

& & & & $C_{54}(5,9,12,13,23,24,27)$, &  \\ 

     & & & & $C_{54}(6,7,9,11,12,25,27)$.  &  T2 \\

201.   & $C_{54}(1,6,12,15,17,18,19)$ &   $C_{54}(6,7,11,12,15,18,25)$ &   $C_{54}(5,6,12,13,15,18,23)$ & $C_{54}(1,6,12,15,17,18,19)$, & \\
\hline \hline
\end{tabular}}
\end{center}
\end{table} 

\begin{table}
	\caption{ {\footnotesize Finding $\theta_{54,3,2}(C_{54}(R))$, $\theta_{54,3,4}(C_{54}(R))$ $\&$ $T1_{54}(C_{54}(R))$ as $\theta_{54,3,6}(C_{54}(R))$ = $C_{54}(R)$.}}
\begin{center}
\scalebox{.75}{
 \begin{tabular}{||c||c|c|c|c|c||}\hline \hline 
  & & & & &  \\
S. No. & $C_{54}(R)$ &  $\theta_{54,3,2}(C_{54}(R))$ = &  $\theta_{54,3,4}(C_{54}(R))$ =  & $C_{54}(S)\in T1_{54}(C_{54}(R))$ & T1   or \\

& & & & & T2  \\

 & & &   & &  \\ \hline \hline

222.   & $C_{54}(1,9,12,15,17,19,21)$ &   $C_{54}(7,9,11,12,15,21,25)$ &   $C_{54}(5,9,12,13,15,21,23)$ & $C_{54}(1,9,12,15,17,19,21)$, & \\

& & & & $C_{54}(3,5,6,9,13,21,23)$, &  \\ 

     & & & & $C_{54}(3,7,9,11,15,24,25)$.  &  T2 \\

223.   & $C_{54}(1,9,12,15,17,19,24)$ &   $C_{54}(7,9,11,12,15,24,25)$ &   $C_{54}(5,9,12,13,15,23,24)$ & $C_{54}(1,9,12,15,17,19,24)$, & \\
\hline \hline
\end{tabular}}
\end{center}
\end{table} 

\begin{table}
	\caption{ {\footnotesize Finding $\theta_{54,3,2}(C_{54}(R))$, $\theta_{54,3,4}(C_{54}(R))$ $\&$ $T1_{54}(C_{54}(R))$ as $\theta_{54,3,6}(C_{54}(R))$ = $C_{54}(R)$.}}
\begin{center}
\scalebox{.75}{
 \begin{tabular}{||c||c|c|c|c|c||}\hline \hline 
  & & & & &  \\
S. No. & $C_{54}(R)$ &  $\theta_{54,3,2}(C_{54}(R))$ = &  $\theta_{54,3,4}(C_{54}(R))$ =  & $C_{54}(S)\in T1_{54}(C_{54}(R))$ & T1   or \\

& & & & & T2  \\

 & & &   & &  \\ \hline \hline

244.   & $C_{54}(1,12,15,17,19,21,24)$ &   $C_{54}(7,11,12,15,21,24,25)$ &   $C_{54}(5,12,13,15,21,23,24)$ & $C_{54}(1,12,15,17,19,21,24)$, & \\

& & & & $C_{54}(3,5,6,13,12,21,23)$, &  \\ 

     & & & & $C_{54}(3,6,7,11,15,24,25)$.  &  T2 \\

245.   & $C_{54}(1,12,15,17,19,21,27)$ &   $C_{54}(7,11,12,15,21,25,27)$ &   $C_{54}(5,12,13,15,21,23,27)$ & $C_{54}(1,12,15,17,19,21,27)$, & \\
\hline \hline
\end{tabular}}
\end{center}
\end{table} 

\begin{table}
	\caption{ {\footnotesize Finding $\theta_{54,3,2}(C_{54}(R))$, $\theta_{54,3,4}(C_{54}(R))$ $\&$ $T1_{54}(C_{54}(R))$ as $\theta_{54,3,6}(C_{54}(R))$ = $C_{54}(R)$.}}
\begin{center}
\scalebox{.75}{
 \begin{tabular}{||c||c|c|c|c|c||}\hline \hline 
  & & & & &  \\
S. No. & $C_{54}(R)$ &  $\theta_{54,3,2}(C_{54}(R))$ = &  $\theta_{54,3,4}(C_{54}(R))$ =  & $C_{54}(S)\in T1_{54}(C_{54}(R))$ & T1   or \\

& & & & & T2  \\

 & & &   & &  \\ \hline \hline

266.  & $C_{54}(1,3,6,9,17,18,19,24)$ &  $C_{54}(3,6,7,9,11,18,24,25)$ &  $C_{54}(3,5,6,9,13,18,23,24)$ & $C_{54}(1,3,6,9,17,18,19,24)$,  & \\   

& & & & $C_{54}(5,9,12,13,15,18,23,24)$, &   \\ 

     & & & & $C_{54}(6,7,9,11,12,18,21,25)$.  &  T2 \\

267.  & $C_{54}(1,3,6,9,17,18,19,27)$ &  $C_{54}(3,6,7,9,11,18,25,27)$ &  $C_{54}(3,5,6,9,13,18,23,27)$ & $C_{54}(1,3,6,9,17,18,19,27)$,  & \\
\hline \hline
\end{tabular}}
\end{center}
\end{table} 

\begin{table}
	\caption{ {\footnotesize Finding $\theta_{54,3,2}(C_{54}(R))$, $\theta_{54,3,4}(C_{54}(R))$ $\&$ $T1_{54}(C_{54}(R))$ as $\theta_{54,3,6}(C_{54}(R))$ = $C_{54}(R)$.}}
\begin{center}
\scalebox{.75}{
 \begin{tabular}{||c||c|c|c|c|c||}\hline \hline 
  & & & & &  \\
S. No. & $C_{54}(R)$ &  $\theta_{54,3,2}(C_{54}(R))$ = &  $\theta_{54,3,4}(C_{54}(R))$ =  & $C_{54}(S)\in T1_{54}(C_{54}(R))$ & T1   or \\

& & & & & T2  \\

 & & &   & &  \\ \hline \hline

288.  & $C_{54}(1,3,6,17,18,19,21,27)$ &  $C_{54}(3,6,7,11,18,21,25,27)$ &  $C_{54}(3,5,6,13,18,21,23,27)$ & $C_{54}(1,3,6,17,18,19,21,27)$,  & \\   

& & & & $C_{54}(3,5,13,15,18,23,24,27)$, &   \\ 

     & & & & $C_{54}(7,11,12,15,18,21,25,27)$.  &  T2 \\

289. & $C_{54}(1,3,6,17,18,19,24,27)$ &  $C_{54}(3,6,7,11,18,24,25,27)$ &  $C_{54}(3,5,6,13,18,23,24,27)$ & $C_{54}(1,3,6,17,18,19,24,27)$,  & \\
\hline \hline
\end{tabular}}
\end{center}
\end{table} 

\begin{table}
	\caption{ {\footnotesize Finding $\theta_{54,3,2}(C_{54}(R))$, $\theta_{54,3,4}(C_{54}(R))$ $\&$ $T1_{54}(C_{54}(R))$ as $\theta_{54,3,6}(C_{54}(R))$ = $C_{54}(R)$.}}
\begin{center}
\scalebox{.73}{
 \begin{tabular}{||c||c|c|c|c|c||}\hline \hline 
  & & & & &  \\
S. No. & $C_{54}(R)$ &  $\theta_{54,3,2}(C_{54}(R))$ = &  $\theta_{54,3,4}(C_{54}(R))$ =  & $C_{54}(S)\in T1_{54}(C_{54}(R))$ & T1   or \\

& & & & & T2  \\

 & & &   & &  \\ \hline \hline

310.  & $C_{54}(1,3,9,17,19,21,24,27)$ &  $C_{54}(3,7,9,11,21,24,25,27)$ &  $C_{54}(3,5,9,13,21,23,24,27)$ & $C_{54}(1,3,9,17,19,21,24,27)$, & \\

& & & & $C_{54}(3,5,9,12,13,15,23,27)$, &   \\ 

     & & & & $C_{54}(6,7,9,11,15,21,25,27)$.  &  T2 \\

311.  & $C_{54}(1,3,12,15,17,18,19,21)$ &  $C_{54}(3,7,11,12,15,18,21,25)$ &  $C_{54}(3,5,12,13,15,18,21,23)$ & $C_{54}(1,3,12,15,17,18,19,21)$, & \\
\hline \hline
\end{tabular}}
\end{center}
\end{table} 

\begin{table}
	\caption{ {\footnotesize Finding $\theta_{54,3,2}(C_{54}(R))$, $\theta_{54,3,4}(C_{54}(R))$ $\&$ $T1_{54}(C_{54}(R))$ as $\theta_{54,3,6}(C_{54}(R))$ = $C_{54}(R)$.}}
\begin{center}
\scalebox{.7}{
 \begin{tabular}{||c||c|c|c|c|c||}\hline \hline 
  & & & & &  \\
S. No. & $C_{54}(R)$ &  $\theta_{54,3,2}(C_{54}(R))$ = &  $\theta_{54,3,4}(C_{54}(R))$ =  & $C_{54}(S)\in T1_{54}(C_{54}(R))$ & T1   or \\

& & & & & T2  \\

 & & &   & &  \\ \hline \hline

332.  & $C_{54}(1,6,9,12,17,18,19,27)$ &  $C_{54}(6,7,9,11,12,18,25,27)$ &  $C_{54}(5,6,9,12,13,18,23,27)$ & $C_{54}(1,6,9,12,17,18,19,27)$,  & \\   

& & & & $C_{54}(5,6,9,13,18,23,24,27)$, &   \\ 

     & & & & $C_{54}(7,9,11,12,18,24,25,27)$.  &  T2 \\

333.  & $C_{54}(1,6,9,12,17,19,21,24)$ &  $C_{54}(6,7,9,11,12,21,24,25)$ &  $C_{54}(5,6,9,12,13,21,23,24)$ & $C_{54}(1,6,9,12,17,19,21,24)$,  & \\
\hline \hline
\end{tabular}}
\end{center}
\end{table} 

\begin{table}
	\caption{ {\footnotesize Finding $\theta_{54,3,2}(C_{54}(R))$, $\theta_{54,3,4}(C_{54}(R))$ $\&$ $T1_{54}(C_{54}(R))$ as $\theta_{54,3,6}(C_{54}(R))$ = $C_{54}(R)$.}}
\begin{center}
\scalebox{.7}{
 \begin{tabular}{||c||c|c|c|c|c||}\hline \hline 
  & & & & &  \\
S. No. & $C_{54}(R)$ &  $\theta_{54,3,2}(C_{54}(R))$ = &  $\theta_{54,3,4}(C_{54}(R))$ =  & $C_{54}(S)\in T1_{54}(C_{54}(R))$ & T1   or \\

& & & & & T2  \\

 & & &   & &  \\ \hline \hline

355.  & $C_{54}(1,6,12,17,19,21,24,27)$ &  $C_{54}(6,7,11,12,21,24,25,27)$ &  $C_{54}(5,6,12,13,21,23,24,27)$ & $C_{54}(1,6,12,17,19,21,24,27)$,  & \\   

& & & & $C_{54}(3,5,6,12,13,23,24,27)$, &   \\ 

     & & & & $C_{54}(6,7,11,12,15,24,25,27)$.  &  T2 \\

356.  & $C_{54}(1,6,15,17,18,19,21,24)$ &  $C_{54}(6,7,11,15,18,21,24,25)$ &  $C_{54}(5,6,13,15,18,21,23,24)$ & $C_{54}(1,6,15,17,18,19,21,24)$,  & \\
\hline \hline
\end{tabular}}
\end{center}
\end{table} 

\begin{table}
	\caption{ {\footnotesize Finding $\theta_{54,3,2}(C_{54}(R))$, $\theta_{54,3,4}(C_{54}(R))$ $\&$ $T1_{54}(C_{54}(R))$ as $\theta_{54,3,6}(C_{54}(R))$ = $C_{54}(R)$.}}
\begin{center}
\scalebox{.68}{
 \begin{tabular}{||c||c|c|c|c|c||}\hline \hline 
  & & & & &  \\
S. No. & $C_{54}(R)$ &  $\theta_{54,3,2}(C_{54}(R))$ = &  $\theta_{54,3,4}(C_{54}(R))$ =  & $C_{54}(S)\in T1_{54}(C_{54}(R))$ & T1   or \\

& & & & & T2  \\

 & & &   & &  \\ \hline \hline

378.  & $C_{54}(1,12,15,17,18,19,24,27)$ &  $C_{54}(7,11,12,15,18,24,25,27)$ &  $C_{54}(5,12,13,15,18,23,24,27)$ & $C_{54}(1,12,15,17,18,19,24,27)$, & \\

& & & & $C_{54}(5,6,12,13,18,21,23,27)$, &   \\ 

     & & & & $C_{54}(3,6,7,11,18,21,25,27)$.  &  T2 \\ 

379.  & $C_{54}(1,12,15,17,19,21,24,27)$ &  $C_{54}(7,11,12,15,21,24,25,27)$ &  $C_{54}(5,12,13,15,21,23,24,27)$ & $C_{54}(1,12,15,17,19,21,24,27)$, & \\
\hline \hline
\end{tabular}}
\end{center}
\end{table} 

\begin{table}
	\caption{ {\footnotesize Finding $\theta_{54,3,2}(C_{54}(R))$, $\theta_{54,3,4}(C_{54}(R))$ $\&$ $T1_{54}(C_{54}(R))$ as $\theta_{54,3,6}(C_{54}(R))$ = $C_{54}(R)$.}}
\begin{center}
\scalebox{.65}{
 \begin{tabular}{||c||c|c|c|c|c||}\hline \hline 
  & & & & &  \\
S. No. & $C_{54}(R)$ &  $\theta_{54,3,2}(C_{54}(R))$ = &  $\theta_{54,3,4}(C_{54}(R))$ =  & $C_{54}(S)\in T1_{54}(C_{54}(R))$ & T1   or \\

& & & & & T2  \\

 & & &   & &  \\ \hline \hline

402.  & $C_{54}(1,3,6,12,15,17,18,19,21)$ &  $C_{54}(3,6,7,11,12,15,18,21,25)$ &  $C_{54}(3,5,6,12,13,15,18,21,23)$ & $C_{54}(1,3,6,12,15,17,18,19,21)$,  & \\   

& & & & $C_{54}(3,5,6,13,15,18,21,23,24)$, &   \\ 

     & & & & $C_{54}(3,7,11,12,15,18,21,24,25)$.  &  T2 \\

403.  & $C_{54}(1,3,6,12,15,17,18,19,24)$ &  $C_{54}(3,6,7,11,12,15,18,24,25)$ &  $C_{54}(3,5,6,12,13,15,18,23,24)$ & $C_{54}(1,3,6,12,15,17,18,19,24)$,  & \\
\hline \hline
\end{tabular}}
\end{center}
\end{table} 

\begin{table}
	\caption{ {\footnotesize Finding $\theta_{54,3,2}(C_{54}(R))$, $\theta_{54,3,4}(C_{54}(R))$ $\&$ $T1_{54}(C_{54}(R))$ as $\theta_{54,3,6}(C_{54}(R))$ = $C_{54}(R)$.}}
\begin{center}
\scalebox{.68}{
 \begin{tabular}{||c||c|c|c|c|c||}\hline \hline 
  & & & & &  \\
S. No. & $C_{54}(R)$ &  $\theta_{54,3,2}(C_{54}(R))$ = &  $\theta_{54,3,4}(C_{54}(R))$ =  & $C_{54}(S)\in T1_{54}(C_{54}(R))$ & T1   or \\

& & & & & T2  \\

 & & &   & &  \\ \hline \hline

427.  & $C_{54}(1,3,9,15,17,18,19,21,24)$ &  $C_{54}(3,7,9,11,15,18,21,24,25)$ &  $C_{54}(3,5,9,13,15,18,21,23,24)$ & $C_{54}(1,3,9,15,17,18,19,21,24)$, & \\

& & & & $C_{54}(3,5,9,12,13,15,18,21,23)$, &   \\ 

     & & & & $C_{54}(3,6,7,9,11,15,18,21,25)$.  &  T2 \\

428.  & $C_{54}(1,3,9,15,17,18,19,21,27)$ &  $C_{54}(3,7,9,11,15,18,21,25,27)$ &  $C_{54}(3,5,9,13,15,18,21,23,27)$ & $C_{54}(1,3,9,15,17,18,19,21,27)$, & \\
\hline \hline
\end{tabular}}
\end{center}
\end{table} 

\begin{table}
	\caption{ {\footnotesize Finding $\theta_{54,3,2}(C_{54}(R))$, $\theta_{54,3,4}(C_{54}(R))$ $\&$ $T1_{54}(C_{54}(R))$ as $\theta_{54,3,6}(C_{54}(R))$ = $C_{54}(R)$.}}
\begin{center}
\scalebox{.68}{
 \begin{tabular}{||c||c|c|c|c|c||}\hline \hline 
  & & & & &  \\
S. No. & $C_{54}(R)$ &  $\theta_{54,3,2}(C_{54}(R))$ = &  $\theta_{54,3,4}(C_{54}(R))$ =  & $C_{54}(S)\in T1_{54}(C_{54}(R))$ & T1   or \\

& & & & & T2  \\

 & & &   & &  \\ \hline \hline

451.  & $C_{54}(1,6,9,15,17,19,21,24,27)$ &  $C_{54}(6,7,9,11,15,21,24,25,27)$ &  $C_{54}(5,6,9,13,15,21,23,24,27)$ & $C_{54}(1,6,9,15,17,19,21,24,27)$,  & \\   

& & & & $C_{54}(3,5,9,12,13,21,23,24,27)$, &   \\ 

     & & & & $C_{54}(3,6,7,9,11,12,15,25,27)$.  &  T2 \\

452.  & $C_{54}(1,6,9,17,18,19,21,24,27)$ &  $C_{54}(6,7,9,11,18,21,24,25,27)$ &  $C_{54}(5,6,9,13,18,21,23,24,27)$ & $C_{54}(1,6,9,17,18,19,21,24,27)$,  & \\
\hline \hline
\end{tabular}}
\end{center}
\end{table} 

\begin{table}
	\caption{ {\footnotesize Finding $\theta_{54,3,2}(C_{54}(R))$, $\theta_{54,3,4}(C_{54}(R))$ $\&$ $T1_{54}(C_{54}(R))$ as $\theta_{54,3,6}(C_{54}(R))$ = $C_{54}(R)$.}}
\begin{center}
\scalebox{.59}{
 \begin{tabular}{||c||c|c|c|c|c||}\hline \hline 
  & & & & &  \\
S. No. & $C_{54}(R)$ &  $\theta_{54,3,2}(C_{54}(R))$ = &  $\theta_{54,3,4}(C_{54}(R))$ =  & $C_{54}(S)\in T1_{54}(C_{54}(R))$ & T1   or \\

& & & & & T2  \\

 & & &   & &  \\ \hline \hline

475.  & $C_{54}(1,3,6,9,12,17,19,21,24,27)$ &  $C_{54}(3,6,7,9,11,12,21,24,25,27)$ &  $C_{54}(3,5,6,9,12,13,21,23,24,27)$ & $C_{54}(1,3,6,9,12,17,19,21,24,27)$,  & \\   

& & & & $C_{54}(3,5,6,9,12,13,15,23,24,27)$, &   \\ 

     & & & & $C_{54}(6,7,9,11,12,15,21,24,25,27)$.  &  T2 \\

476.  & $C_{54}(1,3,6,9,15,17,18,19,21,24)$ &  $C_{54}(3,6,7,9,11,15,18,21,24,25)$ &  $C_{54}(3,5,6,9,13,15,18,21,23,24)$ & $C_{54}(1,3,6,9,15,17,18,19,21,24)$,  & \\
\hline \hline
\end{tabular}}
\end{center}
\end{table} 

\begin{table}
	\caption{ {\footnotesize Finding $\theta_{54,3,2}(C_{54}(R))$, $\theta_{54,3,4}(C_{54}(R))$ $\&$ $T1_{54}(C_{54}(R))$ as $\theta_{54,3,6}(C_{54}(R))$ = $C_{54}(R)$.}}
\begin{center}
\scalebox{.59}{
 \begin{tabular}{||c||c|c|c|c|c||}\hline \hline 
  & & & & &  \\
S. No. & $C_{54}(R)$ &  $\theta_{54,3,2}(C_{54}(R))$ = &  $\theta_{54,3,4}(C_{54}(R))$ =  & $C_{54}(S)\in T1_{54}(C_{54}(R))$ & T1   or \\

& & & & & T2  \\

 & & &   & &  \\ \hline \hline

502.  & $C_{54}(1,3,6,9,12,15,17,18,19,21,24)$ &  $C_{54}(3,6,7,9,11,12,15,18,21,24,25)$ &  $C_{54}(3,5,6,9,12,13,15,18,21,23,24)$ & $C_{54}(1,3,6,9,12,15,17,18,19,21,24)$,  & \\   

& & & & $C_{54}(3,5,6,9,12,13,15,18,21,23,24)$, &   \\ 

     & & & & $C_{54}(3,6,7,9,11,12,15,18,21,24,25)$.  &  T1 \\

503.  & $C_{54}(1,3,6,9,12,15,17,18,19,21,27)$ &  $C_{54}(3,6,7,9,11,12,15,18,21,25,27)$ &  $C_{54}(3,5,6,9,12,13,15,18,21,23,27)$ & $C_{54}(1,3,6,9,12,15,17,18,19,21,27)$,  & \\
\hline \hline
\end{tabular}}
\end{center}
\end{table} 

\begin{table}
	\caption{ {\footnotesize Finding $\theta_{54,3,2}(C_{54}(R))$, $\theta_{54,3,4}(C_{54}(R))$ $\&$ $T1_{54}(C_{54}(R))$ as $\theta_{54,3,6}(C_{54}(R))$ = $C_{54}(R)$.}}
\begin{center}
\scalebox{.78}{
 \begin{tabular}{||c||c|c|c|c|c||}\hline \hline 
  & & & & &  \\
S. No. & $C_{54}(R)$ &  $\theta_{54,3,2}(C_{54}(R))$ = &  $\theta_{54,3,4}(C_{54}(R))$ =  & $C_{54}(S)\in T1_{54}(C_{54}(R))$ & T1 or T2  \\

 & & &   & &  \\ \hline \hline

1.  & $C_{54}(2,3,16,20)$ &  $C_{54}(3,4,14,22)$ &  $C_{54}(3,8,10,26)$ & $C_{54}(2,3,16,20)$, $C_{54}(8,10,15,26)$, &   \\ 

     & & & & $C_{54}(4,14,21,22)$.  &  T2 \\

2.   & $C_{54}(2,6,16,20)$ &   $C_{54}(4,6,14,22)$ &   $C_{54}(6,8,10,26)$ & $C_{54}(2,6,16,20)$, $C_{54}(8,10,24,26)$, &  \\ 

     & & & & $C_{54}(4,12,14,22)$.  &  T2 \\

3.   & $C_{54}(2,9,16,20)$ &   $C_{54}(4,9,14,22)$ &   $C_{54}(8,9,10,26)$ & $C_{54}(2,9,16,20)$, $C_{54}(8,9,10,26)$, &  \\
\hline \hline
\end{tabular}}
\end{center}
\end{table} 

\begin{table}
	\caption{ {\footnotesize Finding $\theta_{54,3,2}(C_{54}(R))$, $\theta_{54,3,4}(C_{54}(R))$ $\&$ $T1_{54}(C_{54}(R))$ as $\theta_{54,3,6}(C_{54}(R))$ = $C_{54}(R)$.}}
\begin{center}
\scalebox{.78}{
 \begin{tabular}{||c||c|c|c|c|c||}\hline \hline 
  & & & & &  \\
S. No. & $C_{54}(R)$ &  $\theta_{54,3,2}(C_{54}(R))$ = &  $\theta_{54,3,4}(C_{54}(R))$ =  & $C_{54}(S)\in T1_{54}(C_{54}(R))$ & T1 or  \\

  & & & & & T2  \\

 & & &   & &  \\ \hline \hline

32.   & $C_{54}(2,12,16,18,20)$ &   $C_{54}(4,12,14,18,22)$ &   $C_{54}(8,10,12,18,26)$ & $C_{54}(2,12,16,18,20)$, $C_{54}(6,8,10,18,26)$,   &  \\ 

     & & & &  $C_{54}(4,14,18,22,24)$.  &   T2 \\

33.   & $C_{54}(2,12,16,20,21)$ &   $C_{54}(4,12,14,21,22)$ &   $C_{54}(8,10,12,21,26)$ & $C_{54}(2,12,16,20,21)$, $C_{54}(3,6,8,10,26)$,   &  \\ 

     & & & &  $C_{54}(4,14,15,22,24)$.  &   T2 \\
\hline \hline
\end{tabular}}
\end{center}
\end{table} 

\begin{table}
	\caption{ {\footnotesize Finding $\theta_{54,3,2}(C_{54}(R))$, $\theta_{54,3,4}(C_{54}(R))$ $\&$ $T1_{54}(C_{54}(R))$ as $\theta_{54,3,6}(C_{54}(R))$ = $C_{54}(R)$.}}
\begin{center}
\scalebox{.8}{
 \begin{tabular}{||c||c|c|c|c|c||}\hline \hline 
  & & & & &  \\
S. No. & $C_{54}(R)$ &  $\theta_{54,3,2}(C_{54}(R))$ = &  $\theta_{54,3,4}(C_{54}(R))$ =  & $C_{54}(S)\in T1_{54}(C_{54}(R))$ & T1 or T2  \\

 & & &   & &  \\ \hline \hline

57.  & $C_{54}(2,3,9,16,20,24)$ &  $C_{54}(3,4,9,14,22,24)$ &  $C_{54}(3,8,9,10,24,26)$ & $C_{54}(2,3,9,16,20,24)$, & \\

& & & & $C_{54}(8,9,10,12,15,26)$, &   \\ 

     & & & & $C_{54}(4,6,9,14,21,22)$.  &  T2 \\

58.  & $C_{54}(2,3,9,16,20,27)$ &  $C_{54}(3,4,9,14,22,27)$ &  $C_{54}(3,8,9,10,26,27)$ & $C_{54}(2,3,9,16,20,27)$, & \\

& & & & $C_{54}(8,9,10,15,26,27)$, &   \\
\hline \hline
\end{tabular}}
\end{center}
\end{table} 

\begin{table}
	\caption{ {\footnotesize Finding $\theta_{54,3,2}(C_{54}(R))$, $\theta_{54,3,4}(C_{54}(R))$ $\&$ $T1_{54}(C_{54}(R))$ as $\theta_{54,3,6}(C_{54}(R))$ = $C_{54}(R)$.}}
\begin{center}
\scalebox{.8}{
 \begin{tabular}{||c||c|c|c|c|c||}\hline \hline 
  & & & & &  \\
S. No. & $C_{54}(R)$ &  $\theta_{54,3,2}(C_{54}(R))$ = &  $\theta_{54,3,4}(C_{54}(R))$ =  & $C_{54}(S)\in T1_{54}(C_{54}(R))$ & T1 or T2  \\

 & & &   & &  \\ \hline \hline

77.   & $C_{54}(2,6,9,16,20,21)$ &   $C_{54}(4,6,9,14,21,22)$ &   $C_{54}(6,8,9,10,21,26)$ & $C_{54}(2,6,9,16,20,21)$, & \\

& & & & $C_{54}(3,8,9,10,24,26)$, &  \\ 

     & & & & $C_{54}(4,9,12,14,15,22)$.  &  T2 \\

78.   & $C_{54}(2,6,9,16,20,24)$ &   $C_{54}(4,6,9,14,22,24)$ &   $C_{54}(6,8,9,10,24,26)$ & $C_{54}(2,6,9,16,20,24)$, & \\

& & & & $C_{54}(8,9,10,12,24,26)$, &  \\
\hline \hline
\end{tabular}}
\end{center}
\end{table} 

\begin{table}
	\caption{ {\footnotesize Finding $\theta_{54,3,2}(C_{54}(R))$, $\theta_{54,3,4}(C_{54}(R))$ $\&$ $T1_{54}(C_{54}(R))$ as $\theta_{54,3,6}(C_{54}(R))$ = $C_{54}(R)$.}}
\begin{center}
\scalebox{.8}{
 \begin{tabular}{||c||c|c|c|c|c||}\hline \hline 
  & & & & &  \\
S. No. & $C_{54}(R)$ &  $\theta_{54,3,2}(C_{54}(R))$ = &  $\theta_{54,3,4}(C_{54}(R))$ =  & $C_{54}(S)\in T1_{54}(C_{54}(R))$ & T1 or T2  \\

 & & &   & &  \\ \hline \hline

98.   & $C_{54}(2,9,12,16,20,24)$ &   $C_{54}(4,9,12,14,22,24)$ &   $C_{54}(8,9,10,12,24,26)$ & $C_{54}(2,9,12,16,20,24)$, & \\

& & & & $C_{54}(6,8,9,10,12,26)$, &  \\ 

     & & & & $C_{54}(4,6,9,14,22,24)$.  &  T2 \\

99.   & $C_{54}(2,9,12,16,20,27)$ &   $C_{54}(4,9,12,14,22,27)$ &   $C_{54}(8,9,10,12,26,27)$ & $C_{54}(2,9,12,16,20,27)$, & \\

& & & & $C_{54}(6,8,9,10,26,27)$, &  \\
\hline \hline
\end{tabular}}
\end{center}
\end{table} 

\begin{table}
	\caption{ {\footnotesize Finding $\theta_{54,3,2}(C_{54}(R))$, $\theta_{54,3,4}(C_{54}(R))$ $\&$ $T1_{54}(C_{54}(R))$ as $\theta_{54,3,6}(C_{54}(R))$ = $C_{54}(R)$.}}
\begin{center}
\scalebox{.8}{
 \begin{tabular}{||c||c|c|c|c|c||}\hline \hline 
  & & & & &  \\
S. No. & $C_{54}(R)$ &  $\theta_{54,3,2}(C_{54}(R))$ = &  $\theta_{54,3,4}(C_{54}(R))$ =  & $C_{54}(S)\in T1_{54}(C_{54}(R))$ & T1 or \\

& & & & & T2  \\

 & & &   & &  \\ \hline \hline

119.   & $C_{54}(2,12,16,20,24,27)$ &   $C_{54}(4,12,14,22,24,27)$ &   $C_{54}(8,10,12,24,26,27)$ & $C_{54}(2,12,16,20,24,27)$, & \\

& & & & $C_{54}(6,8,10,12,26,27)$,   &  \\ 

     & & & &  $C_{54}(4,6,14,22,24,27)$.  &   T2 \\

120.   & $C_{54}(2,15,16,18,20,21)$ &   $C_{54}(4,14,15,18,21,22)$ &   $C_{54}(8,10,15,18,21,26)$ & $C_{54}(2,15,16,18,20,21)$, & \\
\hline \hline
\end{tabular}}
\end{center}
\end{table} 

\begin{table}
	\caption{ {\footnotesize Finding $\theta_{54,3,2}(C_{54}(R))$, $\theta_{54,3,4}(C_{54}(R))$ $\&$ $T1_{54}(C_{54}(R))$ as $\theta_{54,3,6}(C_{54}(R))$ = $C_{54}(R)$.}}
\begin{center}
\scalebox{.8}{
 \begin{tabular}{||c||c|c|c|c|c||}\hline \hline 
  & & & & &  \\
S. No. & $C_{54}(R)$ &  $\theta_{54,3,2}(C_{54}(R))$ = &  $\theta_{54,3,4}(C_{54}(R))$ =  & $C_{54}(S)\in T1_{54}(C_{54}(R))$ & T1 or \\

& & & & & T2  \\

 & & &   & &  \\ \hline \hline

139.  & $C_{54}(2,3,6,12,16,20,24)$ &  $C_{54}(3,4,6,12,14,22,24)$ &  $C_{54}(3,6,8,10,12,24,26)$ & $C_{54}(2,3,6,12,16,20,24)$,  & \\   

& & & & $C_{54}(6,8,10,12,15,24,26)$, &   \\ 

     & & & & $C_{54}(4,6,12,14,21,22,24)$.  &  T2 \\

140.  & $C_{54}(2,3,6,12,16,20,27)$ &  $C_{54}(3,4,6,12,14,22,27)$ &  $C_{54}(3,6,8,10,12,26,27)$ & $C_{54}(2,3,6,12,16,20,27)$,  & \\
\hline \hline
\end{tabular}}
\end{center}
\end{table} 

\begin{table}
	\caption{ {\footnotesize Finding $\theta_{54,3,2}(C_{54}(R))$, $\theta_{54,3,4}(C_{54}(R))$ $\&$ $T1_{54}(C_{54}(R))$ as $\theta_{54,3,6}(C_{54}(R))$ = $C_{54}(R)$.}}
\begin{center}
\scalebox{.78}{
 \begin{tabular}{||c||c|c|c|c|c||}\hline \hline 
  & & & & &  \\
S. No. & $C_{54}(R)$ &  $\theta_{54,3,2}(C_{54}(R))$ = &  $\theta_{54,3,4}(C_{54}(R))$ =  & $C_{54}(S)\in T1_{54}(C_{54}(R))$ & T1   or \\

& & & & & T2  \\

 & & &   & &  \\ \hline \hline

159.  & $C_{54}(2,3,9,15,16,20,27)$ &  $C_{54}(3,4,9,14,15,22,27)$ &  $C_{54}(3,8,9,10,15,26,27)$ & $C_{54}(2,3,9,15,16,20,27)$, & \\

& & & & $C_{54}(8,9,10,15,21,26,27)$, &   \\ 

     & & & & $C_{54}(3,4,9,14,21,22,27)$.  &  T2 \\

160.  & $C_{54}(2,3,9,16,18,20,21)$ &  $C_{54}(3,4,9,14,18,21,22)$ &  $C_{54}(3,8,9,10,18,21,26)$ & $C_{54}(2,3,9,16,18,20,21)$, & \\
\hline \hline
\end{tabular}}
\end{center}
\end{table} 

\begin{table}
	\caption{ {\footnotesize Finding $\theta_{54,3,2}(C_{54}(R))$, $\theta_{54,3,4}(C_{54}(R))$ and $T1_{54}(C_{54}(R))$ as $\theta_{54,3,6}(C_{54}(R))$ = $C_{54}(R)$.}}
\begin{center}
\scalebox{.77}{
 \begin{tabular}{||c||c|c|c|c|c||}\hline \hline 
  & & & & &  \\
S. No. & $C_{54}(R)$ &  $\theta_{54,3,2}(C_{54}(R))$ = &  $\theta_{54,3,4}(C_{54}(R))$ =  & $C_{54}(S)\in T1_{54}(C_{54}(R))$ & T1   or \\

& & & & & T2  \\

 & & &   & &  \\ \hline \hline

179.  & $C_{54}(2,3,15,16,20,21,24)$ &  $C_{54}(3,4,14,15,21,22,24)$ &  $C_{54}(3,8,10,15,21,24,26)$ & $C_{54}(2,3,15,16,20,21,24)$, & \\

& & & & $C_{54}(3,8,10,12,15,21,26)$, &   \\ 

     & & & & $C_{54}(3,4,6,14,15,21,22)$.  &  T2 \\

180.  & $C_{54}(2,3,15,16,20,21,27)$ &  $C_{54}(3,4,14,15,21,22,27)$ &  $C_{54}(3,8,10,15,21,26,27)$ & $C_{54}(2,3,15,16,20,21,27)$, & \\
\hline \hline
\end{tabular}}
\end{center}
\end{table} 

\begin{table}
	\caption{ {\footnotesize Finding $\theta_{54,3,2}(C_{54}(R))$, $\theta_{54,3,4}(C_{54}(R))$ $\&$ $T1_{54}(C_{54}(R))$ as $\theta_{54,3,6}(C_{54}(R))$ = $C_{54}(R)$.}}
\begin{center}
\scalebox{.75}{
 \begin{tabular}{||c||c|c|c|c|c||}\hline \hline 
  & & & & &  \\
S. No. & $C_{54}(R)$ &  $\theta_{54,3,2}(C_{54}(R))$ = &  $\theta_{54,3,4}(C_{54}(R))$ =  & $C_{54}(S)\in T1_{54}(C_{54}(R))$ & T1   or \\

& & & & & T2  \\

 & & &   & &  \\ \hline \hline

200.   & $C_{54}(2,6,9,16,20,24,27)$ &   $C_{54}(4,6,9,14,22,24,27)$ &   $C_{54}(6,8,9,10,24,26,27)$ & $C_{54}(2,6,9,16,20,24,27)$, & \\

& & & & $C_{54}(8,9,10,12,24,26,27)$, &  \\ 

     & & & & $C_{54}(4,6,9,12,14,22,27)$.  &  T2 \\

201.   & $C_{54}(2,6,12,15,16,18,20)$ &   $C_{54}(4,6,12,14,15,18,22)$ &   $C_{54}(6,8,10,12,15,18,26)$ & $C_{54}(2,6,12,15,16,18,20)$, & \\
\hline \hline
\end{tabular}}
\end{center}
\end{table} 

\begin{table}
	\caption{ {\footnotesize Finding $\theta_{54,3,2}(C_{54}(R))$, $\theta_{54,3,4}(C_{54}(R))$ $\&$ $T1_{54}(C_{54}(R))$ as $\theta_{54,3,6}(C_{54}(R))$ = $C_{54}(R)$.}}
\begin{center}
\scalebox{.75}{
 \begin{tabular}{||c||c|c|c|c|c||}\hline \hline 
  & & & & &  \\
S. No. & $C_{54}(R)$ &  $\theta_{54,3,2}(C_{54}(R))$ = &  $\theta_{54,3,4}(C_{54}(R))$ =  & $C_{54}(S)\in T1_{54}(C_{54}(R))$ & T1   or \\

& & & & & T2  \\

 & & &   & &  \\ \hline \hline

222.   & $C_{54}(2,9,12,15,16,20,21)$ &   $C_{54}(4,9,12,14,15,21,22)$ &   $C_{54}(8,9,10,12,15,21,26)$ & $C_{54}(2,9,12,15,16,20,21)$, & \\

& & & & $C_{54}(3,6,8,9,10,21,26)$, &  \\ 

     & & & & $C_{54}(3,4,9,14,15,22,24)$.  &  T2 \\

223.   & $C_{54}(2,9,12,15,16,20,24)$ &   $C_{54}(4,9,12,14,15,22,24)$ &   $C_{54}(8,9,10,12,15,24,26)$ & $C_{54}(2,9,12,15,16,20,24)$, & \\
\hline \hline
\end{tabular}}
\end{center}
\end{table} 

\begin{table}
	\caption{ {\footnotesize Finding $\theta_{54,3,2}(C_{54}(R))$, $\theta_{54,3,4}(C_{54}(R))$ $\&$ $T1_{54}(C_{54}(R))$ as $\theta_{54,3,6}(C_{54}(R))$ = $C_{54}(R)$.}}
\begin{center}
\scalebox{.75}{
 \begin{tabular}{||c||c|c|c|c|c||}\hline \hline 
  & & & & &  \\
S. No. & $C_{54}(R)$ &  $\theta_{54,3,2}(C_{54}(R))$ = &  $\theta_{54,3,4}(C_{54}(R))$ =  & $C_{54}(S)\in T1_{54}(C_{54}(R))$ & T1   or \\

& & & & & T2  \\

 & & &   & &  \\ \hline \hline

244.   & $C_{54}(2,12,15,16,20,21,24)$ &   $C_{54}(4,12,14,15,21,22,24)$ &   $C_{54}(8,10,12,15,21,24,26)$ & $C_{54}(2,12,15,16,20,21,24)$, & \\

& & & & $C_{54}(3,6,8,10,12,21,26)$, &  \\ 

     & & & & $C_{54}(3,4,6,14,15,22,24)$.  &  T2 \\

245.   & $C_{54}(2,12,15,16,20,21,27)$ &   $C_{54}(4,12,14,15,21,22,27)$ &   $C_{54}(8,10,12,15,21,26,27)$ & $C_{54}(2,12,15,16,20,21,27)$, & \\
\hline \hline
\end{tabular}}
\end{center}
\end{table} 

\begin{table}
	\caption{ {\footnotesize Finding $\theta_{54,3,2}(C_{54}(R))$, $\theta_{54,3,4}(C_{54}(R))$ $\&$ $T1_{54}(C_{54}(R))$ as $\theta_{54,3,6}(C_{54}(R))$ = $C_{54}(R)$.}}
\begin{center}
\scalebox{.75}{
 \begin{tabular}{||c||c|c|c|c|c||}\hline \hline 
  & & & & &  \\
S. No. & $C_{54}(R)$ &  $\theta_{54,3,2}(C_{54}(R))$ = &  $\theta_{54,3,4}(C_{54}(R))$ =  & $C_{54}(S)\in T1_{54}(C_{54}(R))$ & T1   or \\

& & & & & T2  \\

 & & &   & &  \\ \hline \hline

266.  & $C_{54}(2,3,6,9,16,18,20,24)$ &  $C_{54}(3,4,6,9,14,18,22,24)$ &  $C_{54}(3,6,8,9,10,18,24,26)$ & $C_{54}(2,3,6,9,16,18,20,24)$,  & \\   

& & & & $C_{54}(8,9,10,12,15,18,24,26)$, &   \\ 

     & & & & $C_{54}(4,6,9,12,14,18,21,22)$.  &  T2 \\

267.  & $C_{54}(2,3,6,9,16,18,20,27)$ &  $C_{54}(3,4,6,9,14,18,22,27)$ &  $C_{54}(3,6,8,9,10,18,26,27)$ & $C_{54}(2,3,6,9,16,18,20,27)$,  & \\
\hline \hline
\end{tabular}}
\end{center}
\end{table} 

\begin{table}
	\caption{ {\footnotesize Finding $\theta_{54,3,2}(C_{54}(R))$, $\theta_{54,3,4}(C_{54}(R))$ $\&$ $T1_{54}(C_{54}(R))$ as $\theta_{54,3,6}(C_{54}(R))$ = $C_{54}(R)$.}}
\begin{center}
\scalebox{.75}{
 \begin{tabular}{||c||c|c|c|c|c||}\hline \hline 
  & & & & &  \\
S. No. & $C_{54}(R)$ &  $\theta_{54,3,2}(C_{54}(R))$ = &  $\theta_{54,3,4}(C_{54}(R))$ =  & $C_{54}(S)\in T1_{54}(C_{54}(R))$ & T1   or \\

& & & & & T2  \\

 & & &   & &  \\ \hline \hline

288.  & $C_{54}(2,3,6,16,18,20,21,27)$ &  $C_{54}(3,4,6,14,18,21,22,27)$ &  $C_{54}(3,6,8,10,18,21,26,27)$ & $C_{54}(2,3,6,16,18,20,21,27)$,  & \\   

& & & & $C_{54}(3,8,10,15,18,24,26,27)$, &   \\ 

     & & & & $C_{54}(4,12,14,15,18,21,22,27)$.  &  T2 \\

289. & $C_{54}(2,3,6,16,18,20,24,27)$ &  $C_{54}(3,4,6,14,18,22,24,27)$ &  $C_{54}(3,6,8,10,18,24,26,27)$ & $C_{54}(2,3,6,16,18,20,24,27)$,  & \\
\hline \hline
\end{tabular}}
\end{center}
\end{table} 

\begin{table}
	\caption{ {\footnotesize Finding $\theta_{54,3,2}(C_{54}(R))$, $\theta_{54,3,4}(C_{54}(R))$ $\&$ $T1_{54}(C_{54}(R))$ as $\theta_{54,3,6}(C_{54}(R))$ = $C_{54}(R)$.}}
\begin{center}
\scalebox{.73}{
 \begin{tabular}{||c||c|c|c|c|c||}\hline \hline 
  & & & & &  \\
S. No. & $C_{54}(R)$ &  $\theta_{54,3,2}(C_{54}(R))$ = &  $\theta_{54,3,4}(C_{54}(R))$ =  & $C_{54}(S)\in T1_{54}(C_{54}(R))$ & T1   or \\

& & & & & T2  \\

 & & &   & &  \\ \hline \hline

310.  & $C_{54}(2,3,9,16,20,21,24,27)$ &  $C_{54}(3,4,9,14,21,22,24,27)$ &  $C_{54}(3,8,9,10,21,24,26,27)$ & $C_{54}(2,3,9,16,20,21,24,27)$, & \\

& & & & $C_{54}(3,8,9,10,12,15,26,27)$, &   \\ 

     & & & & $C_{54}(4,6,9,14,15,21,22,27)$.  &  T2 \\

311.  & $C_{54}(2,3,12,15,16,18,20,21)$ &  $C_{54}(3,4,12,14,15,18,21,22)$ &  $C_{54}(3,8,10,12,15,18,21,26)$ & $C_{54}(2,3,12,15,16,18,20,21)$, & \\
\hline \hline
\end{tabular}}
\end{center}
\end{table} 

\begin{table}
	\caption{ {\footnotesize Finding $\theta_{54,3,2}(C_{54}(R))$, $\theta_{54,3,4}(C_{54}(R))$ $\&$ $T1_{54}(C_{54}(R))$ as $\theta_{54,3,6}(C_{54}(R))$ = $C_{54}(R)$.}}
\begin{center}
\scalebox{.7}{
 \begin{tabular}{||c||c|c|c|c|c||}\hline \hline 
  & & & & &  \\
S. No. & $C_{54}(R)$ &  $\theta_{54,3,2}(C_{54}(R))$ = &  $\theta_{54,3,4}(C_{54}(R))$ =  & $C_{54}(S)\in T1_{54}(C_{54}(R))$ & T1   or \\

& & & & & T2  \\

 & & &   & &  \\ \hline \hline

332.  & $C_{54}(2,6,9,12,16,18,20,27)$ &  $C_{54}(4,6,9,12,14,18,22,27)$ &  $C_{54}(6,8,9,10,12,18,26,27)$ & $C_{54}(2,6,9,12,16,18,20,27)$,  & \\   

& & & & $C_{54}(6,8,9,10,18,24,26,27)$, &   \\ 

     & & & & $C_{54}(4,9,12,14,18,22,24,27)$.  &  T2 \\

333.  & $C_{54}(2,6,9,12,16,20,21,24)$ &  $C_{54}(4,6,9,12,14,21,22,24)$ &  $C_{54}(6,8,9,10,12,21,24,26)$ & $C_{54}(2,6,9,12,16,20,21,24)$,  & \\
\hline \hline
\end{tabular}}
\end{center}
\end{table} 

\begin{table}
	\caption{ {\footnotesize Finding $\theta_{54,3,2}(C_{54}(R))$, $\theta_{54,3,4}(C_{54}(R))$ $\&$ $T1_{54}(C_{54}(R))$ as $\theta_{54,3,6}(C_{54}(R))$ = $C_{54}(R)$.}}
\begin{center}
\scalebox{.7}{
 \begin{tabular}{||c||c|c|c|c|c||}\hline \hline 
  & & & & &  \\
S. No. & $C_{54}(R)$ &  $\theta_{54,3,2}(C_{54}(R))$ = &  $\theta_{54,3,4}(C_{54}(R))$ =  & $C_{54}(S)\in T1_{54}(C_{54}(R))$ & T1   or \\

& & & & & T2  \\

 & & &   & &  \\ \hline \hline

355.  & $C_{54}(2,6,12,16,20,21,24,27)$ &  $C_{54}(4,6,12,14,21,22,24,27)$ &  $C_{54}(6,8,10,12,21,24,26,27)$ & $C_{54}(2,6,12,16,20,21,24,27)$,  & \\   

& & & & $C_{54}(3,6,8,10,12,24,26,27)$, &   \\ 

     & & & & $C_{54}(4,6,12,14,15,22,24,27)$.  &  T2 \\

356.  & $C_{54}(2,6,15,16,18,20,21,24)$ &  $C_{54}(4,6,14,15,18,21,22,24)$ &  $C_{54}(6,8,10,15,18,21,24,26)$ & $C_{54}(2,6,15,16,18,20,21,24)$,  & \\
\hline \hline
\end{tabular}}
\end{center}
\end{table} 

\begin{table}
	\caption{ {\footnotesize Finding $\theta_{54,3,2}(C_{54}(R))$, $\theta_{54,3,4}(C_{54}(R))$ $\&$ $T1_{54}(C_{54}(R))$ as $\theta_{54,3,6}(C_{54}(R))$ = $C_{54}(R)$.}}
\begin{center}
\scalebox{.68}{
 \begin{tabular}{||c||c|c|c|c|c||}\hline \hline 
  & & & & &  \\
S. No. & $C_{54}(R)$ &  $\theta_{54,3,2}(C_{54}(R))$ = &  $\theta_{54,3,4}(C_{54}(R))$ =  & $C_{54}(S)\in T1_{54}(C_{54}(R))$ & T1   or \\

& & & & & T2  \\

 & & &   & &  \\ \hline \hline

378.  & $C_{54}(2,12,15,16,18,20,24,27)$ &  $C_{54}(4,12,14,15,18,22,24,27)$ &  $C_{54}(8,10,12,15,18,24,26,27)$ & $C_{54}(2,12,15,16,18,20,24,27)$, & \\

& & & & $C_{54}(6,8,10,12,18,21,26,27)$, &   \\ 

     & & & & $C_{54}(3,4,6,14,18,22,24,27)$.  &  T2 \\ 

379.  & $C_{54}(2,12,15,16,20,21,24,27)$ &  $C_{54}(4,12,14,15,21,22,24,27)$ &  $C_{54}(8,10,12,15,21,24,26,27)$ & $C_{54}(2,12,15,16,20,21,24,27)$, & \\
\hline \hline
\end{tabular}}
\end{center}
\end{table} 

\begin{table}
	\caption{ {\footnotesize Finding $\theta_{54,3,2}(C_{54}(R))$, $\theta_{54,3,4}(C_{54}(R))$ $\&$ $T1_{54}(C_{54}(R))$ as $\theta_{54,3,6}(C_{54}(R))$ = $C_{54}(R)$.}}
\begin{center}
\scalebox{.65}{
 \begin{tabular}{||c||c|c|c|c|c||}\hline \hline 
  & & & & &  \\
S. No. & $C_{54}(R)$ &  $\theta_{54,3,2}(C_{54}(R))$ = &  $\theta_{54,3,4}(C_{54}(R))$ =  & $C_{54}(S)\in T1_{54}(C_{54}(R))$ & T1   or \\

& & & & & T2  \\

 & & &   & &  \\ \hline \hline

402.  & $C_{54}(2,3,6,12,15,16,18,20,21)$ &  $C_{54}(3,4,6,12,14,15,18,21,22)$ &  $C_{54}(3,6,8,10,12,15,18,21,26)$ & $C_{54}(2,3,6,12,15,16,18,20,21)$,  & \\   

& & & & $C_{54}(3,6,8,10,15,18,21,24,26)$, &   \\ 

     & & & & $C_{54}(3,4,12,14,15,18,21,22,24)$.  &  T2 \\

403.  & $C_{54}(2,3,6,12,15,16,18,20,24)$ &  $C_{54}(3,4,6,12,14,15,18,22,24)$ &  $C_{54}(3,6,8,10,12,15,18,24,26)$ & $C_{54}(2,3,6,12,15,16,18,20,24)$,  & \\
\hline \hline
\end{tabular}}
\end{center}
\end{table} 

\begin{table}
	\caption{ {\footnotesize Finding $\theta_{54,3,2}(C_{54}(R))$, $\theta_{54,3,4}(C_{54}(R))$ $\&$ $T1_{54}(C_{54}(R))$ as $\theta_{54,3,6}(C_{54}(R))$ = $C_{54}(R)$.}}
\begin{center}
\scalebox{.68}{
 \begin{tabular}{||c||c|c|c|c|c||}\hline \hline 
  & & & & &  \\
S. No. & $C_{54}(R)$ &  $\theta_{54,3,2}(C_{54}(R))$ = &  $\theta_{54,3,4}(C_{54}(R))$ =  & $C_{54}(S)\in T1_{54}(C_{54}(R))$ & T1   or \\

& & & & & T2  \\

 & & &   & &  \\ \hline \hline

427.  & $C_{54}(2,3,9,15,16,18,20,21,24)$ &  $C_{54}(3,4,9,14,15,18,21,22,24)$ &  $C_{54}(3,8,9,10,15,18,21,24,26)$ & $C_{54}(2,3,9,15,16,18,20,21,24)$, & \\

& & & & $C_{54}(3,8,9,10,12,15,18,21,26)$, &   \\ 

     & & & & $C_{54}(3,4,6,9,14,15,18,21,22)$.  &  T2 \\

428.  & $C_{54}(2,3,9,15,16,18,20,21,27)$ &  $C_{54}(3,4,9,14,15,18,21,22,27)$ &  $C_{54}(3,8,9,10,15,18,21,26,27)$ & $C_{54}(2,3,9,15,16,18,20,21,27)$, & \\
\hline \hline
\end{tabular}}
\end{center}
\end{table} 

\begin{table}
	\caption{ {\footnotesize Finding $\theta_{54,3,2}(C_{54}(R))$, $\theta_{54,3,4}(C_{54}(R))$ $\&$ $T1_{54}(C_{54}(R))$ as $\theta_{54,3,6}(C_{54}(R))$ = $C_{54}(R)$.}}
\begin{center}
\scalebox{.68}{
 \begin{tabular}{||c||c|c|c|c|c||}\hline \hline 
  & & & & &  \\
S. No. & $C_{54}(R)$ &  $\theta_{54,3,2}(C_{54}(R))$ = &  $\theta_{54,3,4}(C_{54}(R))$ =  & $C_{54}(S)\in T1_{54}(C_{54}(R))$ & T1   or \\

& & & & & T2  \\

 & & &   & &  \\ \hline \hline

451.  & $C_{54}(2,6,9,15,16,20,21,24,27)$ &  $C_{54}(4,6,9,14,15,21,22,24,27)$ &  $C_{54}(6,8,9,10,15,21,24,26,27)$ & $C_{54}(2,6,9,15,16,20,21,24,27)$,  & \\   

& & & & $C_{54}(3,8,9,10,12,21,24,26,27)$, &   \\ 

     & & & & $C_{54}(3,4,6,9,12,14,15,22,27)$.  &  T2 \\

452.  & $C_{54}(2,6,9,16,18,20,21,24,27)$ &  $C_{54}(4,6,9,14,18,21,22,24,27)$ &  $C_{54}(6,8,9,10,18,21,24,26,27)$ & $C_{54}(2,6,9,16,18,20,21,24,27)$,  & \\
\hline \hline
\end{tabular}}
\end{center}
\end{table} 

\begin{table}
	\caption{ {\footnotesize Finding $\theta_{54,3,2}(C_{54}(R))$, $\theta_{54,3,4}(C_{54}(R))$ $\&$ $T1_{54}(C_{54}(R))$ as $\theta_{54,3,6}(C_{54}(R))$ = $C_{54}(R)$.}}
\begin{center}
\scalebox{.59}{
 \begin{tabular}{||c||c|c|c|c|c||}\hline \hline 
  & & & & &  \\
S. No. & $C_{54}(R)$ &  $\theta_{54,3,2}(C_{54}(R))$ = &  $\theta_{54,3,4}(C_{54}(R))$ =  & $C_{54}(S)\in T1_{54}(C_{54}(R))$ & T1   or \\

& & & & & T2  \\

 & & &   & &  \\ \hline \hline

475.  & $C_{54}(2,3,6,9,12,16,20,21,24,27)$ &  $C_{54}(3,4,6,9,12,14,21,22,24,27)$ &  $C_{54}(3,6,8,9,10,12,21,24,26,27)$ & $C_{54}(2,3,6,9,12,16,20,21,24,27)$,  & \\   

& & & & $C_{54}(3,6,8,9,10,12,15,24,26,27)$, &   \\ 

     & & & & $C_{54}(4,6,9,12,14,15,21,22,24,27)$.  &  T2 \\

476.  & $C_{54}(2,3,6,9,15,16,18,20,21,24)$ &  $C_{54}(3,4,6,9,14,15,18,21,22,24)$ &  $C_{54}(3,6,8,9,10,15,18,21,24,26)$ & $C_{54}(2,3,6,9,15,16,18,20,21,24)$,  & \\
\hline \hline
\end{tabular}}
\end{center}
\end{table} 

\begin{table}
	\caption{ {\footnotesize Finding $\theta_{54,3,2}(C_{54}(R))$, $\theta_{54,3,4}(C_{54}(R))$ $\&$ $T1_{54}(C_{54}(R))$ as $\theta_{54,3,6}(C_{54}(R))$ = $C_{54}(R)$.}}
\begin{center}
\scalebox{.59}{
 \begin{tabular}{||c||c|c|c|c|c||}\hline \hline 
  & & & & &  \\
S. No. & $C_{54}(R)$ &  $\theta_{54,3,2}(C_{54}(R))$ = &  $\theta_{54,3,4}(C_{54}(R))$ =  & $C_{54}(S)\in T1_{54}(C_{54}(R))$ & T1   or \\

& & & & & T2  \\

 & & &   & &  \\ \hline \hline

502.  & $C_{54}(2,3,6,9,12,15,16,18,20,21,24)$ &  $C_{54}(3,4,6,9,12,14,15,18,21,22,24)$ &  $C_{54}(3,6,8,9,10,12,15,18,21,24,26)$ & $C_{54}(2,3,6,9,12,15,16,18,20,21,24)$,  & \\   

& & & & $C_{54}(3,6,8,9,10,12,15,18,21,24,26)$, &   \\ 

     & & & & $C_{54}(3,4,6,9,12,14,15,18,21,22,24)$.  &  T1 \\

503.  & $C_{54}(2,3,6,9,12,15,16,18,20,21,27)$ &  $C_{54}(3,4,6,9,12,14,15,18,21,22,27)$ &  $C_{54}(3,6,8,9,10,12,15,18,21,26,27)$ & $C_{54}(2,3,6,9,12,15,16,18,20,21,27)$,  & \\
\hline \hline
\end{tabular}}
\end{center}
\end{table} 

\section{Conclusion}

 In \cite{v2-2} - \cite{v2-4}, the author established that the number of pairs of Type-2 isomorphic circulant graphs of orders 16, 24 and 32 are 8, 32 and 384; and the number of triples of Type-2 isomorphic circulant graphs of order 27 is 12 and presented all. The author feels that a lot of scope is there for further research and proposes the following open problems on this topic.

Here, we propose the following open problems on circulant graphs of order 48.

\begin{oprm} \label{op1} {\rm The following pairs of circulant graphs are non-isomorphic for $s$ = 3, 9, 15, 21.
\begin{enumerate}	\item [\rm (a)]  $C_{48}(1,s,23)$ and $C_{48}(s,11,13)$;  and 
			\item [\rm (b)]  $C_{48}(5,s,19)$ and $C_{48}(s,7,17)$. \hfill $\Box$ 
\end{enumerate} }	
\end{oprm}

\begin{oprm} \label{op2} {\rm Find all pairs of isomorphic circulant graphs of order 48. \hfill $\Box$}	
\end{oprm}

\begin{oprm} \label{op3} {\rm Find all pairs of Type-1 isomorphic circulant graphs of order 48. \hfill $\Box$}	
\end{oprm}

\begin{oprm} \label{op4} {\rm Find all pairs of Type-2 isomorphic circulant graphs of order 48. \hfill $\Box$}	
\end{oprm}
 
\begin{oprm} \label{op5} {\rm For $s$ = 2, 4, 8, 10, 14, 16, 20, 22, 26, each triples of circulant graphs  
 $C_{54}(1,s,17,19)$, $C_{54}(5,s,13,23)$, $C_{54}(s,7,11,25)$ are non-isomorphic.  \hfill $\Box$ }	
\end{oprm}

\begin{oprm} \label{op6} {\rm Find all pairs of isomorphic circulant graphs of order 54. \hfill $\Box$}	
\end{oprm}

\begin{oprm} \label{op7} {\rm Find all pairs of Type-1 isomorphic circulant graphs of order 54. \hfill $\Box$}	
\end{oprm}

\begin{oprm} \label{op8} {\rm Find all pairs of Type-2 isomorphic circulant graphs of order 54. \hfill $\Box$}	
\end{oprm}

\vspace{.1cm}
\noindent
\textbf{Declaration of competing interest}\quad The author declares that he has no conflict of interest.



\begin {thebibliography}{10}

\bibitem {ad67}  
A. Adam, 
{\it Research problem 2-10},  
J. Combinatorial Theory, {\bf 3} (1967), 393.

\bibitem {v2-2-arX} 
V. Vilfred Kamalappan, 
\emph{All Type-2 Isomorphic Circulant Graphs $C_{16}(R)$ and $C_{24}(S)$}, 
arXiv: 2508.09384v1 [math.CO] 12 Aug 2025, 28 pages.

\bibitem {v24} 
V. Vilfred Kamalappan, 
\emph{A study on Type-2 Isomorphic Circulant Graphs and related Abelian Groups}, 
arXiv: 2012.11372v11 [math.CO] (26 Nov. 2024), 183 pages.

\bibitem {v2-1} 
V. Vilfred Kamalappan, 
\emph{A study on Type-2 Isomorphic Circulant Graphs. \\ Part 1: Type-2 isomorphic circulant graphs $C_n(R)$ w.r.t. $m$ = 2}. 
Preprint. 32 pages

\bibitem {v2-2} 
V. Vilfred Kamalappan, 
\emph{A study on Type-2 isomorphic circulant graphs. \\ Part 2: Type-2 isomorphic circulant graphs of orders 16, 24, 27}. 
Preprint. 32 pages

\bibitem {v2-3} 
V. Vilfred Kamalappan, 
\emph{A study on Type-2 isomorphic circulant graphs. \\ Part 3: 384 pairs of Type-2 isomorphic circulant graphs $C_{32}(R)$}. 
Preprint. 42 pages

\bibitem {v2-4} 
V. Vilfred Kamalappan, 
\emph{A study on Type-2 isomorphic circulant graphs. \\ Part 4: 960 triples of Type-2 isomorphic circulant graphs $C_{54}(R)$}. 
Preprint. 76 pages

\bibitem {v2-5} 
V. Vilfred Kamalappan, 
\emph{A study on Type-2 isomorphic circulant graphs. \\ Part 5: Type-2 isomorphic circulant graphs of orders 48, 81, 96}. 
Preprint. 33 pages

\bibitem {v2-6} 
V. Vilfred Kamalappan, 
\emph{A study on Type-2 Isomorphic Circulant Graphs. \\ Part 6: Abelian groups $(T2_{n, m}(C_n(R)), \circ)$ and $(V_{n, m}(C_n(R)), \circ)$}. 
Preprint. 19 pages

\bibitem {v2-7} 
V. Vilfred Kamalappan, 
\emph{A Study on Type-2 Isomorphic Circulant Graphs: Part 7: $C_{432}(R)$, $C_{6750}(S)$ - each has 2 types of Type-2 isomorphic circulant graphs}. 
Preprint. 351 pages

\bibitem {v2-8} 
V. Vilfred Kamalappan, 
\emph{A study on Type-2 Isomorphic Circulant Graphs. \\ Part 8: Isomorphism series, digraph and graph of $C_n(R)$}. 
Preprint. 79 pages

\bibitem {v2-9} 
V. Vilfred Kamalappan and P. Wilson, 
\emph{A study on Type-2 Isomorphic Circulant Graphs. \\ Part 9: Computer program to show Type-1 and -2 isomorphic circulant graphs}. 
Preprint. 21 pages

\bibitem {v2-10} 
V. Vilfred Kamalappan and P. Wilson, 
\emph{A study on Type-2 Isomorphic Circulant Graphs. \\ Part 10: Type-2 isomorphic $C_{np^3}(R)$ w.r.t. $m$ = $p$ and related groups}. 
Preprint. 20 pages

\end{thebibliography}
\end{document}